\documentclass[11pt]{article}

\usepackage{fullwidth}
\usepackage{xy}
\usepackage[top=1.5in, bottom=1.5in, left=1in, right=1in]{geometry}
\usepackage{amsgen}
\usepackage{amsmath}
\usepackage{amstext}
\usepackage{amsbsy}
\usepackage{amsopn}
\usepackage{amsfonts}
\usepackage{amssymb}
\usepackage{comment}
\usepackage{eepic}
\usepackage{pgf,tikz}
\usepackage{graphicx}
\usepackage{epsf}
\usepackage{pstricks}
\usepackage{makeidx}
\xyoption{all}

\def\Box{\square}
\def\edge{\relbar\joinrel\relbar}

\def\mapright#1{\smash{\mathop{\longrightarrow}\limits^{#1}}}
\def\tra#1{\smash{\mathop{\mid\kern
-1pt\joinrel\relbar\joinrel\relbar}\limits^{*}_{#1}}}
\def\longtra#1{\smash{\mathop{\mid\kern
-1pt\joinrel\relbar\joinrel\relbar\joinrel\relbar}\limits^{*}_{#1}}}
\def\vlongtra#1{\smash{\mathop{\mid\kern
-1pt\joinrel\relbar\joinrel\relbar\joinrel\relbar\joinrel\relbar}\limits^{*}_{#1}}}
\def\vvlongtra#1{\smash{\mathop{\mid\kern
-1pt\joinrel\relbar\joinrel\relbar\joinrel\relbar\joinrel\relbar\joinrel\relbar}\limits^{*}_{#1}}}
\def\vvvlongtra#1{\smash{\mathop{\mid\kern
-1pt\joinrel\relbar\joinrel\relbar\joinrel\relbar\joinrel\relbar\joinrel\relbar\joinrel\relbar}\limits^{*}_{#1}}}
\def\etra#1{\smash{\mathop{\mid\kern
-1pt\joinrel\relbar\joinrel\relbar}\limits_{#1}}}

\def\vlongrightarrow{\relbar\joinrel\longrightarrow}

\def\vvlongrightarrow{\relbar\joinrel\vlongrightarrow}

\def\longmapright#1{\smash{\mathop{\vlongrightarrow}\limits^{#1}}}
\def\vlongmapright#1{\smash{\mathop{\vvlongrightarrow}\limits^{#1}}}

\def\X{{\cal{X}}}

\def\Rw{\Rightarrow}
\def\oo{\overline}
\def\un{\underline}

\def\wh{\widehat}
\def\B{{\cal{B}}}

\def\F{{\cal{F}}}
\def\Q{{\cal{Q}}}

\def\L{{\cal{L}}}
\def\M{{\cal{M}}}
\def\N{\mathbb{N}}

\def\S{{\cal{S}}}

\def\cl{\mbox{Cl}}
\def\defe{\mbox{def}}
\def\Defe{\mbox{Def}}

\def\clos{\mbox{Cl}}

\def\SPC{\mbox{SPC}}

\def\fct{\mbox{Fct}}

\def\dim{\mbox{dim}}
\def\codim{\mbox{codim}}

\def\pav{\mbox{Pav}}
\def\bpav{\mbox{BPav}}
\def\tbpav{\mbox{TBPav}}

\def\ker{\mbox{Ker}\,}

\def\min{\mbox{min}}

\def\pure{\mbox{pure}}

\def\gu{\mbox{GU}}
\def\ngu{\mbox{NGU}}
\def\mngu{\mbox{MNGU}}
\def\mgu{\mbox{mGU}}
\def\flatx{\mbox{Fl}}
\def\flats{\mbox{Fl\,}}
\def\P{{\cal{P}}}

\def\G{{\cal{G}}}
\def\H{{\cal{H}}}

\def\tbr{{\cal{TBR}}}
\def\tbp{{\cal{TBP}}}
\def\D{{\cal{D}}}
\def\V{{\cal{V}}}

\def\Ncal{{\cal{N}}}
\def\Z{\mathbb{Z}}

\def\p{\varphi}

\def\inv{^{-1}}

\def\bi{\begin{itemize}}
\def\ei{\end{itemize}}
\def\beq{\begin{equation}}
\def\eeq{\end{equation}}
\def\ba{\begin{array}}
\def\ea{\end{array}}

\def\J{{\cal{J}}}
\def\K{{\cal{K}}}
\def\Y{{\cal{Y}}}

\def\up{^{{\rm up}}}
\def\upp#1{^{{\rm up}^{(#1)}}}

\newtheorem{T}{Theorem}[section]
\newcommand{\bt}{\begin{T}}
\newcommand{\et}{\end{T}}
\newcommand{\ftd}{$\square$\end{T}}

\newtheorem{Proposition}[T]{Proposition}
\newcommand{\bp}{\begin{Proposition}}
\newcommand{\ep}{\end{Proposition}}
\newcommand{\fpd}{$\square$\end{Proposition}}

\newtheorem{Lemma}[T]{Lemma}
\newcommand{\bl}{\begin{Lemma}}
\newcommand{\el}{\end{Lemma}}
\newcommand{\fld}{$\square$\end{Lemma}}

\newtheorem{Corol}[T]{Corollary}
\newcommand{\bc}{\begin{Corol}}
\newcommand{\ec}{\end{Corol}}
\newcommand{\fcd}{$\square$\end{Corol}}

\newtheorem{Conj}[T]{Conjecture}
\newcommand{\bj}{\begin{Conj}}
\newcommand{\ej}{\end{Conj}}
\newcommand{\fjd}{$\square$\end{Conj}}

\newtheorem{Result}[T]{Result}
\newcommand{\br}{\begin{Result}}
\newcommand{\er}{\end{Result}}
\newcommand{\frd}{$\square$\end{Result}}

\newtheorem{Remark}[T]{Remark}
\newcommand{\brem}{\begin{Remark}}
\newcommand{\erem}{\end{Remark}}
\newcommand{\fremd}{$\square$\end{Remark}}

\newtheorem{Example}[T]{Example}
\newcommand{\be}{\begin{Example}}
\newcommand{\ee}{\end{Example}}

\newtheorem{Problem}[T]{Problem}
\newcommand{\bq}{\begin{Problem}}
\newcommand{\eq}{\end{Problem}}

\newcommand{\proof}
   {\par\medbreak\noindent{\bf Proof}.\enspace}

\newcommand{\qed}{
$\Box$
\par\bigbreak}

\def\abstract#1{\par\bigskip
\begingroup\small
\baselineskip=12truept
\begin{center}ABSTRACT\end{center}
\par\medskip\par\noindent
\null\hfill\hbox{\vbox{\hsize=5truein\noindent#1}}
\hfill\null\par\endgroup\par}

\title{Topics in Boolean Representable Simplicial Complexes}
\author{{\bf Stuart Margolis, John Rhodes and Pedro V. Silva}}

\begin{document}
\maketitle

\begin{center}\small


\end{center}

\abstract{We study a number of topics in the theory of Boolean Representable Simplicial Complexes (BRSC). These include various operators on 
BRSC. We look at shellability in higher dimensions and propose a number of new conjectures.}

\newpage

\tableofcontents

\newpage

\section{Introduction}

\section{Preliminaries}

A {\em simplicial complex} is a structure of the form $\H\; = (V,H)$, where $V$ is a finite nonempty set and $H \subseteq 2^V$ satisfies:
\bi
\item
$P_1(V) \subseteq H$;
\item
if $X \subseteq I \in H$, then $X \in H$.
\ei 

\section{The up operator}
\label{going}

In this section we explore the {\em up operator}, which associates a BRSC to each simplicial complex.

Let $\H\; = (V,H)$ be a simplicial complex. We define
$$H\up = H \cup \{ I \cup \{ p \} \mid I \in H,\; p \in V \setminus I \} \quad \mbox{and}\quad \H\up = (V,H\up).$$
We show that
\beq
\label{hup}
H\up = 2^V \setminus \{ X \in 2^V \setminus \{ \emptyset \} \mid X\setminus \{ x \} \notin H \mbox{ for every } x \in X\}.
\eeq

Indeed, $X \in H\up$ if and only if
$$X = \emptyset\mbox{ or $X \setminus \{ x \} \in H$ for some }x \in X$$
if and only if $X$ belongs to the right hand side of (\ref{hup}). Thus (\ref{hup}) holds.

Assume that $\emptyset \neq R \subseteq 2^V$. Let $M(R)$ denote the $R \times V$ boolean matrix $(m_{rv})$ defined as follows: for all $r \in R$ and $v \in V$, let
$$m_{rv} = \left\{
\begin{array}{ll}
0&\mbox{ if }v \in r\\
1&\mbox{ otherwise }
\end{array}
\right.$$
If $R$ contains $V$ and is closed under intersection, it is said to be a {\em Moore family} (see \cite{CIRR}). We shall assume also that $\emptyset \in R$ (since $\emptyset$ is a flat in every simplicial complex). If $R$ is a Moore family, then $(R,\cap)$ constitutes a $\wedge$-semilattice and since $V \in R$ then $R$ has a lattice structure if we consider the determined join. We can view $V$ as a $\vee$-generating set for the lattice $R$ through the mapping
$\p:V \to R$ defined by
$$v\p = \cap\{ r \in R \mid v \in r\}.$$
Then $M(R)$ is the matrix determined by the lattice $(R,V)$ (see \cite[Section 3.3]{RSm}), and so the BRSC defined by $M(R)$ is the complex $\J(R) = (V,J(R))$, where $J(R)$ is the set of transversals of the successive differences for some chain of $R$ (\cite[Theorem 5.4.2]{RSm}). 
It follows from (\cite[Lemma 5.2.1]{RSm}) that
\beq
\label{flaf}
R \subseteq \flatx(V,J(R)).
\eeq

Now let $\H\; = (V,H)$ be a simplicial complex. Since $H$ is closed under taking subsets, it is closed under intersection and so $H \cup \{ V \}$ is a Moore family. We can prove the following result.

\bt
\label{upp}
Let $\H\; = (V,H)$ be a simplicial complex. Then $M(H \cup \{ V \})$ is a boolean matrix representation for $\H\up$.
\et

\proof 
Assume first that $p_1,\ldots, p_n$ is an enumeration of the elements of $I \in H$. Then $I$ is a transversal of the successive differences for the chain 
$$\emptyset \subset p_1 \subset p_1p_2 \subset \ldots \subset p_1\ldots p_n = I$$ 
of $H$, hence $I \in J(H \cup \{ V \})$. Furthermore, if $q \in V \setminus I$, then $I \cup \{ q \}$ is a transversal of the successive differences for the chain 
$$\emptyset \subset p_1 \subset p_1p_2 \subset \ldots \subset p_1\ldots p_n = I \subset V$$
of $H \cup \{ V\}$, hence $I \cup \{ q \} \in J(H \cup \{ V \})$ and so $H\up \subseteq J(H \cup \{ V \})$.

Conversely, let $X \notin H$ be a transversal of the successive differences for a chain 
$$I_0 \subset I_1 \subset \ldots \subset I_n$$ 
of $H \cup \{ V\}$. Since $X \subseteq I_n$, it follows that $I_n = V$ and $I_{n-1} \in H$. Denoting by $p$ the single element of $X \setminus I_{n-1}$, we get $X \setminus \{ p \} \in H$ and $X \in H\up$, hence $J(H \cup \{ V \}) \subseteq H\up$ and so $J(H \cup \{ V \}) = H\up$ as required.
\qed

\bc
\label{fup}
Let $\H\; = (V,H)$ be a simplicial complex. Then $H \cup \{ V \} \subseteq {\rm Fl}(\H\up)$.
\ec

\proof
By \cite[Lemma 5.2.1]{RSm}, 
we have $H \cup \{ V \} \subseteq \flatx(M(H \cup \{ V\})) \subseteq
{\rm Fl}(\H\up)$.
\qed

The following example shows that the inclusion in Corollary \ref{fup} may be strict.

\be
\label{cfup}
Let $V = 1234$ and $H = P_{\leq 1}(V) \cup \{ 12,34\}$. Then ${\rm Fl}(\H\up) = P_{\leq 2}(V) \cup \{ V \}$.
\ee

Indeed, $\H\up = (V,P_{\leq 3}(V))$ and this easily yields the claim.

\bp
\label{mup}
Let $\H\; = (V,H)$ be a matroid. Then $\H\up$ is a matroid.
\ep

\proof
We may assume that $V \notin H$, otherwise $\H\up = \H$.

Let $I',J' \in H\up$ with $|I'| = |J'|+1$. The exchange property holding trivially if $J' \in H$, we may assume that $J' = J \cup \{ p\} \notin H$ with $J \in H$. Write $I' = I \cup \{ q\}$ with $I \in H$ and $q \notin I$. Since $\H$ is a matroid, we have $J \cup \{ i \} \in H$ for some $i \in I\setminus J$. Since $J' \notin H$, we have $i \neq p$, hence $J' \cup \{ i \} \in H\up$ and $i \in I' \setminus J'$. Therefore the exchange property holds in $\H\up$ and so $\H\up$ is a matroid.
\qed

Given a simplicial complex $\H\; = (V,H)$ and $k \geq 0$, the $k$-{\em
  truncation} of $\H$ is the
simplicial complex $\H_k\; = (V,H_k)$
defined by $H_k = H \cap P_{\leq k}(V)$. 

\bp
\label{pup}
Let $\H\; = (V,H)$ be a simplicial complex. Then:
\bi
\item[(i)] If $V \notin H$, then ${\rm dim}(\H\up) = {\rm dim}\H +1$.
\item[(ii)] If $V \notin H$ and $\H \in {\rm Pav}(d)$, then $\H\up \in {\rm BPav}(d+1)$.
\item[(iii)] $(\H\up)_k = (\H_{k-1})\up$ for every $k \geq 0$.
\ei
\ep

\proof
(i) Immediate.

(ii) In view of Theorem \ref{upp}.

(iii) Since
$$\begin{array}{lll}
(H\up)_k&=&H_k \cup \{ I \cup \{ p \} \mid I \in H_{k-1},\; p \in V \setminus I\}\\
&=&H_{k-1} \cup \{ I \cup \{ p \} \mid I \in H_{k-1},\; p \in V \setminus I\} = (H_{k-1})\up.
\end{array}$$ 
\qed

A simplicial complex $\H\; = (V,H)$ is said to be {\em graphic
  boolean} if it can be represented by a boolean matrix $M$ such that:
\bi
\item
$M$ contains all possible rows with one zero;
\item
each row of $M$ has at most two zeroes.
\ei

In \cite[Section 6.2]{RSm}, it is remarked that either $\H\; = U_{2,|V|}$ or
$\H\; \in \bpav(2)$. However, $\H$ is not necessarily a matroid.
We can of course represent $\H$ by a graph where the edges
correspond to the flats of the matrix having precisely two elements. 

We view graphs as simplicial complexes $\Gamma = (V,E)$ with $P_1(V) \subseteq E \subseteq P_{\leq 2}(V)$. Now the following result follows from the definitions.

\bp
\label{gupgb}
The following conditions are equivalent for a simplicial complex $\H$:
\bi
\item[(i)] $\H$ is graphic boolean;
\item[(ii)] $\H = \Gamma\up$ for some graph $\Gamma$.
\ei
\ep

Other interesting examples connected to graphs arise from considering some class $\G$ of graphs closed under removing an edge. Edgeless graphs, forests, triangle-free graphs or graphs containing no cycle of length $\leq \ell$ (for a fixed $\ell \geq 3$) consitute examples of such classes.

Given a graph $\Gamma = (V,E)$, we define a simplicial complex $\H_{\G} = (V,H_{\G})$ where $H_{\G}$ is the set of all $W \subseteq V$ such that the subgraph of $\Gamma$ induced by $W$ belongs to $\G$. 

In general, $\H_{\G}$ is not a BRSC. If $\G$ is the class of all edgeless graphs, then the faces of $\H_{\G}$ are the anticliques of $\Gamma$, and (PR) is not always satisfied by $\H_{\G}$. However, $(\H_{\G})\up$ is a BRSC by Theorem \ref{upp}, and these construction will merit our attention in future sections.

We consider now the iteration of the up operator. Given $m \geq 0$, we denote by $\H\upp{m} = (V,H\upp{m})$ the simplicial complex obtained by applying $m$ times the up operator to $\H\; =(V,H)$. 

\bt
\label{itup}
Let $\H\; = (V,H) \in {\rm Pav}(d)$ and $m \geq 0$. Then
\beq
\label{itup1}
H\upp{m} = P_{\leq d+m+1}(V) \setminus \{ X \in P_{d+m+1}(V) \mid P_{d+1}(X) \cap H = \emptyset\}.
\eeq
\et

\proof
We use induction on $m$. The case $m = 0$ is easy to check, hence we assume that $m > 0$ and (\ref{itup1}) holds for $m-1$. We have $\H\upp{m} \in \bpav(d+m)$ by Proposition \ref{pup}(ii), hence we only need to discuss membership on both sides of (\ref{itup1}) for a fixed $Y \in P_{d+m+1}$.

Assume first that $Y \in H\upp{m}$. Then there exists some $y \in Y$ such that $Y \setminus \{ y \} \in H\upp{m-1}$. By the induction hypothesis, we get 
$$Y \setminus \{ y \} \in P_{\leq d+m}(V) \setminus \{ X \in P_{d+m}(V) \mid P_{d+1}(X) \cap H = \emptyset\}.$$
Hence there exists some $Z \in P_{d+1}(Y \setminus \{ y \}) \cap H \subseteq P_{d+1}(Y) \cap H$. Thus $Y \in P_{\leq d+m+1}(V) \setminus \{ X \in P_{d+m+1}(V) \mid P_{d+1}(X) \cap H = \emptyset\}$.

Suppose now that $Y \notin H\upp{m}$. By (\ref{hup}), we have $Y \setminus \{ y\} \notin H\upp{m-1}$ for every $y \in Y$. By the induction hypothesis, we get  
$P_{d+1}(Y \setminus \{ y \}) \cap H  = \emptyset$ for every $y \in Y$. But then $P_{d+1}(Y) \cap H  = \emptyset$ and so (\ref{itup1}) holds for $m$ as required.
\qed

\bc
\label{citup}
Let $\Gamma\; = (V,E)$ be a graph and $m \geq 0$. Then
$$\Gamma\upp{m} = P_{\leq m+2}(V) \setminus \{ X \in P_{m+2}(V) \mid X\mbox{ is an anticlique of }\Gamma \}.$$
\ec

It follows that $\Gamma\upp{m}$ is an uniform matroid if the greatest anticlique of $\Gamma$ has at most $m+1$ elements.
 
For all $k \leq n$, consider the uniform matroid $U_{k,n} = (V_n,P_{\leq k}(V_n))$. 

By a classical Ramsey theory theorem (see \cite{GRS}), for every $m \geq 2$ there exists some integer $R(m)$ such that every graph with at least $R(m)$ vertices admits an $m$-clique or an $m$-anticlique. This yields a corollary to Corollary \ref{citup}:

\bc
\label{ccitup}
Let $m \geq 2$ and let $\Gamma\; = (V,E)$ be a graph with $|V| \geq R(m)$. Then
$\Gamma\up$ has a restriction isomorphic to $U_{3,m+1}$ or $\Gamma\upp{m-2}$ is not uniform.
\ec

In terms of shellability, it remains an open problem whether or not $\H\up$ is shellable when $\H$ is shellable. For the converse implication, we have the following counterexample:

\be
\label{exs}
Let $\H\; = (V_5,H)$ be the pure complex defined by $H = 2^{123} \cup 2^{345}$. Then $H\up$ is shellable but $\H$ is not.
\ee

Indeed, $\fct\H\; = \{ 123, 345 \}$ and so $\H$ is not shellable. However, $\fct\H\; = \{ 1234, 1235, 1345, 2345 \}$ and this particular enumeration provides a shelling for $\H\up$.

\subsection{Unimodality}

Given a simplicial complex $\H\; = (V,H)$, we define the {\em counting
function} of $\H$ to be the function $\alpha_{\H}:\N \to \N$ defined
by
$$n\alpha_{\H} = |H \cap P_n(V)|.$$
A function $\alpha:\N \to \N$ is {\em unimodal} if there exists some
$m \in \N$ such that
\bi
\item
$i < j \leq m \Rw i\alpha \leq j\alpha$;
\item
$m \geq i < j \Rw i\alpha \geq j\alpha$.
\ei

Giancarlo Rota conjectured that $\alpha_{\H}$ is unimodal for every matroid $\H$. We discuss next this conjecture for BRSCs.

For every $n \geq 1$, write $V_n = \{ 1,\ldots, n \}$.

For all integers $n \geq m \geq k \geq 1$, let $\J_{n,m,k} = (V_n,H)$ be
the simplicial complex defined by
$$H = P_1(V_n) \cup P_{\leq k}(V_m).$$
It follows that $\J_{n,m,k}\up = (V_n,H\up)$ with
\beq
\label{fhup}
H\up = P_{\leq 2}(V_n) \cup P_{\leq k+1}(V_m) \cup \{ I \cup \{ p \}
\mid I \in P_{\leq k}(V_m), \; p \in V_n \setminus V_m \}.
\eeq

\bp
\label{uph}
For all integers $n \geq m \geq k \geq 1$, let $\J_{n,m,k} = (V_n,H)$ and $\J = \J_{n,m,k}\up$. Then ${\rm Fl}\J = H \cup \{ V_n \}$.
\ep

\proof
In view of Corollary \ref{fup} and the definition of flat, it suffices to show that every $X \in 2^{V_n} \setminus (H \cup \{ V_n \})$ contains a facet of $\H$.

Assume first that there exists some distinct $q,r \in X \setminus V_m$. Then $X$ contains the facet $qr$.

Assume now that $|X \setminus V_m| \leq 1$. Then $X \notin H$ yields $|X \cap V_m| > k$. Then any $(k+1)$-subset of $X$ is a facet.
\qed

The following example features a BRSC of dimension 2 with a non unimodal counting function:

\be
\label{cex1}
Let $\J = \J_{16,6,3}$. Then 
$$1\alpha_{\H} = 16,\quad 2\alpha_{\H} = \binom{6}{2} = 15,\quad 3\alpha_{\H} = \binom{6}{3} = 20.$$
\ee

However, $\J_{16,6,3}$ is not connected. We discuss now the connected case.

\bp
\label{ctu}
Let $\H\; = (V,H)$ be a connected BRSC of dimension $\leq 2$. Then $\alpha_{\H}$ is unimodal.
\ep

\proof
We may assume that $\dim\H\; = 2$ and $1\alpha_{\H} = n$. Since the graph $\H_1$ has a 3-cycle, it is a connected graph on $n$ vertices which is not a tree. Hence 
$$2\alpha_{\H} \geq 1\alpha_{\H} > 1 = 0\alpha_{\H}$$
and so $\alpha_{\H}$ is unimodal.
\qed

We discuss now the counting function of $\J = \J_{n,m,3}\up$. In view of (\ref{fhup}), we have
\bi
\item
$2\alpha_{\J} = \binom{n}{2} = \frac{n(n-1)}{2}$;
\item
$3\alpha_{\J} = \binom{m}{3} + \binom{m}{2}(n-m) = \frac{m(m-1)(m-2)}{6} + \frac{m(m-1)(n-m)}{2}$;
\item
$4\alpha_{\J} = \binom{m}{4} + \binom{m}{3}(n-m) = \frac{m(m-1)(m-2)(m-3)}{24} + \frac{m(m-1)(m-2)(n-m)}{6}$.
\ei
To ensure $3\alpha_{\J} < 4\alpha_{\J}$, we must have $6 \leq m \leq n-2$, and it is straightforward to check that for $m = 6$ we have $2\alpha_{\J} > 3\alpha_{\J}$ if and only if $n \geq 26$. In view of Theorem \ref{upp}, we obtain:

\be
\label{cex2}
Let $\J = \J_{26,6,3}\up$. Then $\J$ is a simple BRSC of dimension 3 such that $\alpha_{\J}$ is not unimodal. Moreover, this is the smallest counterexample among complexes of the form $\J(n,m,k)$.
\ee

It is easy to check that $\alpha_{\H\up}$ can be unimodal even if $\alpha_{\H}$ is not, $\H = \J(16,6,3)$ providing a straightforward example. It remains an open question whether $\alpha_{\H}$ unimodal implies $\alpha_{\H\up}$ unimodal. We can show that the implication holds for complexes of the form $\H\; = \J(n,m,3)$, but the general case remains open.

\section{Adding one point}

In this section we discuss alternative constructions designed to add one extra vertex to a complex.

\subsection{The operator $\H + p$}

Given $\H\; = (V,H)$ and $\H' = (V',H')$, we write $$\H \leq \H' \mbox{ if $V \subseteq V'$ and }H \subseteq H'$$
and we define the {\em join} of the simplicial complexes $\H\; = (V,H)$ and $\H' = (V',H')$ to be the simplicial complex
$$(\H \vee \H') = (V\cup V', H \cup H').$$

Given a simplicial complex $\H\; = (V,H)$ and $p \notin V$, we define
$$\H + p = (V \cup \{ p \},H \cup \{ p \}).$$
Then $\H +p = (\H \vee \S_p)$, where $\S_p$ is the unique simplicial complex having $p$ as its unique vertex. 

It is easy to see that $\H + p$ satisfies the point replacement property (PR) if and only if $\dim\H \leq 0$. Since every BRSC satisfies (PR), it is easy to check that $\H +p$ only is a BRSC if $\dim\H \leq 0$.

The next result relates $\H +p$ with the up operator and contraction.

\bp
\label{pupc}
Let $\H\; = (V,H)$ be a simplicial complex and $p \notin V$. Then $(\H +p)\up/p = \H$.
\ep

\proof
We have
$$(H\cup \{ p \})\up = H\up \cup \{ pq \mid q \in V\} \cup \{ I \cup \{ p \} \mid I \in H\},$$
thus the set of faces of $(\H +p)\up/p$ is $P_1(V) \cup H = H$.
\qed

We use Proposition \ref{pupc} in the discussion of the relationship between contractions and the up operator.

\bp
\label{resup}
Let $\H\; = (V,H)$ be a simplicial complex and let $W \subseteq V$. Then $\H\up/W \leq (\H|_{V\setminus W})\up$. Equality does not hold in general.
\ep

\proof
It follows from the definitions that $V\setminus W$ is the vertes set of both complexes. Let $X \subseteq V\setminus W$. If $X$ is a face of $\H\up/W$, then $X \cup W \in H\up$ and so $X \cup W \in H$ or $X \cup W = I \cup \{ p \}$ for some $I \in H$ and $p \in V \setminus I$. In the first case, we get $X \in H|_{V\setminus W}$, in the second $X \setminus \{ p \} \in H|_{V\setminus W}$, so in any case $X$ is a face of $(\H|_{V\setminus W})\up$. Therefore $\H\up/W \leq (\H|_{V\setminus W})\up$.

Suppose now that $V \notin H$ and $p \notin V$. Since $V = (V \cup \{ p \})\setminus \{ p \}$, to complete the proof it is enough to show that
\beq
\label{resup1}
(\H+p)\up/p < ((\H+p)|_{V})\up.
\eeq

By Proposition \ref{pupc}, we have $(\H+p)\up/p = \H$. On the other hand, $((\H+p)|_{V})\up = \H\up$. Since $V \notin H$, we have $\H < \H\up$ and so (\ref{resup1}) holds as required.
\qed

\subsection{The operator $\H \oplus p$}

Let $V$ and $V'$ be disjoint sets. Given $\X \subseteq 2^V$ and $\X' \subseteq 2^{V'}$, we write
$$\X\oplus \X' = \{ X \cup X' \mid X \in \X,\; X' \in \X' \}.$$  

Given simplicial complexes $\H\; = (V,H)$ and $\H' = (V',H')$ with $V \cap V' = \emptyset$, we define the simplicial complex 
$$\H \oplus \H' = (V \cup V',H\oplus H').$$

\bl
\label{pop}
Let $\H\; = (V,H)$ and $\H' = (V',H')$ be simplicial complexes with $V \cap V' = \emptyset$. Then:
\bi
\item[(i)]
${\rm Fl}(\H \oplus \H') = {\rm Fl}\H \oplus {\rm Fl}\H'$;
\item[(ii)]
$\H \oplus \H'$ is boolean representable if and only $\H$ and $\H'$ are both boolean representable.
\ei
\el

\proof
(i) $\Rw$ (ii). Let $C \in \flatx(\H \oplus \H')$. By \cite[Proposition 8.3.3(i)]{RSm}, $C \cap V$ is a flat of the restriction of $\H \oplus \H'$ to $V$, i.e. 
$C\cap V \in \flatx\H$. Similarly, $C\cap V' \in \flatx\H'$. Hence 
$$C = (C\cap V)\cup (C \cap V') \in {\rm Fl}\H \oplus {\rm Fl}\H'$$
and so ${\rm Fl}(\H \oplus \H') \subseteq {\rm Fl}\H \oplus {\rm Fl}\H'$.

Conversely, let $F \in {\rm Fl}\H$ and $F' \in {\rm Fl}\H'$. Let $X \in (H\oplus H') \cap 2^{F\cup F'}$ and $p \in (V\cup V')\setminus (F\cup F')$. We must show that $X \cup \{ p \} \in H\oplus H'$. Without loss of generality, we may assume that $p \in V$.

Write $X = I \cup I'$ with $I \in H$ and $I' \in H'$. Then $I \in H \cap 2^F$ and $p \in V\setminus F$. Since $F \in {\rm Fl}\H$, we get $I \cup \{ p \} \in H$ and so 
$X \cup \{ p \} = (I \cup \{ p\}) \cup I' \in H\oplus H'$ as required.
\qed

Given a simplicial complex $\H\; = (V,H)$ and $p \notin V$, we write 
$$\H \oplus p = (V \cup \{ p \},H \cup \{ I \cup \{ p \} \mid I \in H \}).$$
Thus $\H\oplus p = \H \oplus \S_p$.

\bp
\label{phop}
Let $\H\; = (V,H)$ be a simplicial complex and let $p \notin V$. Then:
\bi
\item[(i)]
${\rm dim}(\H \oplus p) = {\rm dim}\H +1$; 
\item[(ii)]
${\rm Fl}(\H \oplus p) = {\rm Fl}\H \cup \{ F \cup \{ p \} \mid F \in {\rm Fl}\H \}$;
\item[(iii)]
$\H$ is boolean representable if and only $\H \oplus p$ is boolean representable;
\item[(iv)]
$\H$ satisfies (PR) if and only $\H \oplus p$ satisfies (PR);
\item[(v)]
$\displaystyle\H\up = \bigvee_{p \in V} (\H|_{V\setminus\{ p\}} \oplus p).$
\ei
\ep

\proof
(i), (iv). Immediate.

(ii), (iii). By Lemma \ref{pop}.

(v). Both complexes have vertex set $V$. Now $X \subseteq V$ is a face of the right hand side complex if and only if $X$ is a face of $\H|_{V\setminus\{ p\}} \oplus p$ for some $p \in V$. But $X$ is a face of $\H|_{V\setminus\{ p\}} \oplus p$ if and only if $X = I$ or $X = I \cup \{ p\}$ for some $I \in H \cap 2^{V\setminus\{ p\}}$. Since $p$ is arbitrary, we obtain precisely the faces of $\H\up$.
\qed

\subsection{The operator $\H \boxplus p$}

Given a Moore family $R \subseteq 2^V$ and $p \notin V$, we define the Moore family 
$$R\boxplus p = (R \setminus \{ V \}) \cup \{ V \cup \{ p\}, \{ p\} \} \subseteq 2^{V\cup \{ p\}}.$$
This Moore family determines the BRSC $\J(R\boxplus p)$, whose faces are the transversals of the successive differences for some chain of $R\boxplus p$.

If 
$$M(R) = \left(
\begin{matrix}
0&&0&&\ldots&&0\\
a_{11}&&a_{12}&&\ldots&&a_{1n}\\
a_{21}&&a_{22}&&\ldots&&a_{2n}\\
\vdots&&\vdots&&\ddots&&\vdots\\
a_{m1}&&a_{m2}&&\ldots&&a_{mn}\\
\end{matrix}
\right)$$
then
$$M(R\boxplus p) = \left(
\begin{matrix}
0&&0&&\ldots&&0&&0\\
a_{11}&&a_{12}&&\ldots&&a_{1n}&&1\\
a_{21}&&a_{22}&&\ldots&&a_{2n}&&1\\
\vdots&&\vdots&&\ddots&&\vdots&&\vdots\\
a_{m1}&&a_{m2}&&\ldots&&a_{mn}&&1\\
1&&1&&\ldots&&1&&0\\
\end{matrix}
\right)$$

The following example shows that $J(R\boxplus p)$ is not determined by $J(R)$.

\be
\label{jj}
Let $V = 1234$ and consider the Moore families 
$$R = \{ \emptyset, 1,2,123,V\}, \quad R' = \{ \emptyset, 4,14,24,V\}.$$
Then $J(R) = J(R')$ but $J(R\boxplus p) \neq J(R'\boxplus p)$.
\ee

Indeed, it follows from \cite[Proposition 5.7.2]{RSm} that $J(R) = J(R') = P_{\leq 3}(V) \setminus \{ 123 \}$. However, $125 \in J(R\boxplus p) \setminus J(R'\boxplus p)$.

\medskip

Given a BRSC $\H\; = (V,H)$ and $p \notin V$, we define the BRSC
$$\H\boxplus p = \J(\flatx\H \boxplus p) = (V \cup \{ p\}, J(\flatx\H \boxplus p)).$$

\bp
\label{hthird}
Given a BRSC $\H\; = (V,H)$ and $p \notin V$, we have
$$J({\rm Fl}\H \boxplus p) = H \cup \{ I \cup \{ p \} \mid I \in H, \; \oo{I} \neq V, \; q \in I\} \cup \{ vp \mid v \in V\},$$
where $\oo{I}$ denotes the closure of $I$ in ${\rm Fl}\H$.
\ep

\proof
The nonempty elements of $J({\rm Fl}\H \boxplus p)$ are the transversals of the successive differences for chains in $\flatx\H \boxplus p$, where we may assume that $V \cup \{ p\}$ is the biggest element. For a chain of the type
$$F_0 \subset F_1 \subset \ldots \subset F_k = V \cup \{ p \},$$
consider a transversal $x_1x_2\ldots x_{k}$ with $x_i \in F_i\setminus F_{i-1}$ for every $i$.
 
If $F_{k-1} = \{ p\}$, we get as transversals the elements in $P_1(V) \cup \{ vp \mid v \in V\}$ (for $k = 0$ and $k = 1$, respectively). Assume now that $F_{k-1} \neq \{ p\}$. If $x_k \neq p$, then we get as transversals all the nonempty elements of $H$. If $x_k = p$, we get faces of the form $x_1\ldots x_{k-1}p$, with $\oo{x_1\ldots x_{k-1}} \subseteq F_{k-1} \subset V$.
\qed

Note that $\dim(\H\boxplus p) = \dim\H$ unless $\dim\H\; = 0$.

\section{Computing the flats of paving complexes}

We present in this section a method to compute the flats of a paving complex.

We recall the following result.

\bp
\label{cfac}
{\rm \cite[Proposition 4.2.3]{RSm}}
Let $\H\; = (V,H)$ be a simplicial complex. Then $\oo{B} = V$ for every facet of $\H$.
\ep

Let $\H\; = (V,H)$ be a paving simplicial complex of dimension $d \geq 2$. Then  facets may have dimension $d$ or $d-1$. A facet of dimension $d$ is called {\em large}, otherwise it is {\em small}.

We say that $X \in P_{\geq d+1}(V)$ is a {\em long hyperplane} of $\H$ ({\em long line} if $d = 2$) if $X$ contains no facet of $\H$. The set of {\em maximal long hyperplanes} (with respect to inclusion) is denoted by $\L_{\H}$.

\bl
\label{cll}
Let $\H\; = (V,H) \in {\rm Pav}(d)$ with $d \geq 2$. For every $L \in \L_{\H}$, we have either $\oo{L} = L$ or $\oo{L} = V$.
\el

\proof
Assume $\oo{L} \neq V$. By Proposition \ref{cfac}, we have $\oo{B} = V$ for every $B \in \fct\H$. Hence $\oo{L}$ can contain no facet and is thus a long hyperplane containing $L$. Therefore $\oo{L} = L$ by maximality of $L$.
\qed

In view of Lemma \ref{cll}, we consider now a partition 
\beq
\label{partlh1}
\L_{\H} = \L_{\H}^{(1)} \cup \L_{\H}^{(2)} \cup \L_{\H}^{(3)}
\eeq
defined as follows:
\bi
\item
$\L_{\H}^{(1)} = \{ L \in \L_{\H} \cap \flats\H \; {\big{\lvert}} \; |L \cap L'| \leq d-1 \mbox{ for every }L' \in \L_{\H} \setminus \{ L\} \}$;
\item
$\L_{\H}^{(2)} = \{ L \in \L_{\H} \setminus \flats\H \; {\big{\lvert}} \; |L \cap L'| \leq d-1 \mbox{ for every }L' \in \L_{\H} \setminus \{ L\} \}$;
\item
$\L_{\H}^{(3)} = \{ L \in \L_{\H} \setminus \flats\H \; {\big{\lvert}} \; |L \cap L'| \geq d \mbox{ for some }L' \in \L_{\H} \setminus \{ L\} \}$.
\ei
We note that
\beq
\label{partlh}
\oo{L} = V\mbox{ for every }L \in \L_{\H}^{(2)} \cup \L_{\H}^{(3)}.
\eeq
Indeed, since $L \notin \flats\H$, then $\oo{L}$ strictly contains $L$ and thus $\oo{L} = V$ by Lemma \ref{cll}.
Thus (\ref{partlh}) holds.

The following example shows that the partition (\ref{partlh1}) may be nondegenerate, even for dimension $d = 2$.

\be
\label{lhne}
Let $\H\; = (V,H)$ be defined by $V = \{ 0,1,\ldots,9\}$ and 
$$H = P_{\leq 3}(V) \setminus (\{ 123,345,789,890\} \cup \{ 56p \mid p \in V \setminus \{ 5,6\} \}).$$
Then
$$\L_{\H}^{(1)} = \{ 123 \}, \quad \L_{\H}^{(2)} = \{ 345\}, \quad \L_{\H}^{(3)} = \{ 789, 890\}.$$
\ee

Indeed, it is easy to check that these four 3-sets are the unique long hyperplanes of $\H$ (and therefore maximal). It follows immediately that $\L_{\H}^{(3)} = \{ 789, 890\}$ and $123 \in \L_{\H}^{(1)}$. Since $45 \in H$ but $456 \notin H$, we get $6 \in \oo{345}$, hence $\oo{345}$ contains the facet 56 and so $\oo{345} = V$ by Proposition \ref{cfac}. Therefore $345 \in \L_{\H}^{(2)}$ as claimed.

\bq
Can we find such an example with $\H$ boolean representable?
\eq

The next result shows that for some classes of complexes, the partition (\ref{partlh1}) may be degenerate.

\bp
\label{ppur}
Let $\H\; = (V,H) \in {\rm Pav}(d)$ with $d \geq 2$. If $\H$ is pure, then $\L_{\H}^{(2)} = \emptyset$.
\ep

\proof
Let $L \in \L_{\H} \setminus \flatx\H$. Since $L$ contains no facet and $P_{\leq d}(V) \subseteq H$, then there exists some $A \in H \cap P_d(L)$ and $p \in V \setminus L$ such that $A \cup \{ p \} \notin H$. Since $\H$ is pure, it follows that $A \cup \{ p \}$ contains no facet of $\H$ and is consequently a long hyperplane. Thus $A \cup \{ p \} \subseteq L'$ for some $L' \in \L_{\H}$. But $p \in L' \setminus L$ and $A \subseteq L \cap L'$ yields $|L\cap L'| \geq d$, hence $L \in \L_{\H}^{(3)}$ and so $\L_{\H}^{(2)} = \emptyset$.
\qed

Now we can describe all the flats of $\H \in \pav(d)$.

\bt
\label{cfp}
Let $\H\; = (V,H) \in {\rm Pav}(d)$ with $d \geq 2$. Then
\beq
\label{cfp1}
{\rm Fl}\H = P_{\leq d-1}(V) \cup \{ A \in P_d(V) \mid \forall p \in V \setminus A\; A \cup \{ p \} \in H\} \cup \L_{\H}^{(1)} \cup \{ V \}.
\eeq
\et 

\proof
It is easy to check that the right hand side of (\ref{cfp1}) consists solely of flats. Conversely, let $F \in \flatx\H$. We may assume that $d \leq |F| < |V|$. 

Suppose first that $|F| = d$. Since $P_d(V) \subseteq H$, it follows that $A \cup \{ p \} \in H$ for every $p \in V \setminus A$, hence we may assume that $d < |F| < |V|$. 

In view of Proposition \ref{cfac}, $F$ contains no facet of $\H$, hence it is a long hyperplane. Suppose that $p \in V \setminus F$ and take $A \in P_d(F) \subseteq H$. Since $F \in \flatx\H$, we get $A \cup \{ p\} \in H$ (therefore a facet), hence $F \cup \{ p\}$ is not a long hyperplane. Thus $F \in \L_{\H}$ and so $F \in \L_{\H}^{(1)} \cap \L_{\H}^{(3)}$. In view of (\ref{partlh}), we get $F \in \L_{\H}^{(1)}$ and so (\ref{cfp1}) holds as required.
\qed

Note that, by \cite[Proposition 4.2.4]{RSm}, every $F \in \flatx\H$ is of the form $\oo{I}$ for some $I \in H$. Combined with Proposition \ref{cfac}, this implies that every $F \in \flatx\H \setminus \{ V\}$ is of the form $\oo{I}$ for some $I \in H \cap P_{\leq d}(V)$. In the particular case $d = 2$, we generate all the flats (except possibly $V$) as the closure of a set with at most 2 elements.

\brem
\label{tess}
Let $\H\; = (V,H) \in {\rm Pav}(d)$ with $d \geq 2$. Let
$$J = P_{\leq d-1}(V) \cup \{ A \in P_d(V) \mid \forall p \in V \setminus A\; A \cup \{ p \} \in H\}.$$
Then:
\bi
\item[(i)] $J = {\rm Fl}\H \, \cap P_{\leq d}(V)$;
\item[(ii)] $J\up \subseteq H$;
\item[(ii)] if $(V,K)$ is a simplicial complex such that $K\up \subseteq H$, then $K \subseteq J$.
\ei
\erem

\proof
(i) By Theorem \ref{cfp}.

(ii) Let $A \in J$ and let $p \in V$. If $|A| < d$, then $A \cup \{ p \} \in H$ since $P_{d}(V) \subseteq H$. Hence we may assume that $|A| = d$. Then  $A \cup \{ p \} \in H$ if $p \notin A$ by definition of $J$. If $p \in A$, we get $A \cup \{ p \} = A \in P_d(V) \subseteq H$, therefore $A \cup \{ p \} \in H$ in any case as required.

(iii) Let $A \in K$. Since $K\up \subseteq H$, we have $|A| \leq d$. Since $P_{d-1}(V) \subseteq J$, we may assume that $|A| = d$. Let $p \in V \setminus A$. Since $K\up \subseteq H$, we have $A \cup \{ p \} \in H$. Hence $A \in J$ and so $K \subseteq J$.
\qed

\section{Truncation}
\label{stru}

Given a simplicial complex $\H\; = (V,H)$ and $k \geq 1$, the $k$-{\em
  truncation} of $\H$ is the
simplicial complex $\H_k\; = (V,H_k)$
defined by $H_k = H \cap P_{\leq k}(V)$. 

We say that a simplicial complex $\H\; = (V,H)$ is a {\em TBRSC} if $\H\; = \J_k$ for some BRSC $\J$ and $k \geq 1$. For every $d \geq 1$, we denote by $\tbpav(d)$ the class of all paving TBRSCs of dimension $d$.

It is known that not every simplicial complex is a TBRSC \cite[Example 8.2.6]{RSm} and not every TBRSC is a BRSC \cite[Example 8.2.1]{RSm}.

To understand TBRSCs, we need the following definition.

Given a simplicial complex $\H\; = (V,H)$ of dimension $d$, we define
$$T(H) = 
\{ T \subseteq V \mid \forall X \in H_{d} \cap 2^T\; \forall p \in V
\setminus T \hspace{.3cm} X \cup \{ p \} \in H\}.$$

\bl
\label{propt}
{\rm \cite[Lemma 8.2.3]{RSm}}
Let $\H\; = (V,H)$ be a simplicial complex. Then:
\bi
\item[(i)] $T(H)$ is closed under intersection;
\item[(ii)] {\rm Fl}$\H\; \subseteq T(H)$.
\ei
\el

Thus $T(H)$ is a Moore family and $\J(T(H))$ is a BRSC. 

\bt
\label{eqtr}
{\rm \cite[Theorem 8.2.5]{RSm}}
Let $\H\; = (V,H)$ be a simplicial complex of dimension $d$. Then the following
conditions are equivalent:
\bi
\item[(i)] $\H$ is a TBRSC; 
\item[(ii)] $\H = (\J(T(H)))_{d+1}$.
\ei
Furthermore, in this case we have {\rm Fl}$\J(T(H)) = T(H)$.
\et

\subsection{Low dimensions}

\bp
\label{tbone}
Every TBRSC of dimension 1 is boolean representable.
\ep

\proof
Let $\H\; = (V,H)$ be a TBRSC of dimension 1. Let
$$\Gamma\H\; = (V, P_{\leq 1}(V) \cup (P_2(V) \setminus H)).$$ 
By \cite[Proposition 5.3.1]{RSm}, $\H$ is a BRSC if and only if the connected components of $\Gamma\H$ are cliques.

Let $a \edge b \edge c$ be distinct edges in $\Gamma\H$. Suppose that $\Gamma\H$ has no edge $a \edge c$. Then $ac \in H \cap P_2(V)$ and so by Theorem \ref{eqtr} there exists some $T \in T(H)$ such that $|T \cap ac| = 1$. We may assume that $a \in T$, but then $ab \notin H$ yields $b \in T$ and $bc \notin H$ yields $c \in T$, a contradiction. Thus $a \edge c$ is also an edge of $\Gamma\H$ and so 
the connected components of $\Gamma\H$ are cliques. Therefore $\H$ is a BRSC.
\qed

The following example shows that Proposition \ref{tbone} fails for dimension 2, even in the paving case.

\be
\label{btbtwo}
Let $V = 123456$ and 
$$H = P_{\leq 2}(V) \cup \{ X \in P_3(V) \; {\big{\lvert}} \; |X \cap 56| = 1\} \cup \{ 123,124 \}.$$
Then $\H\; = (V,H) \in {\rm TBPav}(2) \setminus {\rm BPav}(2)$.
\ee

Indeed, it is easy to check that $\flatx\H$ is the lattice
$$
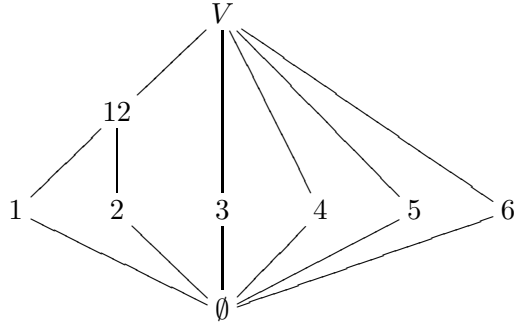
\xymatrix{
&&V \ar@{-}[dl] &&&\\
&12 &&&&\\
1 \ar@{-}[ur] & 2 \ar@{-}[u] & 3 \ar@{-}[uu]
& 4 \ar@{-}[uul] & 5 \ar@{-}[uull] & 6
\ar@{-}[uulll] \\
&&\emptyset \ar@{-}[ull] \ar@{-}[ul] \ar@{-}[u] \ar@{-}[ur] \ar@{-}[urr]
\ar@{-}[urrr] &&&
}$$
and so $\H \notin \bpav(2)$ since 135 is not a transversal of any chain in $\flatx\H$.

On the other hand, $T(H)$ is the lattice
$$
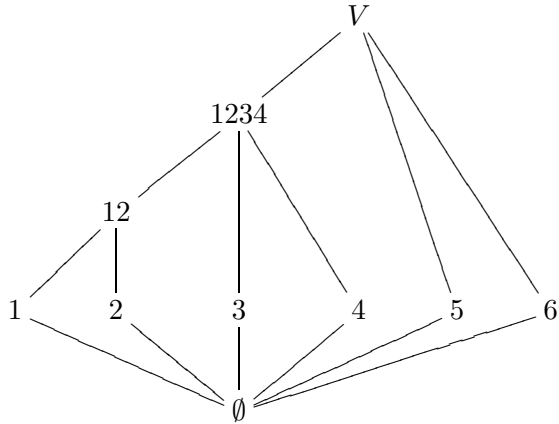
\xymatrix{
&&&V \ar@{-}[dl] &&\\
&&1234 \ar@{-}[dl] &&&\\
&12 &&&&\\
1 \ar@{-}[ur] & 2 \ar@{-}[u] & 3 \ar@{-}[uu]
& 4 \ar@{-}[uul] & 5 \ar@{-}[uuul] & 6
\ar@{-}[uuull] \\
&&\emptyset \ar@{-}[ull] \ar@{-}[ul] \ar@{-}[u] \ar@{-}[ur] \ar@{-}[urr]
\ar@{-}[urrr] &&&
}$$
and it is easy to check that $\H = (\J(T(H)))_{3}$. Therefore $\H \in \tbpav(2)$.

\bigskip

The next result shows that, when it comes to separate BRSCs from TBRSCs, the above example has the minimum number of vertices.

\bp
\label{btbfive}
Every TBRSC with at most 5 vertices is boolean representable.
\ep

\proof
Let $\H\; = (V,H)$ be a TBRSC with $|V| \leq 5$.
In view of Proposition \ref{tbone}, we may assume that $|V| \geq 4$ and $\dim\H \geq 2$.

Suppose first that $|V| = 4$. We may assume that $\dim\H\; = 2$. If $\H\: = \J_3$ for some BRSC $\J\; = (V,J)$ of dimension 4, then $\H$ is the uniform matroid $U_{3,4}$ and is therefore a BRSC.

Thus we may assume that $|V| = 5$. If $\dim\H\; = 3$ and $\H\: = \J_4$ for some BRSC $\J\; = (V,J)$ of dimension 4, then $\H$ is the uniform matroid $U_{4,5}$ and so a BRSC. Hence we may assume that $\dim\H\; = 2$.

Suppose that $\H$ is not a BRSC. Then $\flatx\H \subset T(H)$. Let $T \in T(H) \setminus \flatx\H$. Then $H \cap P_3(T) \neq \emptyset$, hence we may take $a_1a_2a_3 \in H \cap P_3(T)$. Since $a_1a_2a_3 \in J(T(H))$, we may assume that there exists a chain $\emptyset = T_0 \subset T_1 \subset T_2 \subset T_3$ in $T(H)$ such that $a_i \in T_i \setminus T_{i-1}$ for $i = 1,2,3$. By Lemma \ref{propt}(i), $a_1a_2a_3$ is also a transversal of the successive differences for the chain
$$\emptyset \subset T_1 \cap T \subset T_2 \cap T \subset T_3 \cap T$$ 
in $T(H)$, hence there exists a chain 
$$\emptyset = T'_0 \subset T'_1  \subset T'_2 \subset T'_3 \subset T'_4 = V$$ 
in $T(H)$. Since $|V| = 5$, there exists some $j \in 1234$ such that $|T'_j \setminus T'_{j-1}| = 2$. Without loss of generality, we may assume that 
$T'_j = T'_{j-1} \cup 12$. It is easy to check that 
$$abc \in J(T(H)) \mbox{ if } |abc \cap 12| \leq 1,$$
hence the only possible elements of $P_3(V) \setminus H = P_3(V) \setminus (J(T(H)))_3$ are $123, 124, 125$. 

If $12 \notin H$, we have necessarily 
$$H = \{ X \in P_{\leq 3}(V) \mid 12 \not\subseteq X \},$$
hence we have a matroid (therefore boolean representable).

Thus we may assume that we have one of the following four cases:
\bi
\item[(C1)] $H = P_{\leq 3}(V) \setminus \{ 123, 124, 125 \}$;
\item[(C2)] $H = P_{\leq 3}(V) \setminus \{ 123, 124\}$;
\item[(C3)] $H = P_{\leq 3}(V) \setminus \{ 123 \}$;
\item[(C4)] $H = P_{\leq 3}(V)$.
\ei
Now (C3) and (C4) are both matroids (hence BRSCs). We can show that (C1) is a BRSC by checking that $34,35,45$ are flats. Similarly, (C2) is a BRSC because $15, 34,35,45$ are flats. Therefore every TBRSC with 5 vertices is a BRSC.
\qed

\subsection{Union}

Given two simplicial complexes $\H\; = (V,H)$ and $\H' = (V,H')$ we define the {\em union} of $\H$ and $\H'$ as the simplicial complex
$$\H \cup \H' = (V,H \cup H').$$
This is the join of complexes restricted to the case when the vertex sets coincide. 

\bp
\label{unfo}
Let $\H\; = (V,H)$ and $\H' = (V,H')$ be BRSCs with $|V| \leq 4$. Then $\H \cup \H'$ is a BRSC.
\ep

\proof
If $\dim(\H \cup \H') = 1$, we may use Theorem \ref{pavun}. The only other nontrivial case is $\dim(\H \cup \H') = 2$. But it is easy to check that if $|V| = 4$ and $\dim\H\; = 2$, then $\H$ is a BRSC if and only if $|H \cap P_3(V)| \neq 1$. It follows that if $\H \cup \H'$ is not a BRSC, then $\H$ or $\H'$ is not a BRSC.
\qed

The next example shows that neither BRSCs nor TBRSCs are closed under union when we consider 5 vertices (even at dimension $\leq 2$).

\be
\label{nonun}
Let $V = 12345$ and let $\H_1 = (V,P_{\leq 2}(V))$ and 
$$\H_2 = (V, P_{\leq 1}(V) \cup \{ 13,14,23,24,135,145,235,245 \}).$$
Then $\H_1$ and $\H_2$ are both matroids but $\H_1 \cup \H_2$ is not a TBRSC.
\ee

Indeed, it is easy to check that $\H_1$ and $\H_2$ are matroids. We may write $\H_1 \cup \H_2 = (V,H)$ with 
$$H = P_{\leq 2}(V) \cup \{ 135,145,235,245 \}.$$
Let $T \in T(H)$.

If $13 \subseteq T$, then $123 \notin H$ yields $2 \in T$, and $125 \notin H$ yields $5 \in T$.

If $15 \subseteq T$, then $125 \notin H$ yields $2 \in T$, and $123 \notin H$ yields $3 \in T$.

If $35 \subseteq T$, then $345 \notin H$ yields $4 \in T$, and $134 \notin H$ yields $1 \in T$.

It follows that $135 \notin (J(T(H)))_3$ and so $\H_1 \cup \H_2$ is not a TBRSC by Theorem \ref{eqtr}.

\bt
\label{pavun}
Let $d \geq 1$ and let $(V,H), (V,H') \in {\rm TBPav}(d)$. Then $(V,H\cup H') \in {\rm TBPav}(d)$.
\et

\proof
Let 
$$R = \{ T \cap T' \mid T \in T(H),\, T' \in T(H')  \}.$$
In view of Lemma \ref{propt}(i), $R$ is a Moore family and so $\J(R) = (V,J(R))$ is a BRSC. We claim that 
\beq
\label{pavun1}
H \cup H' = (J(R))_{d+1}.
\eeq

Let $X \in H$. By Theorem \ref{eqtr}, there exists a chain
\beq
\label{pavun2}
T_0 \subset T_1 \subset \ldots \subset T_n
\eeq
in $T(H)$ and an enumeration $x_1,\ldots,x_n$ of the elements of $X$ such that $x_i \in T_i \setminus T_{i-1}$ for every $i$. Since $V \in T(H')$, then (\ref{pavun2}) is also a chain in $R$, hence $X \in J(R)$. But $|X| \leq d+1$, thus $H \subseteq (J(R))_{d+1}$ and also $H' \subseteq (J(R))_{d+1}$ by symmetry.

Conversely, let $X \in (J(R))_{d+1}$. Since $H,H' \in \pav(d)$, we may assume that $|X| = d+1$. 
Then there exists a chain
\beq
\label{pavun3}
R_0 \subset R_1 \subset \ldots \subset R_{d+1}
\eeq
in $J(R)$ and an enumeration $x_1,\ldots,x_{d+1}$ of the elements of $X$ such that $x_i \in R_i \setminus R_{i-1}$ for every $i$. 

Write $R_d = T \cap T'$ with $T \in T(H)$ and $T' \in T(H')$. Since $x_{d+1} \notin R_d$, we may assume that $x_{d+1} \notin T$. Since $H \in \tbpav(d)$ yields $P_{\leq d-1}(V) \subseteq T(H)$, we have $P_{\leq d-1}(V)$ and so
$$\emptyset \subset x_1 \subset x_1x_2 \subset \ldots \subset x_1\ldots x_{d-1} \subset T \subset V$$
is a chain in $T(H)$ having $X$ as a transversal of the successive differences. Thus $X \in H$ by Theorem \ref{eqtr} and so (\ref{pavun}) holds. 

Therefore $(V,H\cup H') = (\J(R))_{d+1} \in {\rm TBPav}(d)$.
\qed

The next example shows that we cannot replace $\tbpav(d)$  by $\bpav(d)$ in Theorem \ref{pavun}.

\be
\label{ncu}
Let $V = 123456$, $H = P_{\leq 2}(V) \cup \{ 123, 124, 125, 126\}$ and 
$$H' = P_{\leq 2}(V) \cup 
\{ X \in P_3(V) \; {\big{\lvert}} \; |X \cap 56| = 1\}.$$
Then $(V,H),(V,H') \in {\rm BPav}(d)$ but $(V,H\cup H') \notin {\rm BPav}(d)$.
\ee

Indeed, it is easy to check that
$$\flatx(V,H) = P_{\leq 1}(V) \cup \{ 12, V\}, \quad
\flatx(V,H') = P_{\leq 1}(V) \cup \{ 1234, V\},$$
and it follows easily that 
$(V,H),(V,H') \in {\rm BPav}(d)$. We have seen in Example \ref{btbtwo} that $(V,H\cup H') \notin {\rm BPav}(d)$.

\medskip

Let $V$ be a finite nonempty set and let $L \subseteq V$ be such that $2 \leq d \leq |L| < |V|$. We write
$$\B_d(V,L) = (V,B_d(V,L)) = \J(P_{\leq d-1}(V) \cup \{ V,L\}).$$

\bl
\label{bvl}
Let $V$ be a finite nonempty set and let $L \subseteq V$ be such that $2 \leq d \leq |L| < |V|$. Then $\B_d(V,L) \in {\rm BPav}(d)$.
\el

\proof
We know that the operator $\J$ yields a BRSC. Considering the chains of the form
$$\emptyset \subset a_1 \subset a_1a_2 \subset \ldots a_1\ldots a_{d-1} \subset L \subset V$$
for $a_1,\ldots a_{d-1} \in V$ distinct, we confirm that $P_{\leq d}(V) \subseteq B_d(V,L)$. Since these are the unique maximal chains of maximal length in $P_{\leq d-1}(V) \cup \{ V,L\}$, it follows that $\B_d(V,L) \in \bpav(d)$.
\qed

We can now prove the following result.

\bt
\label{ubvl}
Let $d \geq 2$ and $\H\; = (V,H) \in {\rm Pav}(d)$. Then the following conditions are equivalent:
\bi
\item[(i)] $\H\in {\rm TBPav}(d)$;
\item[(ii)] $\H = \cup \{ \B_d(V,L) \mid L \in \L \}$ for some nonempty $\L \subseteq P_{\geq d} \setminus \{ V\}$.
\ei
\et

\proof
(i) $\Rw$ (ii). Let $\L = (P_{\geq d} \setminus \{ V\}) \cap T(H)$. Since $\dim\H = d$, we have $\L \neq \emptyset$.

Let $X \in H$ be a facet. By Theorem \ref{eqtr}, there exists a chain
$$T_0 \subset T_1 \subset \ldots \subset T_{d+1}$$
in $T(H)$ and an enumeration $a_1,\ldots a_{d+1}$ of the elements of $X$ so that $a_i \in T_i \setminus T_{i-1}$ for $i = 1,\ldots,d+1\}$. Now $a_1\ldots a_{i} \in \flatx\H \subseteq T(H)$ for $i = 0,\ldots,d-1$, hence $X$ is a transversal of the chain
$$\emptyset \subset a_1 \subset a_1a_2 \subset \ldots a_1\ldots a_{d-1} \subset T_d\subset V$$
and so $X \in B_d(V,T_d)$. Since $T_d \in \L$, we get $H \subseteq \cup \{ \B_d(V,L) \mid L \in \L \}$.

The opposite inclusion is immediate.

(ii) $\Rw$ (i). By Lemma \ref{bvl} and Theorem \ref{pavun}.
\qed

\subsection{The six point case}

We identify next, up to isomorphism, all the complexes with 6 points in $\tbpav(2) \setminus \bpav(2)$. We fix $V = \{ 1,\ldots,6\}$ as the set of points and we consider $\H\, = (V,H) \in \tbpav(2) \setminus \bpav(2)$. Then there exists some BRSC $\H' = (V,H')$ such that $\H = \H'_3$. Given $X \subseteq V$, let $\oo{X}$ (respectively $\wh{X}$) the closure of $X$ in $\flatx\H'$ (respectively $\flatx\H$). 

Since $\H\, \notin \bpav(2)$ and $P_{\leq 1}(V) \subseteq \flatx\H$, there exists some $X \in P_3(V) \cap H$ such that 
\begin{equation}
\label{six2}
X \subseteq \wh{X \setminus \{x\}} \mbox{ for every }x \in X.
\end{equation} 
Without loss of generality, we may assume that $X = 345$. On the other hand, since $\H'\, \in \bpav(2)$ and $345 \in H'$, there exists some $x \in 345$ such that $x \notin \oo{345 \setminus \{x\}}$. We may assume that $x = 5$. We claim that
\begin{equation}
\label{six1}
|\oo{34}| = 4.
\end{equation}
Indeed, we know already that $5 \notin \oo{34}$. Suppose that $\oo{34} = 34$. Then $34y \in H'$ (and therefore $34y \in H$) for every $y \in 1256$, yielding $\wh{34} = 34$, contradicting (\ref{six2}). Without loss of generality, we may assume that $34y \notin H'$ for some $y \in 126$, say $y = 1$. Hence $134 \subseteq \oo{34}$. Suppose that $\oo{34} = 134$. Since $134 \notin H$, it is easy to see that this implies $\wh{34} = 134$, contradicting (\ref{six2}). Thus $|\oo{34}| \geq 4$, and we may assume without loss of generality that $1234 \subseteq \oo{34}$. 

Suppose that $1234 \subset \oo{34}$. Since $5 \notin \oo{34}$, we get $\oo{34} = 12346$. It follows that $45z \in H'$ for every $z \in 1236$, hence $\wh{45} = 45$, contradicting (\ref{six2}). Therefore $\oo{34} = 1234$. It follows that $ab5, ab6 \in H$ for all $a,b \in 1234$ distinct. Since $\oo{134} = 1234$, it follows that $\{ 123, 124, 234\} \not\subseteq H$. Together with $134 \notin H$, this implies that the restriction $\H'' = (1234,H'')$ of $\H'$ to $1234$ misses at least two triangles. 

On the one hand, $134 \notin H$ and $\{ 123, 124, 234\} \not\subseteq H$ yield $1234 \subseteq \wh{34}$. 
On the other hand, it follows from (\ref{six2}) that $5 \in \wh{34}$, hence $12345 \subseteq \wh{34}$. Since $ab5, ab6 \in H'$ for all $a,b \in 1234$ distinct, then $1234 \notin \flatx\H$ implies that $1234 \setminus \{c\} \in H$ for some $c \in 1234$. Therefore $\H''$ has exactly one or two triangles. Since $\H''$ is a restriction of the BRSC $\H'$, it follows from \cite[Proposition 8.3.1]{RSm} that $\H''$ is a BRSC. On the other hand, it follows from \cite[Example 5.2.11]{RSm} that a paving BRSC with 4 points cannot have exactly one triangle, hence $\H''$ has exactly two triangles, whose intersection has two points, say $de$. 

Together with $1234 \in \flatx\H'$, this implies that $de \in \flatx\H'$. Since $134 \notin H$, we have $de \in \{ 12,23,24\}$. Since we have not distinguished 3 from 4 so far, we may assume that $de \in \{ 12,23\}$. 

In any case, having $1234 \in \flatx\H'$ determines that $ab5, ab6 \in H$ for all $a,b \in 1234$ distinct (16 elements), and $de \in \flatx\H'$ determines which two elements among the four elements of $P_3(1234)$ belong to $H$. Thus we only need to discuss what happens with $156,256,356,456$. If $356 \in H$, then $35 \in \flatx\H'$ (in view of $1234 \in \flatx\H'$), implying $\wh{35} = 35$ and contradicting (\ref{six2}). Therefore $356 \notin H$. Similarly, $456 \notin H$. It follows that we reduced the discussion to determine whether or not  $156, 256 \in H$, for each choice of $de \in \{ 12,23\}$.

We use now a simplification of the notation $B_2(V,L)$ introduced in Section \ref{sumc}. Given $L \subseteq V$ with $1 < |L| < |V|$, let $\B_2(L) = (V,B_2(L))$ be defined by
$$B_2(L) = P_{\leq 2}(V) \cup \{ X \in P_3(V) : |X \cap L| = 2 \}.$$

If we omit both $156, 256$ from $H$, we get the two cases
\bi
\item[(1)] $H = B_2(1234) \cup B_2(12)$, 
\item[(1')] $H = B_2(1234) \cup B_2(23)$,
\ei
which are clearly isomorphic.

Now adding $156$ (respectively $256$) is the only consequence of adding $15$ (respectively $25$) as a line, and these additions do not interfere with each other. We are then bound to consider the cases: 
\bi
\item[(2)] $H = B_2(1234) \cup B_2(12) \cup B_2(15)$;
\item[(2')] $H = B_2(1234) \cup B_2(12) \cup B_2(25)$;
\item[(3)] $H = B_2(1234) \cup B_2(12) \cup B_2(15) \cup B_2(25)$;
\item[(4)] $H = B_2(1234) \cup B_2(23) \cup B_2(15)$;
\item[(2'')] $H = B_2(1234) \cup B_2(23) \cup B_2(25)$;
\item[(5)] $H = B_2(1234) \cup B_2(23) \cup B_2(15) \cup B_2(25)$. 
\ei

The cases (2), (2') and (2'') are clearly isomorphic. Applying the transposition $(31)$ to $12345$ in the cases (4) and (5), we have reduced our discussion to the cases
\bi
\item[(1)] $H = B_2(1234) \cup B_2(12)$,;
\item[(2)] $H = B_2(1234) \cup B_2(12) \cup B_2(15)$;
\item[(3)] $H = B_2(1234) \cup B_2(12) \cup B_2(15) \cup B_2(25)$;
\item[(4)] $H = B_2(1234) \cup B_2(12) \cup B_2(35)$;
\item[(5)] $H = B_2(1234) \cup B_2(12) \cup B_2(25) \cup B_2(35)$. 
\ei

We list below the triangles missing in each of the cases:
\bi
\item[(1)] 134, 234, 156, 256, 356, 456;
\item[(2)] 134, 234, 256, 356, 456;
\item[(3)] 134, 234, 356, 456;
\item[(4)] 134, 234, 156, 256, 456;
\item[(5)] 134, 234, 156, 456.
\ei

Out of cardinality arguments, we only have to distinguish (2) from (4) and (3) from (5).
Now 1 appears only once in (2), and all points appear more often in (4); 1 and 2 appear only once in (3), but only 2 has a single occurrence in (5). Therefore these complexes (1) -- (5) are nonisomorphic.

By construction, any one of these 5 complexes is in $\tbpav(2)$. We confirm now that neither of them is a BRSC. For the first three cases, we take $345 \in H$.

\bi
\item[(1)] $134 \notin H$, hence $1 \in \wh{34}$; $234 \notin H$, hence $2 \in \wh{34}$; $\wh{34}$ contains the facet $123$, hence $\wh{34} = V$. $356 \notin H$, hence $6 \in \wh{35}$; $456 \notin H$, hence $4 \in \wh{35}$. Similarly, $\wh{45} = V$.
\item[(2)] Same argument as in (1).
\item[(3)] Same argument as in (1).
\ei

For the remaining two cases, we take $145 \in H$.
\bi
\item[(4)] $134 \notin H$, hence $3 \in \wh{14}$; $234 \notin H$, hence $2 \in \wh{14}$; $\wh{14}$ contains the facet $123$, hence $\wh{14} = V$. $156 \notin H$, hence $6 \in \wh{15}$; $456 \notin H$, hence $4 \in \wh{15}$. Similarly, $\wh{45} = V$.
\item[(5)] Same argument as in (4).
\ei

We have therefore proved:

\bp
\label{six}
Up to isomorphism, the complexes with 6 points in $\tbpav(2) \setminus \bpav(2)$ are of the form $(123456,H)$ for:
\bi
\item[(1)] $H = B_2(1234) \cup B_2(12)$;
\item[(2)] $H = B_2(1234) \cup B_2(12) \cup B_2(15)$;
\item[(3)] $H = B_2(1234) \cup B_2(12) \cup B_2(15) \cup B_2(25)$;
\item[(4)] $H = B_2(1234) \cup B_2(12) \cup B_2(35)$;
\item[(5)] $H = B_2(1234) \cup B_2(12) \cup B_2(25) \cup B_2(35)$. 
\ei
Moreover, all the above 5 cases are nonisomorphic.
\ep

\brem
\label{allfive}
We can build the following diagram
$$
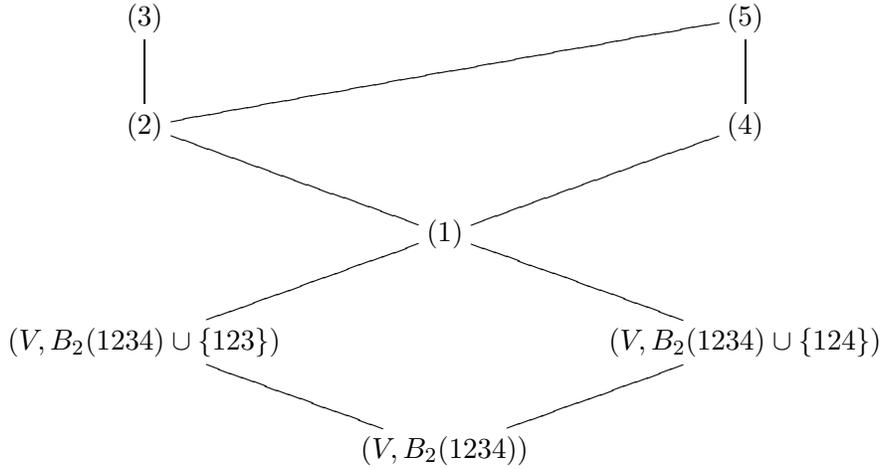
\xymatrix{
(3) \ar@{-}[d]&&(5) \ar@{-}[d] \ar@{-}[dll]\\
(2)\ar@{-}[dr]&&(4) \ar@{-}[dl]\\
&(1) \ar@{-}[dr] \ar@{-}[dl] &\\
(V,B_2(1234) \cup \{123\}) \ar@{-}[dr]&&(V,B_2(1234) \cup \{124\}) \ar@{-}[dl]\\
&(V,B_2(1234))&
}$$
The missing triangles in the three lowest elements are respectively
$$\xymatrix{
{124, 134, 234, 156, 256, 356, 456}&&{123, 134, 234, 156, 256, 356, 456}\\
&{123, 124, 134, 234, 156, 256, 356, 456}&
}$$
hence all the edges correspond to covering relations (recall the previous enumeration of the missing triangles for (1)--(5)).
\erem

We note that:
\bi
\item
$(V,B_2(1234)) \in {\rm BPav}(2)$ by Proposition \ref{bdlbr}.
\item
$(V,B_2(1234) \cup \{123\}) \notin {\rm TBPav}{2}$. Indeed, suppose that there exist $T \in T(B_2(1234) \cup \{123\})$ such that $|T \cap 123| = 2$. Since $124, 134, 234 \notin H$, we successively get $4 \in T$ and $1234 \subseteq T$, a contradiction. In view of Theorem \ref{eqtr}, this implies $(V,B_2(1234) \cup \{123\}) \notin {\rm TBPav}{2}$. 
\item
$(V,B_2(1234) \cup \{124\}) \notin {\rm TBPav}{2}$. Similar to the preceding case.
\item
No simplicial complex isomorphic to (4) embeds in (3). To prove this, recall the missing triangles in (3) and (4). We can check that $3x$ is contained in a missing triangle of (3) for every $x \neq 3$. On the other hand, $4y$ is contained in a missing triangle of (3) for every $x \neq 4$. Suppose that $\varphi \in S_6$ is such that the isomorphic image of (3) through $\varphi$ (call it (3'')) has (4) as subcomplex. Then the missing triangles of (3'') are a proper subset of the missing triangles of (4). Hence $(3\p)x$ is contained in a missing triangle of (4) for every $x \neq 3\p$, and $(4\p)x$ is contained in a missing triangle of (4) for every $x \neq 4\p$. However, only $4$ satisfies this property, yielding $3\p = 4 = 4\p$, a contradiction.
\ei

\subsection{On $\tbpav(d) \setminus \bpav(d)$}

Fix $d \geq 2$. 
Given $L \in P_{\geq d}(V) \setminus \{ V \}$, the complex $\B_d(V,L) = (V,B_d(V,L))$ can be described by
$$B_d(V,L) = P_{\leq d}(V) \cup \{ X \in P_{d+1}(V) \; {\big{\lvert}} \; |X \cap L| = d\}.$$

\bl
\label{flatsbdl}
Let $d \geq 2$ and $L \in P_{\geq d}(V) \setminus \{ V \}$. Then 
$$\flatx(\B_d(V,L)) = \left\{
\begin{array}{ll}
P_{\leq d-1}(V) \cup \{L,V\}&\mbox{ if }|L| < |V|-1\\
P_{\leq d-1}(V) \cup \{L,V\} \cup (P_d(V) \setminus P_d(L)) &\mbox{ if }|L| = |V|-1
\end{array}
\right.$$
\el

\proof
Write $H = B_d(V,L)$. 
Assume first that $|L| < |V|-1$. Let $F \in \flatx(\B_d(V,L))$ and assume that $F \notin P_{\leq d-1}(V) \cup \{V\}$. Then $d \leq |F| < |V|$. Suppose that $F\not\subseteq L$. Then there exists some $a \in F \setminus L$. Fix $b \in V \setminus (L \cup \{ a \})$. 

Suppose that $b \in F$. 
Since $d \leq |F| < |V|$, we can choose some $X \cup P_{d-2}(F \setminus \{ a,b \})$ and $c \in V \setminus F$. Since $X \cup \{ a,b\} \in P_d(F) \subseteq H$ and $c \notin F$, we get $X \cup \{ a,b,c\} \in H$, a contradiction since $a,b \notin L$. 

Hence we may assume that $b \notin F$. Choose $X \cup P_{d-1}(F \setminus \{ a,b \})$. Since $X \cup \{ a\} \in P_d(F) \subseteq H$ and $b \notin F$, we get $X \cup \{ a,b\} \in H$, a contradiction since $a,b \notin L$. 

Thus we may assume that $F\subseteq L$. Suppose that $F \subset L$. Take $X \in P_d(F)$ and $p \in L \setminus F$. Since $X \in P_d(F) \subseteq H$ and $p \notin F$, we get $X \cup \{ p\} \in H$, a contradiction since $X \cup \{ p\} \in P_{d+1}(L)$. 

Therefore $\flatx(\B_d(V,L)) \subseteq P_{\leq d-1}(V) \cup \{L,V\}$. The opposite inclusion is straightforward and we omit it.

Assume now that $|L| = |V|-1$, and denote by $a$ the single element of $V \setminus L$. Let $F \in \flatx(\B_d(V,L))$ and assume that $F \notin P_{\leq d-1}(V) \cup \{V\}$. Then $d \leq |F| < |V|$. 

Suppose that $F \subset L$. Take $X \in P_d(F) \subset H$ and $b \in L \setminus F$. Then $X\cup \{b\} \in H \cap P_{d+1}(L)$, a contradiction. Thus we may assume that $a \in F$.

Suppose that $|F| > d$. Since $a \in F$, it follows that $F$ contains a facet of $\B_d(V,L)$ and therefore $F = V$, a contradiction. It follows that $|F| = d$ and so $F \in P_d(V) \setminus P_d(L)$. 

Therefore $\flatx(\B_d(V,L)) \subseteq P_{\leq d-1}(V) \cup \{L,V\} \cup (P_d(V) \setminus P_d(L))$. The opposite inclusion is straightforward and we omit it.
\qed

\bp
\label{bdlbr}
Let $d \geq 2$ and $\emptyset \neq \L \subseteq P_{\geq d}(V) \setminus \{ V \}$ such that
\beq
\label{bdlbr1}
|L \cap L'| \leq d-1 \mbox{ for all distinct }L,L' \in \L.
\eeq
Then $(V,\displaystyle\bigcup_{L \in \L} B_d(V,L))$ is boolean representable.
\ep

\proof
Write $H = \displaystyle\cup_{L \in \L} B_d(V,L)$ and $\H\, = (V,H)$. Since $P_{\leq d}(V) \subseteq H$, we have $P_{\leq d-1}(V) \subseteq \flatx\H$. Let $K \in \L$ and suppose that $X \in H \cap 2^K$ and $p \in V \setminus K$. Since $P_{\leq d}(V) \subseteq H \subseteq P_{\leq d+1}(V)$, we may assume that $|X| = d$ or $d+1$. 

Suppose that $|X| = d+1$. Since $X \in H = \displaystyle\cup_{L \in \L} B_d(V,L)$, we have $X \in B_d(V,L)$ for some $L \in L$. Thus $|X \cap L| = d$ and so $|K \cap L| \geq d$. In view of (\ref{bdlbr1}), we get $K = L$, hence $X \subseteq L$, a contradiction since $|X| = d+1$ and $|X \cap L| = d$. Therefore
$|X| = d$, hence $X \cup \{ p \} \in B_d(V,K) \subseteq H$ and so $K \in \flatx\H$. 

Let $a_1, \ldots, a_{d-1} \in V$ be distinct. Then
\beq
\label{bdlbr2}
\emptyset \subset a_1 \subset a_1a_2 \subset \ldots \subset a_1\ldots a_{d-1} \subset V
\eeq
is a chain in $\flatx\H$. If $a_1, \ldots, a_{d-1} \in L \in \L$, then (\ref{bdlbr1}) can be refined to
\beq
\label{bdlbr3}
\emptyset \subset a_1 \subset a_1a_2 \subset \ldots \subset a_1\ldots a_{d-1} \subset L \subset V
\eeq
It is easy to check that every $X \in H$ is a partial transversal of the successive differences for some chain of type (\ref{bdlbr2}) or (\ref{bdlbr3}), hence 
$\H$ is boolean representable.
\qed

\bp
\label{makeu}
Let $d \geq 2$, $\H\, = (V,H) \in {\rm BPav}(d)$ and $\L\, = \{ L \in \flatx(V,H_i) \mid d \leq |L| < |V| \}.$ Then $\H\, = (V,\displaystyle\bigcup_{L \in \L} B_d(V,L))$.
\ep

\proof
Since $P_{\leq d-1}(V) \cup \{ V \} \subseteq \flatx\H$ and the maximum length of a chain in $\flatx\H$ is $d+1$, it follows easily that the maximal chains in $\flatx\H$ must be of the form (\ref{bdlbr2}) or (\ref{bdlbr3}), with $L \in \L$. Thus $H = \displaystyle\bigcup_{L \in \L} B_d(V,L)$.
\qed

\bp
\label{tpav}
The following conditions are equivalent for a simplicial complex $\H\, = (V,H)$ of dimension $d \geq 2$:
\bi
\item[(i)] $\H\, \in {\rm TBPav}(d)$;
\item[(ii)] $H = B_d(L_1) \cup \ldots \cup B_d(L_m)$ for some $m \geq 1$ and $L_1,\ldots,L_m \in P_{\geq d}(V) \setminus \{ V \}$.
\ei
\ep

\proof
(i) $\Rw$ (ii). Write $\H\, = \J_{d+1}$ for some BRSC $\J\, = (V,J)$. Since $P_{\leq d}(V) \subseteq H \subseteq J$, we have $P_{\leq d-1}(V) \subseteq \flatx\J$. 
Let $L_1,\ldots,L_m$ be an enumeration of the elements of
$$\L\, = \flatx\J \setminus (P_{\leq d-1}(V) \cup \{ V\}).$$
Note that $\L \neq \emptyset$, otherwise the maximal length of a chain in $\flatx\J$ is $d$ and $\J$ would have dimension $d$, contradicting $\H\, = \J_{d+1}$. We claim that
\beq
\label{tpav1}
H = B_d(L_1) \cup \ldots \cup B_d(L_m).
\eeq

Let $X \in H$. Since $P_{\leq d}(V) \subseteq B_d(L_i)$ by definition, we may assume that $|X| = d+1$. Since $H \subseteq J$, there exists a chain
\beq
\label{tpav2}
F_0 \subset F_1 \subset \ldots \subset F_{d+1}
\eeq
in $\flatx\J$  and an enumeration $x_1,\ldots,x_{d+1}$ of the elements of $X$ such that $x_i \in F_i \setminus F_{i-1}$ for $i = 1,\ldots, d+1$. But then we may replace (\ref{tpav2}) by the chain
$$\emptyset \subset x_1 \subset  \ldots \subset x_1\ldots x_{d-1} \subset F_{d} \subset V.$$
Since $x_1\ldots x_d \subseteq F_d \subset V$, we have $F_d = L_i$ for some $i \in \{ 1,\ldots,m\}$. Thus $X \in B_d(L_i)$ and $H \subseteq B_d(L_1) \cup \ldots \cup B_d(L_m)$. 

Conversely, assume that $X \in B_d(L_i)$. Then $X$ is a partial transversal of the successive differences for some chain of type (\ref{bdlbr2}) or (\ref{bdlbr3}), which is in any case a chain in $\flatx\J$. Thus $X \in J \cap P_{\leq d+1}(V) = H$ and (\ref{tpav1}) holds.

(ii) $\Rw$ (i). Let $R$ denote the Moore family generated by
$$P_{\leq d-1}(V) \cup \{ L_1,\ldots,L_m\}$$
and write $J = J(R)$. Then $\J\, = (V,J)$ is a BRSC. We prove that $\H\, = \J_{d+1}$.

Let $X \in H$. If $|X| \leq d$, say $X = x_1\ldots x_k$, then $X$ is a transversal of the successive differences for some chain in $R \subseteq \flatx\J$, namely
$$\emptyset \subset x_1 \subset  \ldots \subset x_1\ldots x_{i-1} \subset V.$$
Hence $X \in J$.

Thus we may assume that $|X| = d+1$. Since $X \in B_d(L_i)$ for some $i \in \{ 1,\ldots,m\}$, we may assume that $X = x_1\ldots x_dy$ with $x_1,\ldots, x_d \in L_i$ and $y \in V \setminus L_i$. Thus $X$ is a transversal of the successive differences for some chain in $R \subseteq \flatx\J$, namely
$$\emptyset \subset x_1 \subset  \ldots \subset x_1\ldots x_{d-1} \subset L_i \subset V.$$
Hence $X \in J$ and so $H \subseteq J \cap P_{\leq d+1}(V)$.

Conversely, let $X \in J \cap P_{\leq d+1}(V)$. Since $H \supseteq B_d(L_1) \supset P_{d}(V)$, we may assume that $|X| = d+1$. Then $X$ is a transversal of the successive differences for some chain 
\beq
\label{tpav3}
R_0 \subset R_1 \subset \ldots \subset R_{d+1}
\eeq
in $R$, so there exists an enumeration $x_1,\ldots,x_{d+1}$ of the elements of $X$ such that $x_i \in R_i \setminus R_{i-1}$ for $i = 1,\ldots,d+1$. 
But then we can replace (\ref{tpav3}) by
$$\emptyset \subset x_1 \subset  \ldots \subset x_1\ldots x_{d-1} \subset R_d \subset V,$$
which is also a chain in $R$. Since $R_d \in R$, then $R_d$ is necessarily an intersection of some elements of $R$. Since $d \leq |R_d| < |V|$, it follows that $R_d = L_{j_1} \cap \ldots \cap L_{j_k}$ for some $j_1, \ldots,j_k \in \{ 1,\ldots, m\}$. Since $x_{d+1} \notin R_d$, then $x_{d+1} \notin L_{j_s}$ for some $s \in \{ 1,\ldots,k\}$. Thus $X$ is also a transversal of the successive differences for the chain
$$\emptyset \subset x_1 \subset  \ldots \subset x_1\ldots x_{d-1} \subset L_{j_s} \subset V,$$
yielding $X \in B_d(L_{j_s}) \subseteq H$. Therefore $J \cap P_{\leq d+1}(V) = H$ and so $\H\, = \J_{d+1}$, implying $\H \in \tbpav(d)$.
\qed

\bl
\label{bigshort}
Let $d \geq 2$ and let $V$ be a finite set with $|V| \geq d+1$. For every $a \in V$, we have
\beq
\label{bigshort1}
B_d(V, V\setminus \{ a\}) =  \displaystyle\bigcup_{L \in \L_a} B_d(V,L),
\eeq
where $\L_a = \{ L \in P_d(V) \mid a \in L\}$.
\el

\proof
It suffices to show that both sides of (\ref{bigshort1}) contain the same $X \in P_{d+1}(V)$. So let $X \in P_{d+1}(V)$.

Suppose that $X \in B_d(V, V\setminus \{ a\})$. Then $a \in X$. Take $b \in X \setminus \{ a \}$. Then $X \setminus \{ b \} \in \L_a$ and so $$X \in B_d(V, X \setminus \{ b \}) \subseteq \displaystyle\bigcup_{L \in \L_a} B_d(V,L).$$

Conversely, suppose that $X \in B_d(V,L)$ with $L \in \L_a$. Since $|X| = d+1$ and $|L| = d$, we must have $X = L \cup \{ c \}$ for some $c \in V \setminus L$. Hence $a \in L \subset X$ yields $|X \cap (V\setminus \{ a\})| = d$ and $X \in B_d(V, V\setminus \{ a\})$. Therefore (\ref{bigshort1}) holds as required.
\qed

\bp
\label{smallpav}
Let $d \geq 2$ and and let $V$ be a finite set with $|V| \geq d+1$. Let $\emptyset \neq \L \subseteq 2^V$ be such that $|L| \in \{ d,d+1, |V|-1\}$ for every $L \in L$. Then $(V,\displaystyle\bigcup_{L \in \L} B_d(V,L))$ is boolean representable.
\ep

\proof
In view of Lemma \ref{bigshort}, we may assume that $|L| \in \{ d,d+1 \}$ for every $L \in L$. Write $H = \displaystyle\bigcup_{L \in \L} B_d(V,L)$ and $\H\, = (V,H)$. Since $P_{\leq d}(V) \subseteq H$, we have $P_{\leq d-1}(V) \subseteq \flatx\H$, thus we only need to show that every $X \in H \cap P_{d+1}(V)$ is a transversal of the successive differences for some chain in $\flatx\H$. 

Let $X \in H \cap P_{d+1}(V)$. Then $X \in B_d(V,K)$ for some $K \in \L$, hence the elements of $X$ admit an enumeration $x_1,\ldots, x_{d+1}$ such that $X \cap K = x_1\ldots x_d$. If $|K| = d$, then it is easy to see that $K \in \flatx\H$,  and $X$ is a transversal of the successive differences for the chain 
\beq
\label{smallpav1}
\emptyset \subset x_1 \subset x_1x_2 \subset \ldots \subset x_1\ldots x_{d-1} \subset K \subset V
\eeq
in $\flatx\H$. 

Hence we may assume that $|K| = d+1$.
Suppose that $K \notin H$. We claim that $K \in \flatx\H$. 
Let $Y \in H \cap 2^{K}$ and $p \in V \setminus K$. Since $K \notin H$ , we may assume that $|Y| \leq d$. But then $Y \cup \{ p\} \in B_d(V,K) \subseteq H$. Therefore $K \in \flatx\H$ and so $X$ is a transversal of the successive differences for the chain 
(\ref{smallpav1}).

Thus we may assume that $K \in H$. We claim that $X \cap K \in \flatx\H$. Note that $X \cap K \in P_d(V) \subset H$. Let $p \in V \setminus (X \cap K)$. If $p \notin K$, then $(X \cap K) \cup \{ p \} \in B_d(K) \subseteq H$. If $p \in K$, then $(X \cap K) \cup \{ p \} = K \in H$, so $X \cap K \in \flatx\H$. It follows that
$X$ is a transversal of the successive differences for the chain 
$$\emptyset \subset x_1 \subset x_1x_2 \subset \ldots \subset x_1\ldots x_{d-1} \subset X\cap K \subset V$$
in $\flatx\H$. Therefore $\H$ is a BRSC.
\qed

It follows from Proposition \ref{bdlbr} that $(V,B_d(L))$ is a BRSC for every $L \in P_{\geq d}(V) \setminus \{ V\}$. In view of Proposition
 \ref{six}, it follows that there exist $(V,H_1), (V,H_2) \in \bpav(d)$ such that $(V,H_1\cup H_2)$ is not a BRSC.
 
Thus we define, for every $d \geq 2$ and every finite set $V$ with at least $d+2$ elements,
$$\Y(V) = \{ (V,H) \in \bpav(d) \mid (V,H)\mbox{
 has no restriction isomorphic to }U_{d,d+2} \}.$$
 
\bp
\label{uniony}
Let $d \geq 2$ and let $V$ be a finite set with at least $d+2$ elements. Then
$$(V,H_1), (V,H_2) \in \Y(V) \mbox{ implies } (V,H_1 \cup H_2) \in \Y(V).$$
\ep

\proof
Let $i \in \{ 1,2\}$ and $F \in \flatx(V,H_i)$. Suppose that $d+1 < |F| < |V|$. Since $P_d(V) \subseteq H_i$ and no restriction of $(V,H_i)$ to a 4-subset of $F$ is isomorphic to $U_{d,d+2}$, there exists some $X \in P_{d+1}(F) \cap H_i$. But then $F$ contains a facet of $(V,H_i)$ and so $F = V$, a contradiction. Therefore $\flatx(V,H_i) \subseteq P_{\leq d+1}(V) \cup \{ V \}$.

Now let 
$$\L_i = \{ F \in \flatx(V,H_i) \mid d \leq |F| < |V| \}.$$
By the preceding claim, we have
\beq
\label{uniony1}
\L_i = \{ F \in \flatx(V,H_i) \mid |F| \in \{ d,d+1\} \, \}.
\eeq
Since $(V,H_i) \in \bpav(d)$, we have $H_i = \displaystyle\bigcup_{L \in \L_i} B_d(V,L)$ by Proposition \ref{makeu}. Thus $H_1 \cup H_2 = \displaystyle\bigcup_{L \in \L_1 \cup \L_2} B_d(V,L)$ and so $(V,H_1 \cup H_2) \in \bpav(d)$ by (\ref{uniony1}) and Proposition \ref{smallpav}. 

Suppose that there exists some $W \in P_{d+2}(V)$ such that $P_{d+1}(W) \cap (H_1 \cup H_2) = \emptyset$. Then $P_{d+1}(W) \cap H_1 = \emptyset$, contradicting $(V,H_1) \in \Y(V)$. Therefore $(V,H_1 \cup H_2) \in \Y(V)$. 
\qed

Note that an arbitrary $\H\in {\rm TBPav}(d) \setminus {\rm BPav}(d)$ needs not having a restriction isomorphic to $U_{d,d+2}$. The BRSCs featuring Proposition \ref{six} constitute all counterexamples for $d = 2$.

We intend now to show that ${\rm TBPav}(d) \setminus {\rm BPav}(d)$ is in some sense finitely generated. We start with a couple of lemmas.

Let $\tbr$ (respectively $\tbp$) denote the class of all finite truncated boolean representable simplicial complexes (respectively finite paving truncated boolean representable simplicial complexes).

\bl
\label{pret}
The classes $\tbr$ and $\tbp$ are prevarieties of simplicial complexes.
\el

\proof
Let $\H\, = (V,H) \in \tbr$ and let $\emptyset \neq W \subseteq V$. Since $\H\, \in \tbr$, there exist a BRSC $\J\, = (V,J)$ and $m \geq 1$ such that $\H\, = \J_m$. We claim that
\beq
\label{pret1}
\H_W = (\J_W)_m.
\eeq
This is equivalent to the equality 
\beq
\label{pret2}
H \cap 2^W = (J \cap 2^W) \cap P_{\leq m}(W).
\eeq
Now $\H\, = \J_m$ yields $H = J \cap P_{\leq m}(V)$ and so
$$H \cap 2^W = (J \cap P_{\leq m}(V)) \cap 2^W = (J \cap 2^W) \cap P_{\leq m}(W).$$
Hence (\ref{pret2}) and consequently (\ref{pret1}) do hold.

Since BRSCs are closed under restriction by \cite[Proposition 8.3.1(i)]{RSm}, then $\J_W$ isd a BRSC and it follows from (\ref{pret1}) that $\H_W \in \tbr$. Thus $\tbr$ is closed under restriction. Since it is also closed under isomorphism, then $\tbr$ is a prevariety of simplicial complexes.

On the other hand, the class of all finite paving simplicial complexes is a prevariety in view of \cite[Proposition 8.3.1(ii)]{RSm}. Since the intersection of two prevarieties is obviously a prevariety, it follows that $\tbp$ is a prevariety itself.
\qed

Let $\H\, \in \tbpav(d) \setminus \bpav(d)$. By Lemma \ref{pret}, every restriction of $\H$ is in $\tbp$ (with possibly lower dimension). We say that $\H$ is
 {\em minimal} if every proper restriction of $\H$ is boolean representable.
 
\bl
\label{comi}
Let $d \geq 2$. Then the maximum number of vertices for a minimal $\H\, \in \tbpav(d) \setminus \bpav(d)$ is $(d+1)(d+2)$.
\el

\proof
Let $\H\, = (V,H) \in \tbpav(d) \setminus \bpav(d)$ be minimal. Hence $\H\, \notin \bpav(d)$ but every proper restriction of $\H$ is boolean representable. By \cite[Theorem 8.5.2(ii)]{RSm}, we get $|V| \leq (d+1)(d+2)$. 

Now we consider the {\em Swirl}, the simplicial complex defined in the proof of \cite[Theorem 8.5.2(ii)]{RSm}, where it is proved that every proper restriction of this complex is boolean representable, but the Swirl is not. The Swirl is defined as follows:

Let $A = \{ a_0,\ldots, a_d\}$ and $B_i =  \{ b_{i0},\ldots, b_{id}\}$
for $i = 0,\ldots,d$. Write also $A_i = A \setminus \{ a_i
\}$ and 
$$C_i = P_{d+1}(A_i \cup (B_i\setminus \{ b_{i0} \})) \cup \{ B_i \}.$$
We define $$V = A \cup \displaystyle\bigcup_{i=0}^d B_i,\quad
H = P_{\leq d+1}(V) \setminus \displaystyle\bigcup_{i=0}^d C_i.$$
It is easy to check that all the $X \in H \cap P_{d+1}(V)$ fall into four cases (not necessarily disjoint):
\bi
\item[(a)]  
there exist $b_{ij}, b_{k\ell} \in X$ with $i \neq k$;
\item[(b)] 
there exist $b_{i0}, a_j \in X$;
\item[(c)] 
there exist $b_{ij}, a_i \in X$ with $j > 0$;
\item[(d)] 
$X = a_0\ldots a_d$.
\ei
Define
$$\begin{array}{lll}
\L&=&\{ L \in P_d(V) \mid \mbox{ there exist some $b_{ij}, b_{k\ell} \in L$ with $i \neq k$} \}\\
&\cup&\{ L \in P_d(V) \mid \mbox{ there exist some $b_{i0}, a_j \in L$} \}\\
&\cup&\{ L \in P_d(V) \mid \mbox{ there exist some $b_{ij}, a_i \in L$ with $j > 0$} \}\\
&\cup&\{ A_i \cup B_i \mid i = 0,\ldots,d\}.
\end{array}$$
It is straightforward to check that $H = \displaystyle\bigcup_{L\in \L} B_d(L)$, hence $(H,V) \in \tbpav(d)$ by Proposition \ref{tpav}. Since
$|V| = (d+1)(d+2))$, we have found some minimal $\H\, \in \tbpav(d) \setminus \bpav(d)$ with $(d+1)(d+2)$ vertices as required.
\qed

Let $\V$ be a prevariety of simplicial complexes. We say that $\V$ is {\em finitely based} if there exists some $m \geq 1$ such that every simplicial complex not in $\V$ admits a restriction not in $\V$ with at most $m$ vertices.

Given a prevariety $\V$ of simplicial complexes and $d \in
\mathbb{N}$, we define the prevariety
$$\V_d = \{ \H\; \in \V \mid \dim\H\; \leq d \}.$$

Let $\cal{BP}$ denote the class of all finite paving boolean representable simplicial complexes. By \cite[Theorem 8.5.2]{RSm}, ${\cal{BP}}_d$ is finitely based for every $d \geq 1$. Since ${\cal{TBP}}_1 = {\cal{BP}}_1$ by Proposition \ref{tbone}, it follows that ${\cal{TBP}}_1$ is finitely based.

\bt
\label{nfb}
${\cal{TBP}}_2$ is not finitely based.
\et

\proof
It suffices to build arbitrary large simplicial complexes not in ${\cal{TBP}}_2$ with all proper restrictions in ${\cal{TBP}}_2$.

Let $n \geq 6$ and take as vertex set
$$V = \{ x_0, \ldots, x_n, y_0, \ldots, y_6, z_0, \ldots, z_6 \},$$
where we identify
$$x_0 = y_0 = z_6, \quad
x_1 = z_0 = y_6, \quad 
y_1 = z_1 = x_n\}.$$
Let
$$S = \{ x_ix_{i+1}x_{i+2} \mid i = 0,\ldots, n-2\} \cup 
\{ y_iy_{i+1}y_{i+2} \mid 0 \leq i \leq 4 \} \cup \{ z_iz_{i+1}z_{i+2} \mid 0 \leq i \leq 4 \},$$
$H = P_{\leq 3}(V) \setminus S$ 
and $\H\, = (V,H)$. 

Clearly, $\H\, \in \pav(2)$. Hence $P_{\leq 1}(V) \subseteq \flatx\H \subseteq T(H)$. Thus 
\begin{eqnarray}\nonumber
\label{nfb1}
\H\, \in \tbp_2 \mbox{ if and only if, for every $X \in H \cap P_3(V)$,}\hspace{3cm}\\
\hspace{2cm}\mbox{ there exists some $T \in T(H)$ such that }|X \cap T| = 2.
\end{eqnarray}

First, we note that
\beq
\label{nfb2}
\mbox{for every $X \in (H \cap P_3(V)) \setminus \{ x_0x_1y_1\}$, there exists some $F \in \flatx\H$ such that }|X \cap F| = 2.
\eeq
Indeed, such an $X$ contains necessarily some element of $V \setminus \{ x_0x_1y_1\}$. Without loss of generality, we may assume that this element is among $x_2,\ldots,x_{n-1}$ (the other cases follow by symmetry). 

Suppose that $X \subset x_0\ldots x_n$, say $X = x_ix_jx_k$ with $i < j < k$. Since $X \in H$, then $i,j,k$ are not consecutive integers. If $k < n$, then $k > 1$ and $k-i > 2$, hence $x_ix_k \in \flatx\H$ and we are done. Thus we may assume that $k = n$. If $i > 1$, then $k-i > 2$, hence $x_ix_k \in \flatx\H$ and we are done. Thus we may assume that $i \leq 1$. Since $k = n$, this implies $2 \leq j \leq n-1$. Since $n \geq 6$, we get either $k-j > 2$ (yielding $x_jx_k \in \flatx\H$) or $j-i > 2$ (yielding $x_ix_j \in \flatx\H$). 

Hence we may assume that at least one of the other elements of $X$ (say $a$) is not of the form $x_j$. Let $i \in \{ 2,\ldots,n-1\}$ be such that $x_i \in X$. It is easy to check that $x_ia \in \flatx\H$. Therefore (\ref{nfb2}) holds.

Next we show that 
\beq
\label{nfb3}
\mbox{for every $T \in T(H)$, $|T \cap \{ x_0x_1y_1\}| \neq 2$.}
\eeq

Let $T \in T(H)$ and assume that $|T \cap \{ x_0x_1y_1\}| \geq 2$. Assume first that $x_0,x_1 \in T$. Since $x_ix_{i+1}x_{i+2} \notin H$ for $i = 0, \ldots, n-2$, we get successively $x_2 \in T,\ldots, x_n \in T$. Since $x_n = y_1$, we get $x_0x_1y_1 \subseteq T$.

Suppose now that $x_0,y_1 \in T$. Since $x_0 = y_0$, we use the same argument to deduce that $y_2 \in T,\ldots, y_6 \in T$. Since $y_6 = x_1$, we get $x_0x_1y_1 \subseteq T$.

Finally, suppose that $x_1,y_1 \in T$. Since $x_1 = z_0$ and $y_1 = z_1$, we use the same argument to deduce that $z_2 \in T,\ldots, z_6 \in T$. Since $z_6 = x_0$, we get $x_0x_1y_1 \subseteq T$ and (\ref{nfb3}) is proved.

In view of (\ref{nfb1}), it follows from (\ref{nfb3}) that $\H\, \notin \tbp_2$.

Fix now $v \in V$ and write $W = V\setminus \{ v\}$. We must show that $\H|_{W} \in \tbp_2$ (since $\tbp_2$ is closed under restrictions, this implies that $\H|_{W'} \in \tbp_2$ for any $W' \subset V$). 

Since $\H\, \in \pav(2)$, we only need to show that the righthand side of (\ref{nfb1}) holds when we replace $\H$ by $\H|_W$. Let $X \in H \cap P_3(W)$. Suppose first that $X \neq x_0x_1y_1$. By (\ref{nfb2}),
there exists some $F \in \flatx\H$ such that $|X \cap F| = 2$. By \cite[Proposition 8.3.3(i)]{RSm}, $F \cap W \in \flatx(\H|_W)$.  Since $|X \cap (F \cap W)| = 2$, the desired condition is satisfied if $X \neq x_0x_1y_1$. 

Thus we may assume that $X = x_0x_1y_1$. It follows that either $v = x_i$ with $2 \leq i \leq n-1$ or $v = y_j$ or $v = z_j$ with $2 \leq j \leq 5$.

Suppose that $v = x_i$. Let $T = x_0\ldots x_{i-1}$. It is immediate that $T \in T(H \cap 2^W)$ and $|T \cap x_0x_1y_1| = 2$. If $v = y_j$ (respectively $v = z_j$), we take $T = y_0\ldots y_{j-1}$ (respectively $T = z_0\ldots z_{j-1}$). Therefore, in view of (\ref{nfb1}), we get $\H|_{W} \in \tbp_2$ as required.
\qed

\section{The Pure Conjecture}

Let $\H\; = (V,H)$ be a simplicial complex of dimension $d$. We define $\pure(\H) = (V',H')$ by
$$V' = \cup(H \cap P_{d+1}(V)),\quad H' = \cup_ {X \in H \cap P_{d+1}(V)} 2^X.$$
It is immediate that $\pure(\H)$ is the largest pure subcomplex of $\H$.

This section is devoted to the following conjecture, which we call the {\em Pure Conjecture}:

\bj
\label{tpc}
Let $\H$ be a BRSC and let $k \geq 1$. Then {\rm pure}$(\H_k)$ is a BRSC.
\ej

We can disprove the conjecture for $\dim\H\; = 3$ and $k = 3$.

\be
\label{cepc}
Let $\H\; = (V,H)$ with $V = \cup_{i \in \Z_3} \{ i,i',i'' \}$,
$$Z = \cup_{i \in \Z_3} \{ i(i+1)(i+1)', i''(i+1)(i+1)' \}$$
and
$$\begin{array}{lll}
H&=&(P_{\leq 3}(V) \setminus Z) \cup \{ ii''(i+1)p \mid i \in \Z_3,\; p \in V \setminus ii''(i+1)(i+1)' \}\\
&\cup&\{ ii''(i+1)'p \mid i \in \Z_3,\; p \in V \setminus ii''(i+1)(i+1)' \}.
\end{array}$$
Then $\H$ is a BRSC but ${\rm pure}(\H_3)$ is not.
\ee

It is easy to check that $\H$ is indeed a simplicial complex. Clearly, $P_{\leq 1}(V) \subset \flatx\H$. If $X \in P_2(V)$ is not contained in any element of $Z$, then $\oo{X} = X$. Hence, if $abc \in H$ and $ab$ is not contained in any element of $Z$, then $abc$ is a transversal of the successive differences for the chain
$$\emptyset \subset a \subset ab \subset V$$
in $\flatx\H$. On the other hand, it is easy to check that the unique $X \in P_3(V) \cap H$ having all 2-subsets contained in elements of $Z$ is $123$ (see the picture below, where the yellow triangles are the elements of $Z$):

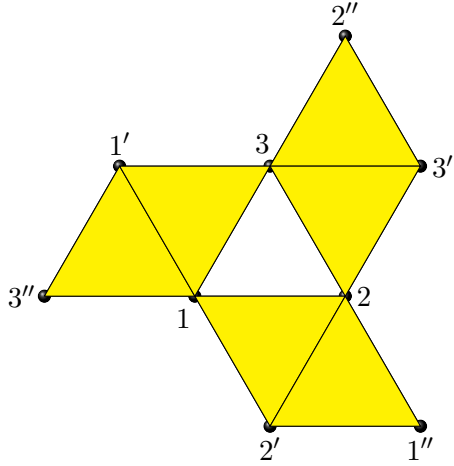
\begin{center}
\begin{tikzpicture}[domain=-4:4]
\draw (0,0) -- (4,0);
\draw (1,1.73) -- (5,1.73);
\draw (3,-1.73) -- (5,-1.73);
\draw (0,0) -- (1,1.73);
\draw (2,0) -- (4,3.46);
\draw (3,-1.73) -- (5,1.73);
\draw (1,1.73) -- (3,-1.73);
\draw (3,1.73) -- (5,-1.73);
\draw (5,1.73) -- (4,3.46);
\node at (1,2.03) {$1'$};
\node at (4,3.76) {$2''$};
\node at (2.9,2.03) {$3$};
\node at (5.3,1.73) {$3'$};
\node at (1.85,-0.3) {$1$};
\node at (-0.3,0) {$3''$};
\node at (4.25,0) {$2$};
\node at (3,-2.03) {$2'$};
\node at (5,-2.03) {$1''$};
\shade[shading=ball, ball color=black] (1,1.73) circle (.09);
\shade[shading=ball, ball color=black] (3,1.73) circle (.09);
\shade[shading=ball, ball color=black] (5,1.73) circle (.09);
\shade[shading=ball, ball color=black] (3,-1.73) circle (.09);
\shade[shading=ball, ball color=black] (5,-1.73) circle (.09);
\shade[shading=ball, ball color=black] (4,3.46) circle (.09);
\shade[shading=ball, ball color=black] (0,0) circle (.09);
\shade[shading=ball, ball color=black] (2,0) circle (.09);
\shade[shading=ball, ball color=black] (4,0) circle (.09);
\draw [fill = yellow] (0,0) -- (2,0) -- (1,1.73) -- (0,0);
\draw [fill = yellow] (3,1.73) -- (2,0) -- (1,1.73) -- (3,1.73);
\draw [fill = yellow] (3,1.73) -- (4,3.46) -- (5,1.73) -- (3,1.73);
\draw [fill = yellow] (4,0) -- (3,1.73) -- (5,1.73) -- (4,0);
\draw [fill = yellow] (2,0) -- (4,0) -- (3,-1.73) -- (2,0);
\draw [fill = yellow] (4,0) -- (3,-1.73) -- (5,-1.73) -- (4,0);
\end{tikzpicture}
\end{center}

Now we check that $ii''(i+1)(i+1)'  \in \flatx\H$ for every $i \in \Z_3$. It follows that
$123$ is a transversal of the successive differences for the chain
$$\emptyset \subset 1 \subset 11''22' \subset V$$
in $\flatx\H$. 

Finally, each facet of the form $ii''(i+1)p$ or $ii''(i+1)'p$ is a transversal of the successive differences for the chain
$$\emptyset \subset i \subset ii'' \subset ii''(i+1)(i+1)' \subset V$$
in $\flatx\H$. Since we have now checked all facets, it follows that $\H$ is a BRSC.

Let $\cl(X)$ denote the closure of $X \subseteq V$ in $\flatx\H_3$. For each $i \in \Z_3$, we have $i(i+1)(i+1)', i''(i+1)(i+1)' \notin H_3$, so we successively get $(i+1)' \in \cl(i(i+1))$ and $i'' \in \cl(i(i+1))$. Thus $\cl(i(i+1))$ contains $ii''(i+1) \in \fct\H_3$, yielding $\cl(i(i+1)) = V$. But then $i \in \cl(123 \setminus \{ i \})$ for every $i \in 123$. Since $123 \in H_3$, it follows from \cite[Corollary 5.2.7]{RSm} that $\H_3$ is not boolean representable.

We remark that $\H$ is not pure since it is straightforward to check that $1'2'2''$ is a facet. But $\H_3$ is pure because there are nine vertices and $|Z| = 6$.

\medskip

Another counterexample is given by the boolean module $\B^{(4)}$: a simplicial complex admitting a $4  \times (2^4-1)$ boolean matrix representation where all columns are distinct and nonzero.

\be
\label{bfour}
The boolean module $\B^{(4)}$ is pure and its truncation to rank 3 is a pure TBRSC which is not a BRSC.
\ee

Let $M$ be such a boolean matrix. Since the columns are all distinct and nonzero, every pair of distinct columns is independent. Now let $X$ be a set of independent columns with $|X| = 2$ or $3$. Let $I \subset 1234$ be such that the square matrix $M[I,X]$ is nonsingular. Let $j \in 1234 \setminus I$ and let $c$ be the column of $M$ having a 1 at row $j$ and 0 elsewhere. Then the permanent of $M[I \cup \{ j \},X \cup \{ c \}]$ equals the permanent of $M[I,X]$, hence $M[I \cup \{ j \},X \cup \{ c \}]$ is nonsingular. and so  $X \cup \{ c \}$ is independent. Thus $\B^{(4)}$ is pure.

Since $\B^{(4)}$ is by definition a BRSC, then $\B^{(4)}_3$ is a TBRSC. Let $\oo{X}$ denote the closure of $X$ in $\flatx \B^{(4)}_3$. Consider the columns of $M$ defined by
$$a = \begin{bmatrix}
1\\ 0\\ 0\\ 0
\end{bmatrix}, \quad
b = \begin{bmatrix}
1\\ 1\\ 1\\ 0
\end{bmatrix}, \quad
c = \begin{bmatrix}
1\\ 1\\ 0\\ 1
\end{bmatrix}.$$
The permanent of the matrix
$$M[134,abc] = \begin{bmatrix}
1 & 1 & 1\\ 0 & 1 & 0\\ 0 & 0 & 1
\end{bmatrix}$$
is 1, hence $abc$ is independent. Define
$$d = \begin{bmatrix}
0\\ 1\\ 1\\ 0
\end{bmatrix}, \quad
e = \begin{bmatrix}
1\\ 0\\ 1\\ 0
\end{bmatrix}, \quad
f = \begin{bmatrix}
0\\ 0\\ 1\\ 1
\end{bmatrix} \quad
g = \begin{bmatrix}
1\\ 0\\ 1\\ 1
\end{bmatrix}.$$
We have
$$M[1234,abd] = \begin{bmatrix}
1 & 1 & 0\\ 0 & 1 & 1\\ 0 & 1 & 1\\ 0 & 0 & 0
\end{bmatrix},\quad 
M[1234,bde] = \begin{bmatrix}
1 & 0 & 1\\ 1 & 1 & 0\\ 1 & 1 & 1\\ 0 & 0 & 0
\end{bmatrix},\quad 
M[1234,abe] = \begin{bmatrix}
1 & 1 & 1\\ 0 & 1 & 0\\ 0 & 1 & 1\\ 0 & 0 & 0
\end{bmatrix}.$$
Since no row of $M[1234,abd]$ has precisely two zeroes, $abd$ is dependent. The same occurs with $bde$. It is immediate that $M[123,abe]$ has permanent 1, hence $abe$ is independent. Thus we successively deduce $d \in \oo{ab}$, $e \in \oo{ab}$ and so $\oo{ab}$ contains the facet $abe$. Therefore $\oo{ab}$ is the full set of vertices. Similarly, so is $\oo{ac}$. 

Now
$$M[1234,bcf] = \begin{bmatrix}
1 & 1 & 0\\ 1 & 1 & 0\\ 1 & 0 & 1\\ 0 & 1 & 1
\end{bmatrix},\quad 
M[1234,bcg] = \begin{bmatrix}
1 & 1 & 1\\ 1 & 1 & 0\\ 1 & 0 & 1\\ 0 & 1 & 1
\end{bmatrix},\quad 
M[1234,bfg] = \begin{bmatrix}
1 & 0 & 1\\ 1 & 0 & 0\\ 1 & 1 & 1\\ 0 & 1 & 1
\end{bmatrix}.$$
Since no row of $M[1234,bcf]$ has precisely two zeroes, $bcf$ is dependent. The same occurs with $bcg$. It is immediate that $M[123,bfg]$ has permanent 1, hence $bfg$ is independent. Thus we successively deduce $f \in \oo{bc}$, $g \in \oo{bc}$ and so $\oo{bc}$ contains the facet $bfg$. Therefore $\oo{bc}$ is the full set of vertices. This proves that $\B^{(4)}_3$ is not a BRSC.

\medskip

However, the Pure Conjecture holds for particular cases as we shall see.

\bl
\label{wool}
Let $\H\; = (V,H)$ be a simplicial complex and let $I,J \in H$ be such that $I \subseteq \oo{J}$. Then there exists some $I' \in H$ such that $I \subseteq I'$ and $\oo{I'} = \oo{J}$.
\el

\proof
Let $I' \in H$ be maximal with respect to $I \subseteq I' \subseteq \oo{J}$. If $\oo{I'} \subset \oo{J}$, we can take $p \in \oo{J} \setminus \oo{I'}$ and get $I' \cup \{ p\} \in H \cap 2^{\oo{J}}$, contradicting the maximality of $I'$. Thus $\oo{I'} = \oo{J}$ and we are done.
\qed

It is well known that a simplicial complex $\H\; = (V,H)$ is a matroid if and only if
\beq
\label{silk}
\oo{a_1} \subset \oo{a_1a_2} \subset \ldots \subset \oo{a_1\ldots a_k}
\eeq
holds for every $a_1\ldots a_k \in H$ (where the enumeration is arbitrary).

We present next a characterization of matroids:

\bp
\label{mrf}
Let $\H\; = (V,H)$ be a simplicial complex. Then the following conditions are equivalent:
\bi
\item[(i)] $\H$ is a matroid;
\item[(ii)] for all $X,Y \in H$, $\oo{X} = \oo{Y}$ implies $|X| = |Y|$.
\ei
\ep

\proof
(i) $\Rw$ (ii). Suppose that $X,Y \in H$ are such that $\oo{X} = \oo{Y}$ and $|X| < |Y|$. By the exchange property, we have $X \cup \{ y\} \in H$ for some $y \in Y\setminus X$. Hence $\oo{X \cup \{ y\}} = \oo{Y} = \oo{X}$, contradicting (\ref{silk}).

(ii) $\Rw$ (i). Let $I,J \in H$ be such that $|I| = |J| +1$. Suppose that $J\cup \{ i \} \notin H$ for every $i \in I \setminus J$. Then $I \subseteq \oo{J}$ and so by Lemma \ref{wool} this contradicts condition (ii). Thus $\H$ satisfies the exchange property and is therefore a matroid.
\qed 

A simplicial complex $\H\; = (V,H)$ is said to be a {\em near-matroid} if 
$$\oo{X} = \oo{Y} \subset V \mbox{ implies }|X| = |Y|$$
for all $X,Y \in H$. In this case we can define a function $\rho: \flatx\H \setminus \{ V \} \to \mathbb{N}$ by
$$F\rho = |X|, \mbox{ where $X \in H$ is such that }\oo{X} = F.$$
Note that  such an $X$ exists by \cite[Proposition 4.2.4]{RSm}.

It follows from Proposition \ref{mrf} that every matroid is a near-matroid. The following result shows that the converse fails.

\bp
\label{bpnm}
Let $d \geq 0$ and $\H\,\in {\rm BPav}(d)$. Then $\H$ is a near-matroid.
\ep

\proof
Write $\H\, = (V,H)$ and suppose that $X,Y \in H$ are such that $\oo{X} = \oo{Y} \subset V$. By 
\cite[Proposition 4.2.3]{RSm}, $X$ and $Y$ are not facets. Suppose that $|X| < d$. Since $P_d(V) \subseteq H$, it follows that $\oo{X} = X$, so in this case we get indeed $Y = X$. Thus we may assume by symmetry that $|X|,|Y| \geq d$. Since $X$ and $Y$ are not facets, then $|X| = d = |Y|$ and so $\H$ is a near-matroid.
\qed

The following example shows that not every near-matroid is boolean representable, even in the paving case.

\be
\label{far}
Let $\H\, = (V,H)$ be the simplicial complex defined by $V = 1234$ and $H = P_{\leq 2}(V) \cup \{ 123\}$. Then $\H$ is a near-matroid which is not boolean representable.
\ee

Indeed, it is easy to check that $\flatx\H\, = P_{\leq 1}(V) \cup \{ V\}$ and $\H$ is a near-matroid. On the other hand, $\H$ is not boolean representable by \cite[Example 5.2.11]{RSm}.

\bl
\label{opre}
Let $\H\, = (V,H)$ be a 
near-matroid and let $F,F' \in {\rm Fl}\H$ be such that $F \subset F' \subset V$. Then $F\rho < F'\rho$.
\el

\proof
Suppose that $F\rho \geq F'\rho$. Then there exist $I,J \in H$ such that $F = \oo{I}$, $F' = \oo{J}$ and $|I| \geq |J|$. Hence $I \subseteq \oo{J}$ and so by Lemma \ref{wool} there exists some $I' \in H$ such that $I \subseteq I'$ and $\oo{I'} = \oo{J}$. But we have then $|I'| > |I| \geq |J|$, a contradiction since $\H$ is a near-matroid. Therefore $F\rho < F'\rho$.
\qed

\bl
\label{step}
Let $\H\, = (V,H)$ be a 
near-matroid and let $F,F' \in {\rm Fl}\H$ be such that $F \subset F' \subset V$. Let $a_1 \in F'\setminus F$ and $k = F'\rho -F\rho$. Then there exist $a_2,\ldots a_k \in V$ such that 
$$F \subset \oo{F\cup a_1} \subset \oo{F\cup a_1a_2} \subset \ldots \subset \oo{F\cup a_1\ldots a_k} = F'.$$
\el

\proof
Write $F = \oo{I}$ with $I \in H$. Since $a_1 \in F'\setminus F$, we have $I \cup a_1 \in H$. Thus 
$$F \subset \oo{I\cup a_1} = \oo{F\cup a_1} \subseteq F'.$$
Moreover, $$\oo{F\cup a_1}\rho = |I \cup a_1| = |I|+1 = F\rho+1.$$
If $\oo{F\cup a_1} = F'$, we can now iterate this argument to produce a chain
$$F \subset \oo{F\cup a_1} \subset \oo{F\cup a_1a_2} \subset \ldots \subset \oo{F\cup a_1\ldots a_s} = F'$$
for some $a_2,\ldots a_s \in V$ such that 
$\oo{F\cup a_1\ldots a_j}\rho = \oo{F\cup a_1\ldots a_{j-1}}\rho+1$ for $j = 1,\ldots,s$. Thus $s = F'\rho -F\rho = k$ and we are done.
\qed  

We can now prove the Pure Conjecture for boolean representable near-matroids.

\bt
\label{brnm}
Let $\H$ be a boolean representable near-matroid and let $k \geq 0$. Then:
\bi
\item[(i)] $\H_k$ is a BRSC;
\item[(ii)] ${\rm pure}(\H_k)$ is a BRSC.
\ei
\et

\proof
(i) Write $\H\, = (V,H)$ and define
$$\F_k = \{ F \in \flatx\H \mid F\rho < k\} \cup \{ V\}.$$
We claim that the matrix $M(\F_k)$ is a boolean representation of $\H_k$, i.e. $H_k = J(\F_k)$.

Let $X \in H_k$ and let $s = |X|$. Then there exists an enumeration $a_1,\ldots,a_s$ of the elements of $X$ such that 
$$\oo{a_1} \subset \oo{a_1a_2} \subset \ldots \subset \oo{a_1\ldots a_s}.$$
Hence $X$ is a transversal of the successive differences for
$$\emptyset \subset \oo{a_1} \subset \oo{a_1a_2} \subset \ldots \subset \oo{a_1\ldots a_{s-1}} \subset V,$$
which is a chain in $\F_k$. Thus $X \in J(\F_k)$.

Conversely, assume that $X \in J(\F_k)$. Since $\F_k \subseteq \flatx\H$ and $M(\flats\H)$ is a boolean representation of $\H$ by \cite[Theorem 5.2.5]{RSm}, it follows that $X \in H$.

Suppose that $|X| > k$. Since $X \in J(\F_k)$, there exist some $F \in F_k$ and $x \in X$ such that $F \cap X = X \setminus \{ x\}$. But $F = \oo{Y}$ for some $Y \in H_{k-1}$. Hence $X \setminus \{ x\} \subseteq \oo{Y}$ and by Lemma \ref{wool} there exists some $Z \in H$ such that $X \setminus \{ x\} \subseteq Z$ and $\oo{Z} = \oo{Y} = F$. But then $|Z| \geq |X \setminus \{ x\}| \geq k > |Y|$, a contradiction since $\H$ is a near-matroid. Thus $X \in H_k$ and so $H_k = J(\F_k)$ as claimed.

(ii) Let $\F'_k$ denote the set of all the flats of $\H$ occurring in chains of the form
$$F_0 \subset F_1 \subset \ldots \subset F_k$$
in $\F_k$.
Note that
\beq
\label{brnm1}
\F_k \mbox{ is closed under intersection}.
\eeq
Indeed, by \cite[Proposition 4.2.2(ii)]{RSm}, $\flatx\H$ is closed under intersection, and the bound for $\rho$ follows easily from Lemma \ref{opre}.

Next we show that
\beq
\label{brnm2}
\F'_k \mbox{ is closed under intersection}.
\eeq
Let $F,F' \in \F'_k$. Then $F \cap F' \in \F_k$ by (\ref{brnm1}). Since $F \in \F'_k$, there exists some $F'' \in \flatx\H$ such that $F \subseteq F''$ and $F''\rho = k-1$. Now we apply Lemma \ref{step} to both inclusions $\emptyset \subseteq F\cap F\subseteq F''$. This ensures that $F\cap F'$ will appear in some chain of flats of length $k$ in $\flatx\H$ of the form
$$\emptyset \subset \ldots \ldots F'' \subset V.$$
Since $F''\rho = k-1$, it follows from Lemma \ref{opre} that this is in fact a chain in $\F_k$ and therefore in $\F'_k$. Thus $F\cap F' \in \F'_k$ and so (\ref{brnm2}) holds.

We claim that the matrix $M(\F'_k)$ is a boolean representation of $\pure(\H_k)$.

Let $X \in H \cap P_k(V)$. Then there exists an enumeration $a_1,\ldots,a_k$ of the elements of $X$ such that 
$$\oo{a_1} \subset \oo{a_1a_2} \subset \ldots \subset \oo{a_1\ldots a_k}.$$
Hence $X$ is a transversal of the successive differences for
$$\emptyset \subset \oo{a_1} \subset \oo{a_1a_2} \subset \ldots \subset \oo{a_1\ldots a_{k-1}} \subset V,$$
which is a chain of length $k$ in $\F_k$. Thus $X \in J(\F'_k)$.

Conversely, assume that $X \in J(\F'_k)$. We may assume that $X$ is a facet of $\J(\F'_k)$. In view of (\ref{brnm2}), we may assume that $X$ is a transversal of the partial differences for some chain in $\F'_r$, say 
\beq
\label{brnm3}
F_0 \subset F_1 \subset \ldots \subset F_s.
\eeq
Thus there exists some enumeration $a_1,\ldots,a_s$ of the elements of $X$ such that $a_i \in F_i \setminus F_{i-1}$ for $i = 1,\ldots,s$.
Since $X$ is a facet, we must have $F_0 = \emptyset$ and $F_s = V$. Suppose that $F_{s-1}\rho = r < k-1$. Since $F_{s-1} \in \F'_k$, then it must occur in some chain of length $k$ in $\F'_k$, hence 
we have some chain 
$$F_{s-1} = F'_0 \subset F'_1 \subset \ldots \subset F'_t \subset F'_{t+1} = V$$ 
in $\F'_k$ for some $t\geq 1$. Since $a_s \in F'_{t+1} \setminus F'_0$, we have $a_s \in F'_{j} \setminus F'_{j-1}$ for some $j \in \{1,\ldots,t+1\}$, hence there exists some $Y \in J(\F'_k) \in P_{s+t}(V)$ containing (strictly) $X$, contradicting $X \in \fct(\J(\F'_k))$. Thus $F_{s-1}\rho = k-1$. 

Now $a_i \in F_i \setminus F_{i-1}$ for $i = 1,\ldots,s-1$ and so we can apply Lemma \ref{step} $s-1$ times to refine (\ref{brnm3}) to a chain of length $k$ in $\flatx\H$ of the form
$$F_0 \subset \oo{F_0 \cup a_1} \subseteq \ldots \subseteq F_1 \subset \oo{F_1 \cup a_2} \subseteq \ldots \subseteq F_{s-1} \subset F_s,$$
which admits a transversal of the successive differences containing $X$. Since $X \in \fct(\J(\F'_k))$, it follows that $s = k$ and so in view of Lemma \ref{opre} we have $X \in H_k\cap P_k(H)$, hence $X$ is a facet of $\pure(\H_k)$. Therefore $M(\F'_k)$ is a boolean representation of $\pure(\H_k)$ as claimed.
\qed

It remains an open problem whether or not $\H$ boolean representable implies that $\pure(\H)$ is boolean representable. But we can settle the question for low dimensions.

\bp
\label{pclow}
Let $\H$ be a BRSC of dimension $\leq 2$. Then {\rm pure}$(\H)$ is a BRSC.
\ep

\proof
Write $\H\; = (V,H)$. The claim holds trivially for dimension $0$ since $\H$ must be itself pure.

Assume next that $\dim\H\, = 1$. By \cite[Proposition 5.3.1]{RSm}, a simplicial complex $\J\; = (W,J)$ of dimension 1 is boolean representable if and only if the connected components of the graph $\Gamma\J\; = (W, P_2(W) \setminus J)$ are cliques. Now $\Gamma(\pure(\H))$ and $\Gamma(\H)$ differ at most on a few isolated points, corresponding to the facets of $\H$ of dimension 0. Thus the claim holds for dimension 1.

Therefore we may assume that $\dim\H\, = 2$. We define a graph 
$$\Gamma = (V, (P_2(V) \cap \fct\H) \cup (P_2(V) \setminus H)).$$
We show that
\beq
\label{pclow2}
\mbox{the connected components of $\Gamma$ are cliques}.
\eeq

Let $a,b,c \in V$ be distinct and assume that $ab,ac$ are edges of $\Gamma$. We must show that $bc$ is also an edge of $\Gamma$, i.e. there is no $d \in V \setminus bc$ such that $bcd \in H$. We split the discussion into three cases:

\smallskip
\noindent
\underline{{\bf Case 1}}: $ab,ac \in \fct\H$.

\smallskip
\noindent
Suppose that $bcd \in H$. Hence $d \neq a$ since $ab \in \fct\H$. In view of\cite[Theorem 5.2.6]{RSm}, and by symmetry, we may assume that $d \notin \oo{bc}$ or $c \notin \oo{bd}$. 

However, $ab,ac \in \fct\H$ imply $abc, acd \notin H$, hence we successively get $a \in \oo{bc}$ and $d \in \oo{abc} = \oo{bc}$. 

On the other hand, $abd, acd \notin H$ successively yield $a \in \oo{bd}$ and $c \in \oo{abd} = \oo{bd}$. 

Therefore we reach a contradiction.

\smallskip
\noindent
\underline{{\bf Case 2}}: $ab,ac \notin H$.

\smallskip
\noindent
If $ab,ac \notin H$ then $\oo{b} = \oo{a} = \oo{c}$ and so $bc \notin H$ in view of \cite[Theorem 5.2.6]{RSm}. Therefore $bc \notin H$.

\smallskip
\noindent
\underline{{\bf Case 3}}: $ab \in \fct\H$ and $ac \notin H$.

\smallskip
\noindent
Suppose that $bcd \in H$. Hence $d \neq a$ since $ac \notin H$. In view of\cite[Theorem 5.2.6]{RSm}, we may assume that $d \notin \oo{bc}$ or $c \notin \oo{bd}$ or $b \notin \oo{cd}$. 

However, $ac \notin H$ implies $a \in \oo{bc}$ and $abd \notin H$ implies $d \in \oo{abc} = \oo{bc}$.

On the other hand, $abd \notin H$ yields $a \in \oo{bd}$ and $ac \notin H$ yields $c \in \oo{abd} = \oo{bd}$. 

Finally, $ac \notin H$ implies $a \in \oo{cd}$ and $abd \notin H$ implies $b \in \oo{acd} = \oo{cd}$.

Thus we reach a contradiction in all cases.

Therefore we cannot have $bcd \in H$ and so $bc$ is an edge of $\Gamma$, completing the proof of (\ref{pclow2}).

Now write $\pure(\H) = (V',H')$. We show that $\pure(\H)$ is boolean representable. Let $abc$ be a facet of $\pure(\H)$. Since $abc \in H$ and $\H$ is a BRSC, we may assume that the closure $\oo{ab}$ of $ab$ in $\flatx\H$ does not contain $c$. Let $C_a \subseteq V$ denote the connected component of $a$ in $\Gamma$. We claim that
\beq
\label{pclow3}
\emptyset \subset C_a \cap V' \subset \oo{ab} \cap V' \subset V'
\eeq
is a chain in $\flatx(\pure(\H))$.

Indeed, $\emptyset, V' \in \flatx(\pure(\H))$ trivially. Let $d \in C_a \cap V'$ and $e \in V' \setminus C_a$. 
Then $de$ is not an edge of $\Gamma$ and so $de \in H$. Thus $C_a \cap V' \in \flatx(\pure(\H))$.

Next we show that $C_a \cap V' \subset \oo{ab} \cap V'$. Let $d \in C_a \cap V'$. We may assume that $d \neq a$. If $ad \notin H$, then $\oo{d} = \oo{a}$ and so $d \in \oo{ab}$. Hence we may assume that $ad \in \fct\H$. Since $ab \in H$ and $abd \notin H$, we get $d \in \oo{ab}$ also in this case. Thus 
$C_a \cap V' \subseteq \oo{ab}$. Since there is no edge $ab$ in $\Gamma$, it follows from (\ref{pclow2}) that $b \notin C_a$ and so the inclusion is strict.

Finally, we show that $\oo{ab} \cap V'  \in \flatx(\pure(\H))$. Let $X \in H' \cap 2^{\oo{ab}}$ and 
$p \in V' \setminus \oo{ab}$. We may assume that $|X| \geq 1$ . 

Suppose that $X = d \in V'$. Since $\oo{a} \neq \oo{b}$, we must have either $\oo{d} \neq \oo{a}$ or $\oo{d} \neq \oo{b}$, hence $ad \in H$ or $bd \in H$. It follows that $adp \in H$ or $bdp \in H$, hence $dp = X \cup \{ p\} \in H'$.

If $|X| = 2$, then $X \cup \{ p\} \in H \cap P_3(V)$, hence $X \cup \{ p\} \in H'$. 

If $|X| > 2$, then $X \cup \{ p\} \in H \cap P_{\geq 4}(V)$, contradicting $\dim\H\, = 2$. Thus $\oo{ab} \cap V'  \in \flatx(\pure(\H))$. Since $c \in V' \setminus \oo{ab}$, it follows that (\ref{pclow3}) is a chain in $\flatx(\pure(\H))$. Since $abc$ is a transversal of the successive differences for this chain, and $abc$ is an arbitrary facet of a pure complex, then {\rm pure}$(\H)$ is a BRSC.
\qed

\subsection{Second version of the Pure Conjecture}

The second version of the Pure Conjecture may be stated as follows:

\bj
\label{tpct}
Let $\H$ be a BRSC and let $k \geq 1$. Then {\rm pure}$(\H_k)$ is a TBRSC.
\ej

Note that this is equivalent to the statement:
\beq
\label{tpct1}
\mbox{Let $\H$ be a TBRSC and let $k \geq 1$. Then {\rm pure}$(\H_k)$ is a TBRSC.}
\eeq
For the nontrivial implication, let $\H$ be a TBRSC of rank $r$ and assume that Conjecture \ref{tpct} holds. Since $\H$ is a TBRSC, we have $\H = \J_r$ for some BRSC $\J$. 

Suppose first that $k \geq r$. Then $\H_k = \H$, hence we must show that  {\rm pure}$(\H)$ is a TBRSC. Since $\H = \J_r$, our goal follows from applying Conjecture \ref{tpct} to $\J$ and $r$.

Assume now that $k < r$. It is easy to check that $\H_k = \J_k$. Since Conjecture \ref{tpct} implies that {\rm pure}$(\J_k)$ is a TBRSC, then {\rm pure}$(\H_k)$ is a TBRSC and so (\ref{tpct1}) holds.

Therefore (\ref{tpct1}) is equivalent to Conjecture \ref{tpct}.

\be
\label{cepct}
Let $V = 12345678$ and let $\H = (V,H)$ be the BRSC defined by the lattice
$$
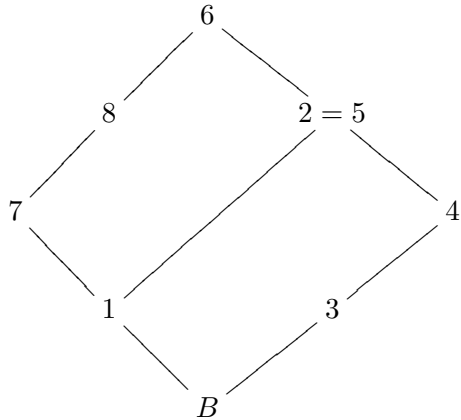
\xymatrix{
&&6 \ar@{-}[dl]  \ar@{-}[dr] &&\\
&8 \ar@{-}[dl] &&{2=5} \ar@{-}[ddll]  \ar@{-}[dr] &\\
7 \ar@{-}[dr] &&&&4 \ar@{-}[dl] \\
&1 \ar@{-}[dr] &&3 \ar@{-}[dl] &\\
&&B&&
}$$
where we associate the points of $V$ to a $\vee$-generating set. Then {\rm pure}$(\H)$ is not a TBRSC.
\ee

Write {\rm pure}$(\H) = \H' = (V,H')$. The maximal chains in the lattice yield the facets
$$1782, 1783, 1784, 1785, 1786, 3416, 3417, 3418, 3426, 3427, 3428, 3456, 3457, 3458,$$
hence $H'$ consists of these facets and their subsets.

Suppose that $\H'$ is a TBRSC.

Consider
$$T(H') = 
\{ T \subseteq V \mid \forall X \in H'_3 \cap 2^T\; \forall p \in V
\setminus T \hspace{.3cm} X \cup \{ p \} \in H'\}.$$
Let $\J = \J(T(H'))$ be defined by the transversals of the successive differences for chains in $T(H')$. By
Theorem \ref{eqtr}, we have $\H' = \J_4$.

We have $3456 \in H'$. Since $3456$ is a transversal of the successive differences for some chain in $T(H')$, there exists some $T \in T(H')$ such that $|T \cap 3456| = 3$. We consider now the four possible cases for $T \cap 3456$, and reach a contradiction in any one of them.

\smallskip

\noindent
\underline{Case 1}: $T \cap 3456 = 345$.

\smallskip

\noindent
Since $345 \subseteq T$, $345 \subset 3456 \in H$ but $3451, 3452 \notin H$, we have $1,2 \in T$. 
Since $12 \subseteq T$, $12 \subset 1782 \in H$ but $126 \notin H$, we have $6 \in T$. 
Hence $3456 \subseteq T$, a contradiction.

\smallskip

\noindent
\underline{Case 2}: $T \cap 3456 = 346$.

\smallskip

\noindent
Since $346 \subseteq T$, $346 \subset 3426 \in H$ but $3467, 3468 \notin H$, we have $7,8 \in T$. 
Since $786 \subseteq T$, $786 \subset 1786 \in H$ but $7865 \notin H$, we have $5 \in T$. 
Hence $3456 \subseteq T$, a contradiction.

\smallskip

\noindent
\underline{Case 3}: $T \cap 3456 = 356$.

\smallskip

\noindent
Since $356 \subseteq T$, $356 \subset 3456 \in H$ but $3567, 3568 \notin H$, we have $7,8 \in T$. 
Since $786 \subseteq T$, $786 \subset 1786 \in H$ but $7864 \notin H$, we have $4 \in T$. 
Hence $3456 \subseteq T$, a contradiction.

\smallskip

\noindent
\underline{Case 4}: $T \cap 3456 = 456$.

\smallskip

\noindent
Since $456 \subseteq T$, $456 \subset 3456 \in H$ but $4567, 4568 \notin H$, we have $7,8 \in T$. 
Since $786 \subseteq T$, $786 \subset 1786 \in H$ but $7863 \notin H$, we have $3 \in T$. 
Hence $3456 \subseteq T$, and we have reached a contradiction in all four cases.

Therefore {\rm pure}$(\H)$ is not a TBRSC.

\bc
\label{spcfails}
Conjecture \ref{tpct} fails for the BRSC of Example \ref{cepct} and $k \geq 4$.
\ec

Next we show that Conjecture \ref{tpct} holds for $k \leq 3$.

\bt
\label{spch}
Let $\H$ be a BRSC and let $1 \leq k \leq 3$. Then {\rm pure}$(\H_k)$ is a TBRSC.
\et

\proof
Suppose first that $k \leq 2$. By Proposition \ref{tbone}, $\H_k$ is a BRSC, therefore {\rm pure}$(\H_k)$ is a BRSC by Proposition \ref{pclow}.

Thus we may assume that $k = 3$. Let 
$$V' = \{ a \in V \mid abc \in H \cap P_3(V) \mbox{ for some } b,c \in V\}.$$
Then $V'$ is the vertex set of {\rm pure}$(\H_3)$. Let $\H^R$ denote the restriction of $\H$ to $V'$. Then  {\rm pure}$(\H_3) = {\rm pure}(\H^R_3)$. Since 
BRSCs are closed under restriction by \cite[Proposition 8.3.1]{RSm}, $\H^R$ is also a BRSC. Therefore we may assume that $V' = V$.

We show next how we can reduce the problem to the simple case. Indeed, suppose that $\H = (V,H)$ is not simple. Let $\H^S$ denote the simplification of $\H$. Following \cite{MRS}, we can provide an easy description of $\H^S$. Assume that $\H$ is represented by the $R \times V$ boolean matrix $M$. Since we are assuming that $P_1(V) \subseteq H$, every column of $\H$ is nonzero, and $ab \notin H$ if and only if the columns corresponding to $a$ and $b$ are equal. Then $\H_S$ is the BRSC represented by the matrix obtained by removing repeated columns from $M$. More precisely, we define an idempotent mapping $\alpha:V \to V$ such that $a\alpha = b\alpha$ if and only if the columns of $a$ and $b$ are equal, and $\H^S = (V\alpha, H \cap 2^{V\alpha})$ is the restriction of $\H$ to $V\alpha$. Moreover, for every $X \subseteq V$, we have
\beq
\label{spch1}
X \in H \mbox{ if and only if $\alpha|_X$ in injective and }X\alpha \in H.
\eeq

Note that $V = V'$ implies that also every vertex of $V\alpha$ occurs in some triangle of $\H^S$.

Write {\rm pure}$(\H_3) = (V,H')$ and {\rm pure}$(\H^S_3) = (V\alpha,H'')$. Note that the vertex set of  the latter complex must be indeed $V\alpha$, in view of (\ref{spch1}).
We claim that, for every $X \subseteq V$:
\beq
\label{spch2}
X \in H' \mbox{ if and only if $\alpha|_X$ in injective and }X\alpha \in H''.
\eeq

Suppose that $X \in H'$. Then $X \subseteq Y$ for some $Y \in H \cap P_3(V)$. It follows from (\ref{spch1}) that $\alpha|_Y$ in injective and $Y\alpha \in H$. Hence $Y\alpha \in H \cap P_3(V\alpha)$ and so $Y\alpha \in H''$, yielding  $X\alpha \subseteq Y\alpha \in H''$ as well.

Conversely, assume that $\alpha|_X$ in injective and $X\alpha \in H''$. Then there exists some $Y \in P_3(V)$ containing $X$ such that $\alpha|_Y$ in injective and $Y\alpha \in H'' \subseteq H$. By (\ref{spch1}), we get $Y \in H$, hence $Y \in H'$ and so $X \in H'$ as well. Thus (\ref{spch2}) holds.

Now if the theorem holds for the simple case, we apply it to the simplification $\H^S$ (which is a BRSC by \cite{MRS}), to deduce that {\rm pure}$(\H^S_3)$ is a TBRSC. We may assume that {\rm pure}$(\H^S_3) = \J_3$ for some BRSC $\J = (V,J)$, represented by some $R' \times V\alpha$ boolean matrix $N$. Let $N'$ be the $R' \times V$ boolean matrix where the column $p \in V$ equals the column $p\alpha$ of $N$. Given $X \in P_{\leq 3}(V)$, it is easy to see that $X$ is recognized by $N'$ if and only if $\alpha|_X$ in injective and $X\alpha$ is recognized by $N$, that is, if and only if $\alpha|_X$ in injective and $X\alpha \in J$. Since {\rm pure}$(\H^S_3) = \J_3$, we can replace $X\alpha \in J$ by $X\alpha \in H''$, and in view of (\ref{spch2}) we get that
\beq
\label{spch3}
\mbox{$X$ is recognized by $N'$ if and only if $X \in H'$}.
\eeq
Let $\J'$ be the BRSC represented by the boolean matrix $N'$. It follows from (\ref{spch3}) that {\rm pure}$(\H_3) = \J'_3$, hence {\rm pure}$(\H_3)$ is a TBRSC as required.

Therefore we only need to deal with the simple case (and recall that we are also assuming that each $p \in V$ occurs in some $X \in H \cap P_3(V)$). Assume so that $\H\; = (V,H)$ is a simple BRSC and {\rm pure}$(\H_3) = (V,H')$. For every $X \subseteq V$, let $\oo{X}$ denote the closure of $X$ in $\flatx\H$. We show the following:
\beq
\label{spch4}
\mbox{if $ab$ is a facet of $\H$, then $\oo{ap} = \oo{bp}$ for every $p \in V \setminus ab$.}
\eeq
Indeed, suppose that $b \notin \oo{ap}$. Since $\H$ is simple, then $ap \in H$ and so $b \notin \oo{ap}$ yields $abp \in H$, contradicting $ab \in fct\H$. Hence $b \in \oo{ap}$ and so $\oo{ap} = \oo{abp}$. By symmetry, we get $\oo{ap} = \oo{abp} = \oo{bp}$ and so (\ref{spch4}) holds.

We define a relation $\tau$ on $V$ by:
$$a \tau b \mbox{ if} \hspace{1cm} a = b \mbox{ or $ab$ is a facet of }\H.$$
This relation is obviously reflexive and symmetric. To show it is transitive, it is enough to consider the case $ab,bc \in \fct\H$ with $a \neq c$. Suppose that $ac \notin \fct\H$. Then $acp \in H$ for some $p \in V \setminus ac$. Since $ab \in \fct\H$, we have $p \neq b$.

Since $\H$ is a BRSC, $acp \in H$ implies one of the cases
$$\mbox{$p \notin \oo{ac}$ or $c \notin \oo{ap}$ or $a \notin \oo{cp}$.}$$
Suppose that $p \notin \oo{ac}$. Since the closure of a facet is always $V$, then $b \notin \oo{ac}$. But $\H$ is simple, hence $ac \in H$ yields $abc \in H$, contradicting $ab \in \fct\H$. Thus, out of symmetry, we may assume that $c \notin \oo{ap}$. But now, since $ab,bc \in \fct\H$, (\ref{spch4}) yields $\oo{ap} = \oo{bp} = \oo{cp}$ and $c \in \oo{ap}$, also a contradiction. Thus $ac \in \fct\H$ and so $\tau$ is an equivalence relation on $V$.

Next we distinguish an element inside each $\tau$-class of $V$. More precisely, we fix an idempotent mapping $\beta:V \to V$ such that $\ker\beta = \tau$. 
We claim that
\beq
\label{spch5}
\mbox{$abc \in H \cap P_3(V)$ if and only if $a'bc \in H \cap P_3(V)$}
\eeq
holds for all $a,a',b,c \in V$ such that $a \tau a'$.

We may assume that $a \neq a'$. Assume that $abc \in H \cap P_3(V)$. If $a' \in bc$, then $abc$ contains two elements in the same $\tau$-class, a contradic tion since two $\tau$-related vertices must constitute a facet. Thus $|a'bc| = 3$. Now $abc \in H \cap P_3(V)$ implies one of the options
$$\mbox{$a \notin \oo{bc}$ or $b \notin \oo{ac}$ or $c \notin \oo{ab}$.}$$ 

Suppose first that $a \notin \oo{bc}$. Since $aa' \in \fct\H$, it follows from (\ref{spch4}) that $\oo{a'b} = \oo{ab}$. But then $a' \in \oo{bc}$ implies $\oo{bc} = \oo{a'bc} = \oo{abc}$, contradicting $a \notin \oo{bc}$. Thus $a' \notin \oo{bc}$ and $bc \in H$ ($\H$ is simple) yields $a'bc \in H$.

Suppose next that $b \notin \oo{ac}$. Since $aa' \in \fct\H$, it follows from (\ref{spch4}) that $\oo{a'c} = \oo{ac}$. But then $b \notin \oo{a'c}$ and $a'c \in H$ yields $a'bc \in H$. 

The case $c \notin \oo{ab}$ follows by symmetry, hence the direct implication of (\ref{spch5}) holds, and so does the opposite implication (also out of symmetry). 

Assume now that $M$ is an $R \times V$ boolean matrix representing $\H$. Let $N$ be the $R \times V$ boolean matrix obtained from $M$ by replacing the column $p$ by the column $p\beta$, for every $p \in V$. Let $\J = (V,J)$ be the BRSC represented by the matrix $N$. We prove that
\beq
\label{spch6}
\mbox{pure$(\H_3) = \J_3$.}
\eeq
 
Let $abc \in H \cap P_3(V)$. Using three times (\ref{spch5}), we obtain $(a\beta)(b\beta)(c\beta) \in H \cap P_3(V)$. Write $X = (a\beta)(b\beta)(c\beta)$. Then $M[R,X]$ contains a nonsingular square submatrix of size 3. Since $N[R,X] = M[R,X]$, the same applies to $N[R,X]$. But $N[R,abc]$ equals $N[R,X]$ up to permutation of columns, so also $N[R,abc]$ contains a nonsingular square submatrix of size 3. Thus $abc \in J$. Therefore $H' \subseteq J$.

Suppose now that $abc \in J \cap P_3(V)$. Then $N[R,abc]$ contains a nonsingular square submatrix of size 3. Write $X = (a\beta)(b\beta)(c\beta)$. 
If $|X| < 3$, then $abc$ contains two $\tau$-related elements, which would constitute a facet, contradiction. Thus $|X| = 3$.
Since $N[R,abc]$ equals $N[R,X]$ up to permutation of columns, so also $N[R,X]$ contains a nonsingular square submatrix of size 3. Since $M[R,X] = N[R,X]$, the same applies to $M[R,X]$, hence $X \in H \cap P_3(V)$. Using three times (\ref{spch5}), we obtain $abc \in H \cap P_3(V)$. 
Therefore $J \cap P_3(V) \subseteq H'$.

Suppose now that $ab \in J \cap P_2(V)$. Then the columns $a$ and $b$ are different in $N$, hence $(a,b) \notin \tau$. Similarly to the preceding case, we get $(a\beta)(b\beta)  \in J \cap P_2(V)$ and $(a\beta)(b\beta)  \in H$. Since $(a,b) \notin \tau$, $ab \notin \fct\H$, hence $abc \in H$ for some $c \in V \setminus ab$. But then $abc \in H'$ and so $ab \in H'$. Since $V$ has been established to be the vertex set of pure$(\H_3)$ anyway, we have shown that $J \cap P_{\leq 3}(V) \subseteq H'$. 

Therefore (\ref{spch6}) holds and so pure$(\H_3)$ is a TBRSC as required.
\qed

\section{Sum of complexes}
\label{sumc}

Let $\H\, = (V,H)$ and $\H' = (V,H')$ be simplicial complexes. The {\em sum} $\H + \H' = (V,H+H')$ is defined by
$$H+H' = \{ I \cup I' \mid I \in H,\; I' \in H'\}.$$

\be
\label{exsum}
If $n \geq 2m \geq 2$, then $U_{m,n} + U_{m,n} = U_{2m,n}$.
\ee

Given $L \in P_{\geq 2}(V) \setminus \{ V \}$, the complex $\B_2(V,L) = (V,B_2(V,L))$ can be described by
$$B_2(V,L) = P_{\leq 2}(V) \cup \{ X \in P_3(V) \; {\big{\lvert}} \; |X \cap L| = 2\}.$$
Indeed, the maximal chains in $P_{\leq 1}(V) \cup \{L,V\}$ are of the form
$$\emptyset \subset a \subset V \quad \mbox{or} \quad \emptyset \subset b \subset L \subset V,$$
where $a \in V \setminus L$.

\bt
\label{sumtwo}
Let $L,L' \in P_{\geq 2}(V) \setminus \{ V \}$. Then the following conditions are equivalent:
\bi
\item[(i)] $\B_2(V,L) + \B_2(V,L')$ is a BRSC;
\item[(ii)] $\B_2(V,L) + \B_2(V,L')$ is a TBRSC;
\item[(iii)] $|L \setminus L'| \leq 3$ or $|L' \setminus L| \leq 3$.
\ei
\et

\proof
(i) $\Rw$ (ii). Trivial.

(ii) $\Rw$ (iii). Write $H = B_2(V,L)$, $H' = B_2(V,L')$, $H'' = H+H'$ and $\H'' = B_2(V,L) + \B_2(V,L')$.

Suppose that $0123 \subseteq L \setminus L'$ and $4567 \subseteq L' \setminus L$. Then $014 \in H$ and $256 \in H'$. Hence $012456 \in H''$. Since $\H''$ is a TBRSC, it follows from Theorem \ref{eqtr} that there exists some $T \in T(H'')$ such that $|T \cap 012456| = 5$. By symmetry, we may assume that $T \cap 012456 = 01245$. Now we successively deduce the following:

Since $01245 \in H'' \cap 2^T$ but $012345 \notin T$, then $3 \in T$.

Since $01234 \in H'' \cap 2^T$ but $012346 \notin T$, then $6 \in T$.

This contradicts $T \cap 012456 = 01245$, hence $|L \setminus L'| \leq 3$ or $|L' \setminus L| \leq 3$.

(iii) $\Rw$ (i). Without loss of generality, we may assume that $|L \setminus L'| \leq 3$. Let $\cl(X)$ (respectively $\cl'(X)$, $\cl''(X)$) denote the closure of $X \subseteq V$ in the lattice of flats of $\B_2(V,L)$ (respectively $\B_2(V,L')$, $\H''$). Since $P_{\leq 2}(V) \subseteq H \cap H'$, we have $P_{\leq 4}(V) \subseteq H''$ and $P_{\leq 3}(V) \subseteq \flatx\H''$. Thus it suffices to show that:
\beq
\label{sumtwo1}
\mbox{if $X \in H'' \cap P_5(V)$, then there exists some $x \in X$ such that $x \notin {\rm Cl}''(X \setminus \{ x\})$};
\eeq
\beq
\label{sumtwo2}
\mbox{if $X \in H'' \cap P_6(V)$, then there exists some $x \in X$ such that $x \notin {\rm Cl}''(X \setminus \{ x\})$};
\eeq

Assume that $X = abcde \in H'' \cap P_5(V)$. By symmetry, we may assume that $abc \in H \cup H'$. But then $abcdy \in H''$ for every $y \in V \setminus abcd$, hence $\cl''(X \setminus \{ e\}) = abcd$ and so (\ref{sumtwo1}) holds.

Before proving (\ref{sumtwo2}), we establish several preliminary results.
\beq
\label{sumtwo3}
\mbox{If $ab \subseteq L \setminus L'$, $c \in L \cap L'$ and $de \in L' \setminus L$, then $abcde \in \flatx\H''$.}
\eeq
First, note that $abcde = abd \cup ce \in H+H' = H''$, so we need to show that $abcdef \in H''$ for every $f \in V \setminus abcde$. 

If $f \in L \cap L'$, we have $abcdef = acd \cup efb \in H +H'$. 
If $f \notin L$, we get $abcdef = acf \cup deb \in H +H'$. If $f \notin L'$, we get $abcdef = abd \cup cef \in H +H'$.
Therefore (\ref{sumtwo3}) holds.

%
%
%
%

\beq
\label{sumtwo4}
L,L' \in \flatx\H''.
\eeq

Let $I \in H'' \cap 2^L$. Then $|I| \leq 4$. Since $P_4(V) \subseteq H''$, we may assume that $I = abcd$. Given $p \in V \setminus L$, we have $I \cup \{ p\} = abp \cup cd \in H + H'$, hence $L \in \flatx\H''$. Similarly, $L' \in \flatx\H''$ and so (\ref{sumtwo4}) holds.

\beq
\label{sumtwo5}
\mbox{If $p \in V \setminus (L\cup L')$, then $(L \cap L') \cup \{ p\} \in \flatx\H''$.}
\eeq

Let $I \subseteq (L \cap L') \cup \{ p\}$ satisfy
$I \in H''$ and let $q \in V \setminus ((L \cap L') \cup \{ p\})$. We need to show that $I \cup \{ q\} \in H''$. It is straightforward to see that we may assume that $|I| = 5$ and $p \in I$. Write $I = abcdp$. 

If $q \notin L$, we have $I \cup \{ q\} = abq \cup cdp \in H +H'$. If $q \notin L'$, we have $I \cup \{ q\} = abp \cup cdq \in H +H'$. Therefore (\ref{sumtwo5}) holds.

\beq
\label{sumtwo6}
\mbox{If $abc \subseteq L \setminus L'$ and $de \subseteq L' \setminus L$, then $abcde \in \flatx\H''$.}
\eeq

First, note that $abcde = abd \cup ce \in H+H' = H''$, so we need to show that $abcdef \in H''$ for every $f \in V \setminus abcde$. 

If $f \notin L$, then $abcdef = abf \cup dec \in H+H'$. If $f \in L$, then $f \in L \cap L'$ in view of $|L \setminus L'| \leq 3$ and we get 
$abcdef = abd \cup efc \in H+H'$. Therefore (\ref{sumtwo6}) holds.

Back to the proof of (\ref{sumtwo2}), we assume now that $X \in H'' \cap P_6(V)$. We may write $X = abc \cup def$, where $a,b \in L$, $c \notin L$, $d,e \in L'$ and $f \notin L'$.

We consider several cases.

\smallskip
\noindent
\underline{Case 1:} $cf \not\subseteq L \cup L'$.
\smallskip

\noindent
\underline{Subcase 1.1:} $abde \subseteq L \cap L'$.
\smallskip

By (\ref{sumtwo5}), we have $(L \cap L') \cup \{ p\} \in \flatx\H''$ for $p \in cf\setminus(L\cup L')$. Taking $x \in cf\setminus \{ p\}$, we get $x \notin {\rm Cl}''(X \setminus \{ x\})$ as required.

\smallskip
\noindent
\underline{Subcase 1.2:} $abde \not\subseteq L \cap L'$.
\smallskip

Take $y \in cf\setminus (L \cup L')$. It suffices to show that $abdey \in \flatx\H$. Let $p \in V \setminus abdey$. If $p \notin L$, we get $abdeyp = abp \cup dey \in H+H'$. If $p \notin L'$, we get $abdeyp = aby \cup dep \in H+H'$. Hence we may assume that $p \in L\cap L'$. Out of symmetry, we can now reduce the discussion to two cases: $a \notin L'$ and $d \notin L$.

If $a \notin L'$, we get $abdeyp = bpy \cup dea \in H+H'$. 
If $d \notin L$, we get $abdeyp = abd \cup epy \in H+H'$. 

\smallskip
\noindent
\underline{Case 2:} $cf \subseteq L \cup L'$.

\smallskip
\noindent
\underline{Subcase 2.1:} $|L \setminus L'| \leq 1$ or $|L' \setminus L| \leq 1$.
\smallskip

We assume that $|L \setminus L'| \leq 1$, the case $|L' \setminus L| \leq 1$ being analogous. Since $cf \subseteq L \cup L'$, we have $X \subseteq L \cup L'$ and so $|X \cap L'| \geq 5$. But $X \subseteq L'$ is obviously impossible, so there exists some $x \in X \setminus L'$ and so $X \setminus \{ x\} \subseteq L' \in \flatx\H''$ by (\ref{sumtwo4}). Therefore $x \notin {\rm Cl}''(X \setminus \{ x\})$.

\smallskip
\noindent
\underline{Subcase 2.2:} $|L \setminus L'| \geq 2$ and $|L' \setminus L| \geq 2$.
\smallskip

Suppose first that $L \cap L' \neq \emptyset$. Then there exists some $x \in X$ such that $X \setminus \{ x\} \in \flatx\H''$ by (\ref{sumtwo3}).

On the other hand, if $L \cap L' \neq \emptyset$, then we must necessarily have $|L \setminus L'| = |L' \setminus L| = 3$. Taking $x \in X \cap L'$, we get $X \setminus \{ x\} \in \flatx\H''$ by (\ref{sumtwo6}). Thus in any case $x \notin {\rm Cl}''(X \setminus \{ x\})$. 

Therefore (\ref{sumtwo2}) holds as required.
\qed

It follows that $\B_2(V,L) + \B_2(V,L')$ is a BRSC if $|V| \leq 7$ and there is a unique counterexample (up to isomorphism) for $|V| = 8$.

\section{Going from dimension 2 to dimension 3}

In this section we explore a systematic way of producing BRSCs of dimension 3 from TBRSCs of dimension 2. 

Let $\H\, = (V,H)$ be a simplicial complex of dimension $d$. An {\em extension} of $\H$ is a simplicial complex $\H' = (V,H')$ such that $\H = \H'_{d+1}$ (i.e. $\H$ is a truncation of $\H'$). If $\H$ is a BRSC of dimension 3, then $\H$ is an extension of $\H_3$, which is a TBRSC of dimension 2. So we can aim at classifying BRSCs of dimension 3 by discussing the extensions of dimension 3 of TBRSCs of dimension 2. 

Recall the definition of $T(H)$ in Section \ref{stru}. 

\bp
\label{ttf}
Let $\H\, = (V,H)$ be a simplicial complex of dimension $\geq 3$. Then:
\bi
\item[(i)] {\rm Fl}$\H \subseteq T(H_3)$;
\item[(ii)] $T(H_3) \cap P_{\leq 2}(V) \subseteq {\rm Fl}\H$;
\item[(iii)] $(T(H_3) \cap P_{3}(V))\setminus H \subseteq {\rm Fl}\H$.
\ei
\ep

\proof
(i) Let $F \in \flatx\H$. Suppose that $X \in (H_{3})_2 \cap 2^F$ and $p \in V
\setminus F$.  Since $(H_3)_2 = H_2 \subseteq H$ and $F \in \flatx\H$, we get $X \cup \{ p \} \in H$ and so $X \cup \{ p \} \in H_3$. Thus $F \in T(H_3)$.

(ii) Let $T \in T(H_3) \cap P_{\leq 2}(V)$. Suppose that $X \in H \cap 2^T$ and $p \in V
\setminus T$.  Since $X \in H_2 = (H_3)_2$ and $T \in T(H_3)$, we get $X \cup \{ p \} \in H_3 \subseteq H$ and so $T \in \flatx\H$.

(iii) Let $T \in (T(H_3) \cap P_{3}(V))\setminus H$. Suppose that $X \in H \cap 2^T$ and $p \in V
\setminus T$.  Since $T \notin H$, we also get $X \in H_2 = (H_3)_2$, yielding $X \cup \{ p \} \in H_3 \subseteq H$ and $T \in \flatx\H$.
\qed

\subsection{The Desargues complex}
\label{Dc}

In this subsection, let $K_5 = (V,E)$ denote the complete graph on $V = 12345$. We discuss the {\em Desargues complex} $\D = (E,D)$, where $D$ consists of all the subforests of $K_5$ with at most 3 edges (including the empty forest). The {\em lines} of $\D$ are the elements of $P_3(E)$ of the form
$$P_2(abc) = \{ ab, ac, bc\},\hspace{1cm}\mbox{where $a,b,c \in V$ are distinct}.$$
Note that these are the 10 lines of the famous Desargues configuration \cite{Wdc}. As subsets of edges of $K_5$, the lines of $\D$ are in fact triangles. If $\L$ denotes the set of lines of $\D$, then
\beq
\label{llines}
D = P_{\leq 3}(E) \setminus \L.
\eeq

We say that a graph is {\em transitive} if its edge set (viewed as a binary relation on the vertex set) is transitive. Equivalently, a transitive graph is a disjoint union of complete graphs (cliques).

\bl
\label{desc}
\bi
\item[(i)] $T(D) = \{ \mbox{transitive subgraphs of }K_5\}$ and is isomorphic to the partition lattice of 12345;
\item[(ii)] ${\rm Fl}\D = \{ \mbox{transitive subgraphs of $K_5$ with at most 3 edges}\} \cup \{E\}$;
\item[(iii)] $\D$ and $\J(T(D))$ are matroids.
\ei
\el

\proof
(i) Let $T \in T(D)$. Suppose that $ab,bc \in T$ are distinct and $ac \notin T$. Then $\{ ab,bc\} \in D_2 \cap 2^T$ and $ac \in E \setminus T$. Since  $T \in T(D)$, we get $\{ ab,ac,bc\} \in D$, a contradiction. Hence $ab,bc \in T$ implies $ac \in T$ and so $T$ is transitive.

Suppose now that $T$ is not transitive. Then there exist distinct $a,b,c \in V$ such that $ab,bc \in T$ but $ac \notin T$. It follows that 
$\{ ab,bc\} \in D_2 \cap 2^T$, $ac \in E \setminus T$. Since $\{ ab,ac,bc\} \notin D$, then $T \notin T(D)$.

It is easy to see that transitive subgraphs of $K_5$ correspond precisely to partitions of 12345: we identify the transitive subgraph $T$ with the partition induced by the cliques of $T$, and inclusion is preserved.

(ii) Let $F \in \flatx\D$. By Lemma \ref{propt}(ii) and part (i), $F$ is transitive. Suppose that $3 < |F| < 10$. Out of symmetry, we may assume that
$$F = P_2(1234) \quad \mbox{or} \quad F = P_2(123) \cup \{ 45 \}.$$

Suppose that $F = P_2(1234)$. Then $\{ 12,23,34\} \in D \cap 2^F$ and $45 \in E \setminus F$, hence $F \in \flatx\D$ yields 
$\{ 12,23,34,45\} \in D$, a contradiction.

Suppose now that $F = P_2(123) \cup \{ 45\}$. Then $\{ 12,23,45\} \in D \cap 2^F$ and $34 \in E \setminus F$, hence $F \in \flatx\D$ yields 
$\{ 12,23,34,45\} \in D$, once again a contradiction.

Therefore $F = E$ or $|F| \leq 3$.

The opposite inclusion follows from Proposition \ref{ttf}.

(iii) Since $\D$ is paving of dimension 2, it suffices to check the exchange property for $I,J \in D$ such that $|J| + 1 = |I| = 3$. In view of (\ref{llines}), we may assume out of symmetry that $J = \{ 12,23\}$. Since $\{ 12,13,23\} \notin D$, we may assume that 4 or 5 occurs in some $ab \in I$, hence $J \cup \{ ab\} \in D$ and so $\D$ is a matroid.

Next we show that
\beq
\label{desc1}
J(T(D)) = \{\mbox{subforests of }K_5\}.
\eeq

Assume that $e_1,\ldots,e_m$ is a transversal of the successive differences for some chain 
$$T_0 \subset T_1 \subset \ldots \subset T_m$$
in $T(D)$ (so that $e_i \in T_i \setminus T_{i-1}$ for $i = 1,\ldots,m$). If this transversal contains a cycle, then there exists some $j > 2$ such that $e_j \subset e_1 \cup \ldots \cup e_{j-1}$ (as subsets of $V$). But then $e_1, \ldots,e_{j-1} \in T_{j-1}$, which is a transitive graph, and so also $e_j \in T_{j-1}$, a contradiction. Therefore our transversal must be a forest.

Conversely, every forest can be enumerated as $e_1,\ldots,e_m$ with $e_i \not\subseteq e_1 \cup \ldots \cup e_{i-1}$ for every $i$. It is immediate that $e_i$ does not belong to the transitive closure $T_i$ of $e_1,\ldots,e_i$. Thus $e_1,\ldots,e_m$ is a transversal of the successive differences for the chain 
$T_0 \subset \ldots \subset T_m$, and so (\ref{desc1}) holds. It follows that $\J(T(D))$ is the graphic matroid defined by $K_5$.
\qed

Note that $\flatx\D$ consists of:
\bi
\item
the empty set and the full set $E$;
\item 
singletons and lines;
\item
pairs of two disjoint edges.
\ei
The latter type (which we may call {\em short lines}) correspond to the edges of the {\em Petersen graph} (viewed as a Kneser graph) \cite{Wpg}.

There exist thousands of boolean representable extensions of $\D$ with dimension 3. We prove next that only one of them is a matroid.

\bt
\label{med}
$\J(T(D))$ is the only proper matroid extension of $\D$.
\et

\proof
Suppose that $\H$ is a proper matroid extension of $\D$. Then $\D = \H_3$. 
By Proposition \ref{ttf}, we have 
\beq
\label{desc3}
T(D) \cap (P_{\leq 2}(E) \cup (P_{3}(E) \setminus H)) \subseteq \flatx\H \subseteq T(D).
\eeq
It follows from Lemma \ref{desc}(i) that $\flatx\H$ contains all the transitive subgraphs of $K_5$ with at most 3 edges, and $E \in \flatx\H$ trivially. We show next that 
\beq
\label{desc2}
\mbox{if $T \in T(D)$ satisfies $3 < |T| <10$, then $T \in \flatx\H$.}
\eeq

Note that $T(D)$ is a geometric lattice by Lemma \ref{desc}(i), and the height of an element may be any number from 0 to 4:
\bi
\item
height 0: the empty graph
\item
height 1: one edge
\item
height 2: two disjoint edges, or one triangle
\item
height 3: a 4-clique, or the disjoint union of a triangle with an edge
\item
height 4: the full graph
 \ei

By Lemma \ref{desc}, $\flatx\D$ contains all the elements of $T(D)$ except those of height 3 (which are the transitive graphs $T$ satisfying $3 < |T| <10$). Since $\H$ is a proper matroid extension of $\D$, then $\flatx\H$ must contain some $T_0 \in T(D)$ of height 3. Let $T \in T(D)$ be another element of height 3. It is easy to check that $T \cap T_0 \neq \emptyset$ in every possible case.

Assume first that $|T \cap T_0| > 1$ (so has height 2). Suppose that $T \notin \flatx\H$. Let $S \in T(D) \setminus \{ T\cap T_0\}$ be covered by $T$ (every element of height 3 in $T(D)$ covers at least 4 elements of height 2). We claim that
$$
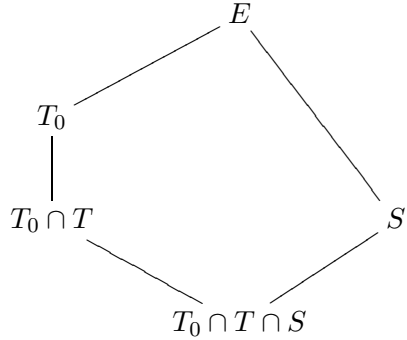
\xymatrix{
\ &E \ar@{-}[dl] \ar@{-}[ddr]&\\
T_0 \ar@{-}[d] &&\\
T_0 \cap T \ar@{-}[dr] && S \ar@{-}[dl]\\
&T_0 \cap T \cap S&
}$$
is a sublattice of $\flatx\H$.

Indeed, it suffices to check that 
\beq
\label{med1}
T_0 \cap S = T_0 \cap T \cap S \quad \mbox{and} \quad (T_0 \cap T) \vee S = E
\eeq
hold in $\flatx\H$. 

Since $S \subset T$, we have $T_0 \cap S \subseteq T_0 \cap T$ and so $T_0 \cap S = T_0 \cap T \cap S$. On the other hand, $(T_0 \cap T) \vee S = T$ in $T(D)$ since $T$ covers both $S$ and $T_0 \cap T$, which are distinct elements of height 2. It follows that the join of $(T_0 \cap T) \vee S$ in $\flatx\H$ is $\geq T$. But $T$ is a co-atom of $T(D)$ which is not in $\flatx\H$, hence $(T_0 \cap T) \vee S = E$ holds in $\flatx\H$. 

Thus (\ref{med1}) holds and so $\flatx\H$ is not semimodular, and consequently not geometric. This contradicts the fact that $\H$ is a matroid, whence  $T \in \flatx\H$. 

We must now deal with the case $|T \cap T_0| = 1$. Let $T_1$ be a 4-clique. It is easy to check that $|T_1 \cap T'| > 1$ for every $T' \in T(D)$ of height 3, so we can use the preceding case to yield successively $T_1 \in \flatx\H$ and $T \in \flatx\H$. Therefore (\ref{desc2}) holds.

Together with (\ref{desc3}), this yields $\flatx\H = T(D)$. Hence the faces of both $\H$ and $\J(T(D))$ are the transversals of the successive differences for chains in $\flatx\H = T(D)$. Therefore $\H = \J(T(D))$.
\qed

\subsection{The non Desargues complex}

In this subsection, we keep the notation of Subsection \ref{Dc} and discuss the {\em non Desargues complex} $\Ncal = (E,N)$. We fix $L_0 = \{ 34,35,45\}$ and set $N = D \cup L_0$.

\bl
\label{ndesc}
\bi
\item[(i)] $T(N) = \{$subgraphs of $K_5$ where each connected component is either a clique or a 2-subset of $L_0 \}$;
\item[(ii)] ${\rm Fl}\Ncal = ({\rm Fl}\D \setminus \{ L_0 \} ) \cup P_2(L_0)$;
\item[(iii)] $\Ncal$ is a matroid.
\ei
\el

\proof
(i) Let $T \in T(N)$. If $abc \in P_3(V) \setminus \{ 345 \}$, we claim that
\beq
\label{ndesc1}
ab,bc \in T \; \mbox{ implies } \; ac \in T.
\eeq
Indeed, we have $\{ ab,bc \} \in N \cap 2^T$, hence $ac \in E \setminus T$ would imply $\{ ab,bc, ac \} \in N$ a contradiction. Hence (\ref{ndesc1}) holds. Let $C$ be a connected component of $T$ (seen as a subgraph of $K_5$). It follows from (\ref{ndesc1}) that $C$ is a clique unless $|C \cap L_0| = 2$. What can $C$ be then? Out of symmetry, we may assume that $C \cap L_0 = \{ 34, 45\}$. Suppose that $C \not\subseteq L_0$. Since $C$ is connected, we may assume that $ab \in C$ with $a \in 12$ and $b \in 345$. Using successively (\ref{ndesc1}), we deduce that $a3, a4, a5 \in T$, which imply then that $35 \in T$, a contradiction. Therefore $C$ is either a clique or a 2-subset of $L_0$. 

Now let $X \in 2^E \setminus T(N)$. Then there exist $ab,cd \in X$ and $ef \in E \setminus X$ such that $\{ ab,cd,ef \} \notin N$. Then 
$\{ ab,cd,ef \} \in \L \setminus \{ L_0\}$ and we may assume without loss of generality that $d = b$ and $ef = ac$. But then the connected component of $X$ containing $ab$ and $bc$ is neither a clique nor a 2-subset of $L_0$.

(ii) By \cite[Proposition 4.2.3]{RSm}, no $F \in \flatx\Ncal \setminus \{ E\}$ can contain a facet of dimension 2. It follows easily that $\flatx\Ncal$ contains $E$ and the elements of $T(N)$ which contain no facet of dimension 2. With respect to the characterization of $T(N)$ in part (i), this excludes:
\bi
\item
all 4-cliques;
\item
the 3-clique $L_0$;
\item
the simultaneous presence of a 3-clique and a 2-clique;
\item
the simultaneous presence of a 2-clique and a 2-subset of $L_0$.
\ei
Straightforward checking shows that we are left precisely with $({\rm Fl}\D \setminus \{ L_0 \} ) \cup P_2(L_0)$.

(iii) Since $\Ncal$ is paving of dimension 2, it suffices to check the exchange property for $I,J \in D$ such that $|J| + 1 = |I| = 3$. I
We may assume that $J$ is connected, otherwise $J$ is a flat and $J \cup \{ i\} \in N$ for every $i \in I \setminus J$. Assume then that $J = \{ ab, bc\}$. We may also assume that $J \not\subset I$, but in this case some $k \in 12345 \setminus abc$ occurs in some $i \in I$. Thus $J \cup \{ i\} \in N$ and so $\Ncal$ is a matroid.
\qed

There exist thousands of boolean representable extensions of $\Ncal$ with dimension 3. We prove next that none of them is a matroid.

\bt
\label{men}
There exists no proper matroid extension of $\Ncal$.
\et

\proof
Suppose that $\H$ is a proper matroid extension of $\Ncal$. Then $\Ncal = \H_3$. Since $\H$ is a proper extension, we have $\dim\H \geq 3$.

Note that $\flatx\H \subseteq T(N)$ by Proposition \ref{ttf}(i). Now
$$\emptyset \subset \{ 12 \} \subset P_2(123) \subset P_2(1234) \subset E$$
is a maximal chain in $T(N)$ in view of Lemma \ref{ndesc}(i). Since $\H$ is a matroid, then $\flatx\H$ is a geometric lattice, so in particular all maximal chains in $\flatx\H$ have the same length (every semimodular lattice satisfies the Jordan-Dedekind property). Thus $\dim\H\, = 3$.

By Proposition \ref{ttf}, 
\beq
\label{men1}
\emptyset \subset \{ 34 \}  \subset \{ 34,45 \}  \subset E 
\eeq
is a chain of length 3 in $\flatx\H \subseteq T(N)$. By the Jordan-Dedekind property, this chain may be refined to a chain of length 4, which is then a maximal chain. It follows from Lemma \ref{ndesc}(i) that the only possible candidates are:
\bi
\item[(A)]
$P_2(i345)$ with $i = 1$ or $2$;
\item[(B)]
$L_0$;
\item[(C)]
$L_0 \cup \{ 12\}$;
\item[(D)]
$\{ 34,45,12\}$.
\ei
We claim that, for each one of these cases, $\flatx\H$ must contain one of the following sublattices:
$$
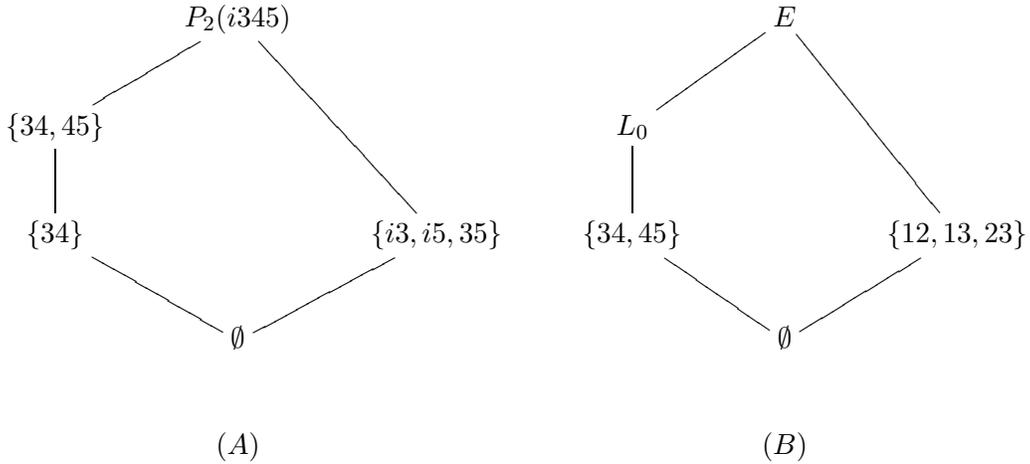
\xymatrix{
&P_2(i345) \ar@{-}[dl] \ar@{-}[ddr]&&&E \ar@{-}[dl] \ar@{-}[ddr]&\\
\{ 34,45 \} \ar@{-}[d] &&&L_0 \ar@{-}[d] &&\\
\{ 34 \} \ar@{-}[dr] && \{ i3,i5,35\} \ar@{-}[dl]&\{ 34,45\} \ar@{-}[dr] && \{ 12,13,23\} \ar@{-}[dl]\\
&\emptyset&&&\emptyset&\\
&(A)&&&(B)&
}$$
$$\xymatrix{
&E \ar@{-}[dl] \ar@{-}[ddr]&&&E \ar@{-}[dl] \ar@{-}[ddr]&\\
L_0 \cup \{ 12\} \ar@{-}[d] &&&\{ 34,45,12\} \ar@{-}[d] &&\\
\{ 34,45 \} \ar@{-}[dr] && \{ 13,24\} \ar@{-}[dl]&\{ 34,45\} \ar@{-}[dr] && \{ 13,24\} \ar@{-}[dl]\\
&\emptyset&&&\emptyset&\\
&(C)&&&(D)&
}$$
Indeed, we always get empty intersections for elements from distinct sides. In case (A), we have $\{ 34 \} \vee \{ i3,i5,35\} = P_2(i345)$ because $P_2(i345)$ is the smallest clique of $K_5$ containing the edges $\{ 34,i3,i5,35\}$ (and 2-subsets $L_0$ won't do either!). In case (B), we have $\{ 34,45 \} \vee \{ 12,13,23\} = E$ because there exist no proper clique of $K_5$ containing the edges $\{ 34,45,12,13,23\}$ (let alone 2-subsets of $L_0$). In cases (C) and (D), we have $\{ 34,45 \} \vee \{ 13,24\} = E$ because there exist no proper clique of $K_5$ containing the edges $\{ 34,45,13,24\}$.  

Thus in any possible case we have shown that $\flatx\H$ is not semimodular, and consequently not geometric. This contradicts the fact that $\H$ is a matroid, therefore there exists no proper matroid extension of $\Ncal$.
\qed

The non Desargues complex shows that things are much more complicated when we go from dimension 2 to dimension 3, with respect to the transition from dimension 1 to dimension 2.

Consider also the problem of existence of a proper matroid extension for a simplicial complex $\H\, = (V,H)$ of dimension $d$. Since matroids are closed under truncation (the exchange property is trivially inherited), there exists a proper matroid extension of $\H$ if and only if there exists a proper matroid extension of $\H$ of dimension $d+1$. Moreover, in this case $\H$ must be a matroid itself, being the truncation of a matroid. But there exists another obvious necessary condition: if $\H$ admits a proper matroid extension, then there exists some $X \in P_{d+2}(V)$ such that $P_{d+1}(X) \subset H$. The next result shows that these conditions suffice to characterize the situation at dimension 1. Recall that a simplicial complex of dimension 1 can be viewed as a graph with at least one edge.

\bp
\label{caseone}
Let $\H\, = (H,V)$ be a matroid of dimension 1. Then the following conditions are equivalent:
\bi
\item[(i)] $\H$ admits a proper matroid extension;
\item[(ii)] $\H$ is not a complete bipartite graph;
\item[(iii)] there exists some $X \in P_{3}(V)$ such that $P_{2}(X) \subset H$.
\ei
\ep

\proof
It is known that a simplicial complex $\H\, = (V,H)$ of dimension 1 is a matroid if and only it is boolean representable if and only if $V$ admits a nontrivial partition $V = V_1 \cup \ldots \cup V_m$ such that 
\beq
\label{caseone1}
H \cap P_2(V) = P_2(V) \setminus (\cup_{i=1}^m P_2(V_i))
\eeq
(so we have a complete $m$-partite graph). Indeed, the second and the third conditions are equivalent by \cite[Proposition 5.3.1]{RSm}, and every matroid is a BRSC by \cite[Theorem 5.2.10]{RSm}. It is immediate that (\ref{caseone1}) implies the exchange property, so the three conditions above are indeed equivalent.

We may assume then that $V$ admits a nontrivial partition $V = V_1 \cup \ldots \cup V_m$ such that 
(\ref{caseone1}) holds. 

(i) $\Rw$ (iii). Let $\J\, = (V,J)$ be a proper extension of $\H$. Then $\dim\J\, \geq 2$ and $\J_2 = \H$. Let $X \in J \cap P_3(V)$. Then $P_2(X) \subseteq J \cap P_{\leq 2}(V) = H$ and (iii) holds.

(iii) $\Rw$ (ii). Suppose that $\H$ is a complete bipartite graph. Then $m = 2$. Let $X \in P_{3}(V)$. Then there exists distinct $x,y \in X$ and $i \in \{ 1,2\}$ such that $x,y \in V_i$. Hence $xy \notin H$ and so $P_{2}(X) \not\subset H$. Therefore (iii) fails.
 
(ii) $\Rw$ (i). If $\H$ is not a complete bipartite graph, then $m > 2$. Let $\J \, = (V,J)$ be defined by
$$J = \{ Y \in P_{\leq 3}(V) \mid \mbox{no two elements of $Y$ belong to the same }V_i\}.$$
Then $\H\, = \J_2$ and $H \subset J$ since $m > 2$. It remains to show that $\J$ is a matroid. Since $\J_2 = \H$ is a matroid itself, we only need to check the exchange property for $X,Y \in J$ such that $|X| = |Y|+1 = 3$.

Now the two elements of $Y$ belong to two distinct $V_j$, say $V_{j_1}$ and $V_{j_2}$, and the three elements of $X$ belong to three distinct $V_i$, say $V_{i_1}$, $V_{i_2}$ and $V_{i_3}$. Without loss of generality, we may assume that $i_3 \notin \{ j_1,j_2\}$. If $x \in X \cap V_{i_3}$, we get $Y \cup \{ x\} \in J$ and $x \in X\setminus Y$, thus the exchange property holds as required.
\qed

Now $\Ncal$ stands as a counterexample for a possible generalization of Proposition \ref{caseone}: $\Ncal \, = (E,N)$ is a matroid of dimension 2, there exists $X \in P_{4}(E)$ such that $P_{3}(X) \subset N$, and yet $\Ncal$ admits no proper matroid extension.

\section{Matroids of codimension 1}

Assume that $\H\, = (V,H)$ is a TBRSC of dimension $d$. Then $\H\; = \H'_{d+1}$ for some BRSC $\H' = (V,H')$. By Theorem \ref{eqtr}, we have $\H\, = (\J(T(H)))_{d+1}$. We claim that
\beq
\label{trex}
\flatx\H' \subseteq T(H) \quad\mbox{and}\quad H' \subseteq J(T(H)).
\eeq
Indeed, let $F \in \flatx\H'$. Let $X \in H_d \cap 2^F$ and $p \in V \setminus F$. Since $H_d \subset H'$ and $F \in \flatx\H'$, we get $X \cup \{ p\} \in H'$. Hence $X \cup \{ p\} \in H'_{d+1} = H$ and so $F \in T(H)$. Thus $\flatx\H' \subseteq T(H)$. Now $H' = J(\flatx\H')$ since $\H'$ is a BRSC, and so $\flatx\H' \subseteq T(H)$ yields $H' = J(\flatx\H') \subseteq J(T(H))$. Therefore (\ref{trex}) holds.

It follows from (\ref{trex}) that $\dim\H' \leq \dim \J(T(H))$ whenever a BRSC $\H'$ is a extension of $\H$. The {\em codimension} of $\H$ is defined as
$$\codim\H = \dim\J(T(H)) - \dim\H.$$

It is easy to see that a TBRSC of codimension 0 must be a BRSC and admits no proper extensions. Thus we turn our attention into codimension 1. 

\bt
\label{codone}
Let $\H\, = (V,H)$ be a matroid of codimension 1. Then $\H$ has at most one proper matroid extension, which will then be $\J(T(H))$.
\et

\proof
Assume that $\dim\H\, = d$. Suppose that $\H' = (V,H')$ is a proper matroid extension of $\H$. Since matroids are boolean representable, it follows from (\ref{trex}) that $H' \subseteq J(T(H))$. It remains to be proved that $J(T(H)) \subseteq H'$. We denote by $\clos:2^V \to \flatx\H'$ the closure operator on $2^V$ induced by $flatx\H'$ (so that $\clos(X) = \cap\{ F \in \flatx\H' \mid X \subseteq F\}$).  

Let $X \in J(T(H))$. Suppose first that $|X| \leq d+1$. By Theorem \ref{eqtr}, we have $X \in H$ and so $X \in H'$. Since $\dim \J(T(H)) = d+1$, we may assume then that $|X| = d+2$. Then there exists some chain
$$T_0 \subset T_1 \subset \ldots \subset T_{d+2}$$
in $T(H)$ and an enumeration $x_1,\ldots,x_{d+2}$ of the elements of $X$ such that $x_i \in T_i \setminus T_{i-1}$ for $i = 1,\ldots,d+2$. Let $X' = X \setminus \{ x_{d+2} \}$. By the first case, we get $X' \in H \subseteq H'$. Now let $Y \in H' \cap P_{d+2}$ (it exists since $\H'$ is a proper extension of $\H$). Since $\H'$ is a matroid, we may apply the exchange property to $X',Y \in H'$, hence $X' \cup \{ y\} \in H'$ for some $y \in Y \setminus X'$. Since $\H'$ is a matroid, any enumeration of the elements of $X' \cup \{ y\}$ can be a transversal of the successive differences for some chain in $\flatx\H'$ \cite{Oxl}. It follows that $y \notin \clos(X')$. Write $C = \clos(X')$.  

Suppose that $x_{d+2} \in C$. Then $x_1,\ldots,x_{d+2},y$ becomes a 
transversal of the successive differences for 
$$T_0 \cap C \subset T_1 \cap C \subset \ldots \subset T_{d+2}\cap C \subseteq V,$$
which is a chain in $T(H)$, in view of (\ref{trex}) and Lemma \ref{propt}(i). Thus $X \cup \{ y\} \in J(T(H))$, contradicting $\dim \J(T(H)) = d+1$. 

Therefore $x_{d+2} \notin C \in \flatx\H'$. Since $X' \in H' \cap 2^C$, we get $X \in H'$ and so $J(T(H)) \subseteq H'$ as required.
\qed

It follows easily from (\ref{desc1}) that the Desargues complex $\D$ has codimension 1. Note that $\J(T(D))$ is a matroid by Lemma \ref{desc}(iii), so Theorem \ref{codone} provides an alternative proof for Theorem \ref{med}. 

The next example shows that proper matroid extensions do not always exist at codimension 1.

\be
\label{triang}
Let $\H\, = (V,H)$ be defined by $V = 123456$ and $H = P_{\leq 3}(V) \setminus \{ 124,135,236\}$. Then:
\bi
\item[(i)] $\H$ is a matroid of dimension 2 and codimension 1;
\item[(ii)] $\H$ admits no proper matroid extension.
\ei
\ee

Indeed, let $I,J \in H$ with $|I| = |J| +1$. We may assume that $|I| = 3$ and $J \not\subseteq I$. Then $|I \setminus J| \geq 2$. Since any two of the three forbidden 3-subsets will intersect at just one point, it follows that $J \cup \{ i \} \in H$ for some $i \in I \setminus J$. Thus $\H$ is a matroid (of dimension 2).

Write $\L = \{ 124,135,236 \}$ (lines). We claim that
\beq
\label{triang1}
T(H) = \{ X \subseteq V \mid |X \cap L| \leq 1 \mbox{ for every }L \in \L\} \cup \L \cup \{ V\}.
\eeq

Indeed, let $T \in T(H)$. We may assume that $|T \cap L| \geq 2$ for some $L \in \L$. Since $L \notin H$, it follows easily that $L \subseteq T$. Now we may assume that $L \subset T$. Then $|T \cap L'| \geq 2$ for some $L' \in \L \setminus \{ L\}$. Since $L' \notin H$, we get $L' \subseteq T$. The same argument shows that the third line is contained in T, hence $T = V$.

The proof for the opposite inclusion is straightforward, therefore (\ref{triang1}) holds.

Suppose that $\J(T(H))$ is a matroid. It is easy to check that 
\beq
\label{triang2}
\emptyset \subset 4 \subset 45 \subset 456 \subset V
\eeq
is a chain in $T(H)$ of maximum length, hence $\codim\H\, = 1$. Moreover, $4,5,6,1$ is a transversal of the successive differences for (\ref{triang2}), hence $I = 1456 \in J(T(H))$. On the other hand, $1,2,3$ is a transversal of the successive differences for another chain in $T(H)$, e.g.
$$\emptyset \subset 1 \subset 124 \subset V,$$
hence $J = 123 \in J(T(H))$. 

By the exchange property, $J \cup \{ i\} \in J(T(H))$ for some $i \in I \setminus J$. Out of symmetry, we may assume that $i = 4$. But then there exists some $T \in T(H)$ such that $|T \cap 1234| = 3$, which implies $T = 123$ in view of (\ref{triang1}). But then $T$ must contain some $T' \in T(H)$ such that $|T' \cap 123| = 2$, contradicting (\ref{triang1}). Thus $\J(T(H))$ is not a matroid. By Theorem \ref{codone}, $\H$ admits no proper matroid extension.

\medskip

We consider next codimension 2. We start by returning to the non Desargues complex $\Ncal$. In view of Lemma \ref{ndesc}(i), it is easy to see that
$$\emptyset \subset \{ 34 \} \subset \{ 34, 45 \} \subset \{ 34,35,45 \} \subset \{ 12,34,35,45 \} \subset E$$
is a chain in $T(N)$, and there are no longer chains. Hence $\dim\J(T(N)) = 5$ and so $\codim\Ncal = 2$. Thus Theorem \ref{men} shows that it is possible that a matroid of codimension 2 admits no proper matroid extension. 

The next example shows that a matroid of dimension 2 and codimension 2 can admit several proper matroid extensions of dimension 3.

\be
\label{sme}
Let $\H\, = (V,H)$ be defined by $V = 123456$ and $H = P_{\leq 3}(V) \setminus \{ 456\}$. Then:
\bi
\item[(i)] $\H$ is a matroid of dimension 2;
\item[(ii)] $\J(T(H))$ is a matroid of dimension 4;
\item[(iiii)] if $Q = \{ X \in P_{\leq 4}(V) \mid 456 \not\subseteq X\}$, then $\Q_k = (V, Q \setminus \{123k\})$ is a matroid extension of $\H$ of dimension 3 for $k = 4,5,6$.
\ei
\ee

To check the exchange property for $\H$, it is enough to consider $I,J \in H$ with $|I| = |J|+1 = 3$. We may assume that $J \not\subseteq I$, hence $|I\setminus J| \geq 2$. Since 456 is the only 3-subset of $V$ which is not in $H$, it follows that $J \cup \{ i\} \in H$ for some $i \in I \setminus J$. Thus $\H$ is a matroid (of dimension 2).

It is straightforward to check that
$$T(H) = \{ T \subseteq V : |T \cap 456| \neq 2\}.$$
We claim that 
\beq
\label{sme1}
J(T(H)) = \{ X \subseteq V : |X \cap 456| \leq 2\}.
\eeq

Indeed, the direct inclusion follows from the fact that every $T \in T(H)$ containing two elements from 456 must contain the three of them. For the opposite inclusion, and out of symmetry, it suffices to check that $12345 \in J(T(H))$. This follows from $1,2,3,4,5$ being a transversal of the successive differences for the chain
$$\emptyset \subset 1 \subset 12 \subset 123 \subset 1234 \subset V$$
in $T(H)$. Therefore (\ref{sme1}) holds.

We check now the exchange property for $\J(T(H))$. Let $I,J \in J(T(H))$ with $|I| = |J|+1$. We may assume that $|J \cap 456| = 2$, otherwise we may pick any element of $I \setminus $ to get $I \cup \{ j\} \in J(T(H))$. But then $|I \setminus 456| > |J \setminus 456|$ and $I \cup \{ j\} \in J(T(H))$ holds for some $i \in I \setminus (J \cup 456)$. Therefore $\J(T(H))$ is a matroid, and has dimension 4 since 12345 is a facet of maximum dimension.

Now fix $k \in 456$ and write $Q_k = Q \setminus \{123k\}$. Let $I,J \in Q_k$ with $|I| = |J|+1$. Since $P_2(V) \setminus Q_k$, we may assume that $|J| \geq 2$. Suppose that $|J| = 2$. If $J \not\subseteq 456$, then $J \cup \{ i\} \in Q_k$ for any $i \in I \setminus J$, hence we may assume that $J \subseteq 456$. But then $I \not\subseteq 456$ and we get $J \cup \{ i\} \in Q_k$ for some $i \in I \setminus 456$.

Finally, we are left with the case $|J| = 3$. Suppose first that $|J \cap 456| = 0$. Then $J = 123$. Since $I \neq 123k$, we get $J \cup \{ i\} \in Q_k$ for any $i \in I \setminus 123k$. Suppose now that $|J \cap 456| = 2$. Since $|I \cap 123| \geq 2$, we get $J \cup \{ i\} \in Q_k$ for any $i \in I \setminus (J \cup 456)$. Finally, we assume that $|J \cap 456| = 1$. If $I \not\subseteq J \cup 123$, then we get $J \cup \{ i\} \in Q_k$ for any $i \in I \setminus (J \cup 123)$. If $I \subseteq J \cup 123$, then $|J \cap 456| = |I \cap 456| = 1$ yields $J \subset I$ and the exchange property holds in every case. Therefore $\Q_k$ is a matroid, indeed a matroid extension of $\H$ of dimension 3.

\medskip

However, part of Theorem \ref{codone} can still be true:

\bq
\label{ctnm}
Is there a matroid $\H\, = (V,H)$ of codimension 2 such that  $\J(T(H))$ is not a matroid but $\H$ admits a proper matroid extension?
\eq

\section{On the Dowling and Rhodes matroids}

\subsection{The reduced Rhodes matroid}
\label{rrm}

We recall the definition of the {\em reduced Rhodes matroid} from \cite{MRS4}.
Let $G$ be a finite nontrivial group and fix an integer $n \geq 2$. Let $\un{n} = \{ 1,\ldots,n\}$. Given $I \subseteq \un{n}$, we denote by $F(I,G)$ the collection of all functions
$f:I \rightarrow G$. Given a partition $\pi$ of $I \subseteq \un{n}$ and $f,h \in F(I,G)$, we write 
$$f \sim_{\pi} h \quad \mbox{ if $f|_{\pi_i} \in G(h|_{\pi_i})$ for each block $\pi_i$ of $\pi$.}$$
Then $\sim_{\pi}$ is an equivalence relation on $F(I,G)$. Let $[f]_{\pi}$ denote the $\sim_pi$ class of $f \in F(I,G)$. We define $\SPC(\un{n},G)$ as the set of triples of the form $(I,\pi,[f]_{\pi})$, where $I \subseteq \un{n}$, $\pi$ is a partition of $I$ and $f\in F(I,G)$.

Let $V$ denote the set of all $(I,\pi,[f]_{\pi}) \in \SPC(\un{n},G)$ such that $|I| = 2$ and $\pi$ is the trivial partition (one single block). Then $V$ is the set of atoms of the reduced Rhodes lattice defined by $G$ and $n$ (see  \cite{MRS4}). If $I = \{ i,j\}$ and $(f(i))\inv f(j) = g \in G$, we define $\Gamma(I,\pi,[f]_{\pi})$ to be the directed labeled graph 
\beq
\label{fever}
\xymatrix{
i \ar@/^/[rr]^{g}&& j \ar@/^/[ll]^{g\inv}
}
\eeq
Assuming the necessity of the inverse edges, it is of course enough to represent the above graph by the one-edge graph $i \mapright{g} j$ (or $j \mapright{g\inv} i$). To simplify notation, we shall denote by $(i,g,j)$ (or $(j,g\inv,i)$) the unique $a \in V$ such that $\Gamma(a)$ is $i \mapright{g} j$.

Now each $Z \subseteq V$ can be represented by a directed labeled graph $\Gamma(Z)$. We denote by $\Gamma_0(Z)$ the undirected multigraph with the same vertex set as $\Gamma(Z)$ and an edge $i \edge j$ for each pair of inverse edges (\ref{fever}) in $\Gamma(Z)$.

We say that a multigraph $\Gamma$ is {\em unicyclic} if:
\bi
\item
$\Gamma$ is connected;
\item
$\Gamma$ has no loops;
\item
the number $v$ of vertices equals the number $e$ of edges in $\Gamma$.
\ei 
Since finite trees can be characterized as connected graphs satisfying $e = v-1$ \cite[Theorem 5.1.2]{CM}, it follows that unicyclic graphs are precisely those graphs which can be obtained from a tree by adding a new edge connecting two distinct vertices. 

If $Z \subseteq V$, we say that $\Gamma(Z)$ is $G$-{\em trivial} if, for every cycle
$$q_0 \mapright{g_1} q_1 \mapright{g_2} \ldots \mapright{g_m} q_m = q_0$$
in $\Gamma(Z)$, the equality $g_1g_2\ldots g_m = 1$ holds in $G$. Otherwise, $\Gamma(Z)$ is $G$-{\em nontrivial}. It is a simple exercise to show that $\Gamma(Z)$ is trivial whenever $\Gamma_0(Z)$ is a tree. On the other hand, if $\Gamma_0(Z)$ is unicyclic, $G$-triviality is determined by the unique nontrivial cycle of $\Gamma(Z)$ (we call a cycle of the form (\ref{fever})
trivial).  

Let $H$ denote the set of all $Z \subseteq V$ such that
\bi
\item
each connected component of $\Gamma_0(Z)$ is either a tree or a unicyclic graph:
\item
at most one connected component of $\Gamma_0(Z)$ is unicyclic;
\item
if $\Gamma_0(Z)$ has a unicyclic connected component, then $\Gamma(Z)$ is $G$-nontrivial.
\ei
Then $\H\, = (V,H)$ is reduced Rhodes matroid defined by $G$ and $n$. Note that $\H$ has dimension $n-1$ and every facet has $n$ elements.

Given a partition $\pi$ of $I \subseteq \oo{n}$, let $C(\pi)$ denote the set of all 
$(\{ i,j\}, \pi,[f]_{\pi}) \in V$ such that $i,j$ belong to the same $\pi$ class. 

\bl
\label{theflatsT}
The flats of $\H$ are the sets of the following two types:
\bi
\item[(i)] 
$C(\pi)$, where $\pi$ is a partition of $I \subseteq \un{n}$;
\item[(ii)] 
$Z \subseteq V$ such that $\Gamma_0(Z)$ is a union of cliques and $\Gamma(Z)$ is $G$-trivial.
\ei
\el

Recall that $\H_k = (V,H_k)$ denotes the truncation defined by $H_k = H \cap \P_{\leq k}$. Hence 
$$T(H_k) = \{ T \subseteq V \mid \forall X \in H_{k-1} \cap 2^T\; \forall p \in V
\setminus T \hspace{.3cm} X \cup \{ p \} \in H_k\}.$$

\bl
\label{intrun}
Let $\H$ be a BRSC and let $k \geq 0$. Then $\flatx\H \, \subseteq T(H_k)$.
\el

\proof
We have $\flatx\H \, \subseteq T(H)$ by Lemma \ref{propt}(ii). We may assume that $k \leq \dim\H$, otherwise $H = H_k$. We claim that $T(H) \subseteq T(H_k)$. Indeed, let $T \in T(H)$. Let $X \in H_{k-1} \cap 2^T$ and $p \in V \setminus T$. Since $k-1 < \dim\H$ and $T \in T(H)$, it follows that $X \cup \{ p\} \in H$. Since $|X| \leq k-1$, we get $X \cup \{ p\} \in H_k$ and so $T \in T(H_k)$. Thus $\flatx\H \, \subseteq T(H) \subseteq T(H_k)$ and so $\flatx\H \, \subseteq T(H_k)$ holds in any case. 
\qed

\bt
\label{spell}
If $k \geq 4$, then $\H = \J(T(H_k))$.
\et

\proof
We start by showing that
\beq
\label{spell1}
\flatx\H \, = T(H_k).
\eeq

The direct inclusion follows from Lemma \ref{intrun}.  

Conversely, let $Z \in T(H_k)$. Suppose that there exist edges $i \mapright{g} j \mapright{h} \ell$ in $\Gamma(Z)$ for some $g,h \in G$ but no edge $i \mapright{gh} \ell$. Then $\{ (i,g,j),(j,h,\ell),(i,gh,\ell) \} \notin H$ and so $\{ (i,g,j),(j,h,\ell),(i,gh,\ell) \} \notin H_k$. However, $\{ (i,g,j),(j,h,\ell) \}\in H_{k-1} \cap 2^Z$, yielding $Z \notin T(H_k)$. Thus we may assume that
\beq
\label{spell2}
i \mapright{g} j,\; j \mapright{h} \ell \in \Gamma(Z) \quad \mbox{implies} \quad i \mapright{gh} \ell \in \Gamma(Z)
\eeq
holds for all $i,j,\ell,g,h$. 

Suppose first that $\Gamma_0(Z)$ has parallel edges. It follows that there exist edges
$i \mapright{g} j$ and $i \mapright{h} j$ in $\Gamma(Z)$ for some distinct $i,j \in \un{n}$, with $g,h \in G$ distinct. Let $I$ be the set of all vertices in $\Gamma_0(Z)$ and let $\pi$ be the partition of $I$ defined by the connected components of $\Gamma_0(Z)$. Clearly, $Z \subseteq C(\pi)$. We claim that equality holds. Let $i',j' \in I$ be distinct elements of the same block of $\pi$ and let $g' \in G$. We must show that $(i',g',j') \in Z$. 

We consider two subcases. Assume first that $\{ i',j'\} = \{ i,j \}$. Then we may assume that $i' = i$ and $j' = j$, as well as $g' \notin \{ g,h\}$. Then 
$\{ (i,g,j),(i,h,j),(i,g',j) \} \notin H$ and so $\{ (i,g,j),(i,h,j)$, $(i,g',j) \} \notin H_k$. However, $\{ (i,g,j),(i,h,j) \}\in H_{k-1} \cap 2^Z$. Since $Z \in T(H_k)$, this yields $(i,g',j) \in Z$ as required.

Assume now that $\{ i',j'\} \neq \{ i,j \}$. 
By definition of $\pi$, we have some path
$$i' = i'_0 \mapright{g'_1} i'_1 \mapright{g'_2} \ldots \mapright{g'_s} i'_s = j'$$
in $\Gamma(Z)$. Let $h' = g'_1\ldots g'_s$. In view of (\ref{spell2}), we have $(i',h',j') \in Z$.  Then 
$\{ (i,g,j),(i,h,j),(i',g',j')$, $(i',h',j') \} \notin H$ and so $\{ (i,g,j),(i,h,j),(i',g',j'),(i',h',j') \} \notin H_k$. However, $\{ (i,g,j),(i,h,j)$, $(i',h',j') \}\in H_{k-1} \cap 2^Z$. Since $Z \in T(H_k)$, this yields $(i',g',j') \in Z$. 

Thus $Z = C(\pi)$ when $\Gamma_0(Z)$ has parallel edges and so $Z \in \flatx\H$ by Lemma \ref{theflatsT}.

Hence we may assume that $\Gamma_0(Z)$ is a graph (without parallel edges!). In view of (\ref{spell2}), $\Gamma_0(Z)$ is a union of cliques. Suppose that 
$\Gamma(Z)$ is $G$-nontrivial. We claim that $\Gamma(Z)$ has a $G$-nontrivial triangle. Indeed, let 
$$i_0 \mapright{g_1} i_1 \mapright{g_2} \ldots  \mapright{g_m} i_m = i_0$$
be a $G$-nontrivial cycle of $\Gamma(Z)$ of shortest length. Suppose that $m > 3$. By (\ref{spell2}), 
$$i_0 \mapright{g_1} i_1 \mapright{g_2} \ldots  \longmapright{g_{m-2}} i_{m-2} \vlongmapright{g_{m-1}g_m} i_m = i_0$$
is a shorter $G$-nontrivial cycle of $\Gamma(Z)$, a contradiction. Therefore $m = 3$ and so $\Gamma(Z)$ has a $G$-nontrivial triangle
$$i \mapright{g} j \mapright{h} \ell \mapright{s} i.$$
Hence $s \neq (gh)\inv$. Now $\{ (i,g,j),(j,h,\ell),(\ell,s,i) \} \in H_{k-1} \cap 2^Z$. On the other hand, since $\Gamma_0(Z)$ has no parallel edges, we have $(i,gh,\ell) \notin Z$. Since $Z \in T(H_k)$, we get $\{ (i,g,j),(j,h,\ell),(\ell,s,i)$, $(i,gh,\ell) \} \in H_{k} \subseteq H$, a contradiction. 

Therefore $\Gamma(Z)$ is $G$-trivial and so $Z \in \flatx\H$ by Lemma \ref{theflatsT}. This establishes (\ref{spell1}).

Since $\H$ is a matroid, the elements of $H$ are precisely the transversals of the successive differences for chains in $\flatx\H = T(H_k)$. Therefore $H = J(T(H_k))$ and so $\H\, = \J(T(H_k))$.
\qed

We show next that the value $k \geq 4$ in Theorem \ref{spell} is optimal.
Assume that $n \geq 3$.
Write $(1,G,2) = \{ (1,g,2) \mid g \in G\}$. It is straightforward to check that
$$\emptyset \subset \{ (1,1,2) \} \subset (1,G,2) \subset (1,G,2) \cup \{ (2,1,3) \} \subset V$$ 
is a chain in $T(H_3)$, hence $\dim\J(T(H_3)) \geq 3$. However, $\dim\H\, = 2$.

\subsection{The Dowling matroid}

Let $G$ be a finite nontrivial group and fix an integer $n \geq 2$. In this subsection
$V$ denotes the set of all $(I,\pi,[f]_{\pi}) \in \SPC(\un{n},G)$ of the following two types:
\bi
\item[(V1)]
$|I| = n-1$ and $\pi$ is the identity partition;
\item[(V2)]
$I = \un{n}$ and $\pi$ has precisely $n-1$ blocks.
\ei
Thus $V$ is the set of atoms of the Dowling lattice $Q_n(G)$. Also in this case, the elements of $V$ admit a graph-theoretical description. For this, we fix some element $y \in G \setminus \{ 1 \}$. If $(I,\pi,[f]_{\pi})$ is of type (V1) and $I = \un{n} \setminus \{ i \}$, we define $\Gamma(I,\pi,[f]_{\pi})$ to be the directed labeled graph 
$$\xymatrix{
i \ar@(ur,r)^{y}
}
$$
We shall use the simplified notation $(i,y,i)$ to denote the above $(I,\pi,[f]_{\pi})$.

If $(I,\pi,[f]_{\pi})$ is of type (V2), $\{ i,j\}$ is the unique nonsingular block of $\pi$ and $(f(i))\inv f(j) = g \in G$, we define $\Gamma(I,\pi,[f]_{\pi})$ to be the directed labeled graph 
\beq
\label{arcade1}
\xymatrix{
i \ar@/^/[rr]^{g}&& j \ar@/^/[ll]^{g\inv}
}
\eeq
Once again, it suffices to represent the above graph by the one-edge graph $i \mapright{g} j$ (or $j \mapright{g\inv} i$). To simplify notation, we shall denote by $(i,g,j)$ (or $(j,g\inv,i)$) the above $(I,\pi,[f]_{\pi})$.

Now each $Z \subseteq V$ can be represented by a directed labeled graph $\Gamma(Z)$. We denote by $\Gamma_0(Z)$ the undirected multigraph with the same vertex set as $\Gamma(Z)$ and an edge $i \edge j$ for each pair of inverse edges (\ref{arcade1}) in $\Gamma(Z)$.

We say that a multigraph $\Gamma$ is {\em unicyclic with loops} if:
\bi
\item
$\Gamma$ is connected;
\item
the number $v$ of vertices equals the number $e$ of edges in $\Gamma$.
\ei 
Since finite trees can be characterized as connected graphs satisfying $e = v-1$, it follows that unicyclic graphs are precisely those graphs which can be obtained from a tree by adding a new edge. 

The concept of $G$-trivial cycle or graph is inherited from Subsection \ref{rrm}.

Let $H$ denote the set of all $Z \subseteq V$ such that
\bi
\item
each connected component of $\Gamma_0(Z)$ is either a tree or a unicyclic graph:
\item
every unicyclic connected component of $\Gamma_0(Z)$ arises from a  $G$-nontrivial cycle of $\Gamma(Z)$.
\ei
Then $\H\, = (V,H)$ is the Dowling matroid defined by $G$ and $n$. Note that $\H$ has dimension $n-1$.

Since the Dowling lattice $Q_n(G)$ is geometric and generates $\H$, then $Q_n(G)$ is isomorphic to $\flatx\H$. So the flats of $\H$ correspond to the vertices of $Q_n(G)$, which are precisely the elements of $\SPC(n,G)$. More precisely, the flat of $\H$ determined by $(I,\pi,[f]_{\pi}) \in \SPC(n,G)$ is the set of elements of $V$ (i.e. atoms of $Q_n(G)$) lying below $(I,\pi,[f]_{\pi})$ for the Dowling order. These are elements of the following three types:
\bi
\item
$(i,y,i)$ with $i \in \un{n} \setminus I$;
\item
$(i,g,j)$ with $i,j \in \un{n} \setminus I$ distinct and $g \in G$;
\item
$(i,g,j)$ with $i,j \in I$ distinct in the same block of $\pi$ and $g = (f(i))\inv f(j)$.
\ei
In graph theoretical terms, this amounts to say that there is (at most) one connected component in $\Gamma_0(Z)$ which contains all possible edges (one loop at each vertex plus $|G|$ parallel edges between any pair of distinct vertices), while the remaining connected components are complete graphs (no loops, no parallel edges) and arise from $G$-trivial components of $\Gamma(Z)$. 

\bt
\label{dowt}
If $k \geq 3$, then $\H = \J(T(H_k))$.
\et

\proof
Similarly to the proof of Theorem \ref{spell}, it suffices to show that
\beq
\label{dowt1}
\flatx\H \, = T(H_k).
\eeq

The direct inclusion follows from Lemma \ref{intrun}.  

Conversely, let $Z \in T(H_k)$. Suppose that there exist edges $i \mapright{g} j \mapright{h} \ell$ in $\Gamma(Z)$ for some $g,h \in G$ but no edge $i \mapright{gh} \ell$. Then $\{ (i,g,j),(j,h,\ell),(i,gh,\ell) \} \notin H$ and so $\{ (i,g,j),(j,h,\ell),(i,gh,\ell) \} \notin H_k$. However, $\{ (i,g,j),(j,h,\ell) \}\in H_{k-1} \cap 2^Z$, yielding $Z \notin T(H_k)$. Thus we may assume that
\beq
\label{dowt2}
i \mapright{g} j,\; j \mapright{h} \ell \in \Gamma(Z) \quad \mbox{implies} \quad i \mapright{gh} \ell \in \Gamma(Z)
\eeq
holds for all $i,j,\ell,g,h$. 

Suppose now that some connected component of $\Gamma_0(Z)$ arises from a $G$-nontrivial component $C$ of $\Gamma(Z)$. We claim that this same component must contain a loop. Indeed, let 
$$i_0 \mapright{g_1} i_1 \mapright{g_2} \ldots  \mapright{g_m} i_m = i_0$$
be a $G$-nontrivial cycle of $\Gamma(Z)$ of shortest length. Suppose that $m > 2$. By (\ref{dowt2}), 
$$i_0 \mapright{g_1} i_1 \mapright{g_2} \ldots  \longmapright{g_{m-2}} i_{m-2} \vlongmapright{g_{m-1}g_m} i_m = i_0$$
is a shorter $G$-nontrivial cycle of $C$, a contradiction. Therefore $m \leq 2$. Suppose that $m = 2$. Then $C$ contains edges
$$
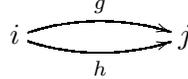
\xymatrix{
i \ar@/^/[rr]^{g} \ar@/_/[rr]_{h} && j
}$$
for some distinct $g,h \in G$. Note that $\{ (i,g,j), (i,h,j) \} \in H_{k-1} \cap Z$ but $\{ (i,g,j), (i,h,j), (i,y,i) \} \notin H$ (and therefore $\notin H_k$). Thus $(i,y,i)$ and so the component $C$ must contain some loop. 

Next we show that at most one connected component of $\Gamma_0(Z)$ can arise from a $G$-nontrivial component of $\Gamma(Z)$. Indeed, suppose there exist two such components. By the preceding claim, each of these two components contains a loop. Let $(i,y,i), (j,y,j)$ correspond to loops in different components. Then $\{ (i,y,i), (j,y,j) \} \in H_{k-1} \cap 2^{Z}$ but $(i,1,j) \notin Z$. Since $Z \in T(H_k)$, we get $\{ (i,y,i), (j,y,j), (i,1,j) \} \in H$, a contradiction. Thus at most one connected component of $\Gamma_0(Z)$ can arise from a $G$-nontrivial component of $\Gamma(Z)$.

Note that in this $G$-nontrivial component $C$ each pair of vertices are connected by some edge in view of (\ref{dowt2}). Suppose that $(i,y,i)$ produces a loop in $C$ and $j$ is another vertex of $C$. Then $(i,g,j) \in Z$ for some $g \in G$. Then $\{ (i,y,i), (i,g,j) \} \in H_{k-1} \cap 2^{Z}$ but $\{ (i,y,i), (i,g,j), (j,y,j) \} \notin H$. Since $Z \in T(H_k)$, it follows that $(j,y,j) \in Z$. Therefore $C$ has loops at every vertex. 

Finally, suppose that $i,j$ are distinct vertices of $C$ and $g \in G$. Then $\{ (i,y,i), (j,y,j) \} \in H_{k-1} \cap 2^{Z}$ but $\{ (i,y,i), (i,g,j), (j,y,j) \} \notin H$. Since $Z \in T(H_k)$, it follows that $(i,g,j) \in Z$. Therefore the unique $G$-nontrivial component of $\Gamma(Z)$
contains all possible edges (one loop at each vertex plus $|G|$ parallel edges between any pair of distinct vertices), while the remaining connected components are complete graphs (in view of (\ref{dowt2})) and arise from $G$-trivial components of $\Gamma(Z)$. This implies that $Z \in \flatx\H$ and so (\ref{dowt1})) holds as required.
\qed

\section{Shellability for high dimensions}

In this section, we discuss shellablity in paving boolean representable simplicial complexes.

Let $\H\, = (V,H) \in \bpav(d)$ for $d \geq 2$. We say that $F \in \flatx\H$ is a {\em line} of $\H$ if $d \leq |F| < |V|$. Let $\L_{\H}$ denote the set of lines of $\H$. Given $L \in \L_{\H}$, write
$$L\mu = \{ I \cup \{ p \} \mid I \in P_d(L),\; p \in V \setminus L \}.$$

\bl
\label{big}
Let $\H\, = (V,H) \in {\rm BPav}(d)$. Then:
\bi
\item[(i)] ${\rm Fl}\H\, = P_{\leq d-1}(V) \cup \L_{\H} \cup \{ V\}$;
\item[(ii)] if $L,L' \in \L_{\H}$ are distinct, then $|L \cap L'| \leq d-1$;
\item[(iii)] $H = P_{\leq d}(V) \cup (\displaystyle\bigcup_{L \in \L_{\H}} L\mu)$.
\ei
\el

\proof
(i) By \cite[Lemma 6.1.1]{RSm}.

(ii) Suppose that $|L \cap L'| \geq d$. We may assume that $L \not\subseteq L'$ and take distinct elements $a_1,\ldots,a_{d-1} \in L \cap L'$. Since $\flatx\H$ is closed under intersection by \cite[Proposition 4.2.2(ii)]{RSm}, then
$$\emptyset \subset a_1 \subset a_1a_2 \subset \ldots \subset a_1\ldots a_{d-1} \subset L \cap L' \subset L \subset V$$
is a chain in $\flatx\H$, which admits a transversal of the successive differences with $d+2$ elements, contradicting $\H\,\in \bpav(d)$. 

(iii) Since $\H \in \bpav(d)$, we have $P_{\leq d}(V) \subset H \subseteq P_{\leq d+1}(V)$. Now it follows from the definition of flat that 
$\displaystyle\bigcup_{L \in \L_{\H}} L\mu \subseteq H$. 

Finally, let $I \in H \cap P_{d+1}(V)$. Then $I$ is a transversal of the successive differences for some chain 
$$F_0 \subset F_1 \subset \ldots \subset F_{d+1}$$
in $\flatx\H$. Since this chain is necessarily maximal, we have $F_{d+1} = V$, and it follows from part (i) that $F_d \in \L_{\H}$. Hence $I \in F_d\mu$ as required.
\qed

We denote by $\fct_i\H$ the set of facets of $\H$ of dimension $i$. Since $\H \in \bpav(d)$, we have $\fct\H\, = \fct_{d}\H \, \cup \fct_{d-1}\H$, and it follows from Lemma \ref{big}(iii) that 
\beq
\label{big1}
\fct_d\H\, = \displaystyle\bigcup_{L \in \L_{\H}} L\mu.
\eeq

Now define a simplicial complex $\H^* = (V^*,H^*)$ of dimension $d-1$ by
$$V^* = \bigcup \L_{\H}, \quad H^* = \bigcup_{L \in \L_{\H}} P_{\leq d}(L).$$
We tried to determine if there is any relationship between shellability of $\H \in \bpav(d)$ and shellability of $\H^*$. The following two examples show that neither of the two implications holds in general.

\be
\label{boom}
Let $V = \un{6}$ and let $H = P_{\leq 3}(V) \cup \{ X \in P_4(V) \mid \mbox{$X$ contains three consecutive numbers }\}$. Then $\H \in {\rm BPav}(3)$ is not shellable but $\H^*$ it is so.
\ee

Indeed, it is straightforward to check that 
$$\flatx\H\, = P_{\leq 2}(V) \cup \{ 123, 234, 345, 456, V\},$$
and it follows easily that $\H$ is a BRSC. Moreover, $\L_{\H} = \{ 123, 234, 345, 456 \}$ and $\H^* = (V,H^*)$ with
$$H^* = P_{\leq 1}(V) \cup \{ 12,13,23,24,34,35,45,46,56 \} \cup \L_{\H}.$$
Thus $\fct\H^* = \L_{\H}$ and $123, 234, 345, 456$ constitutes a shelling of $\H^*$.

Since $16 \in H$, we can consider the contraction 
$$\H/16 = (2345, P_{\leq 1}(2345) \cup \{ 23,45\}).$$
It is easy to check that a graph is shellable if and only if it has at most one nontrivial connected component. Hence $\H/16$ is not shellable. Since the class of shellable simplicial complexes is closed under contraction
\cite[Proposition 7.1.5]{RSm}, it follows that $\H$ is not shellable either.

\be
\label{tracks}
Let $V = \un{7}$ and let $H = P_{\leq 2}(V) \cup \{ X \in P_3(V) \mid \mbox{$X$ contains 12, 23, 45, 56 or 67 }\}$. Then $\H \in {\rm BPav}(2)$ is shellable but $\H^*$ is not.
\ee

Indeed, it is straightforward to check that 
$$\flatx\H\, = P_{\leq 1}(V) \cup \{ 12, 23, 45, 45, 56, 67, V\},$$
and it follows easily that $\H$ is a BRSC. Moreover, $\L_{\H} = \{ 12, 23, 45, 45, 56, 67 \}$ and $\H^* = (V,H^*)$ with
$$H^* = P_{\leq 1}(V) \cup \L_{\H}.$$
Since the graph $\H^*$ has two nontrivial connected components, then $\H^*$ is not shellable.

Now $\H^*$ coincides with the {\em graph of flats} defined in \cite[Section 6.4]{RSm}. Therefore $\H$ is shellable by \cite[Theorem 7.2.8]{RSm}.

\medskip

The best we can do is to show that the implications just disproved hold in some particular cases:

\bp
\label{rega}
If $\H \in {\rm BPav}(2)$ and $\H^*$ is shellable, so is $\H$.
\ep

\proof
We have already remarked that $\H^*$ is the graph of flats of $\H$, which has at least an edge since $\H^*$ has dimension 1. Thus $\H^*$ shellable is equivalent to say that the graph of flats of $\H$ has a single nontrivial connected component. Therefore $\H$ is shellable by \cite[Theorem 7.2.8]{RSm}.
\qed

\bp
\label{iso}
Let $\H \in {\rm BPav}(d)$ be shellable. If $V^* \subset V$, then $\H^*$ is shellable.
\ep

\proof
Let $z \in V \setminus V^*$. Since shellability is preserved under contraction, $\H/z$ is shellable as well. Clearly, 
$$\fct\H^* = \bigcup_{L \in \L_{\H}} P_{d}(L).$$
We claim that
\beq
\label{iso1}
\fct \H^* = \fct_{d-1} \H/z.
\eeq
Let  $X \in \fct \H^*$. Then $X \in P_d(L)$ for some $L \in \L_{\H}$. Since $z \in V \setminus V^*$, we have $X \cup \{ z \} \in L\mu$ and so $X \cup \{ z \} \in H$ by Lemma \ref{big}(iii). Thus $X \cup H/z$. Since $|X| = d$ and $\H/z$ has dimension $d-1$, we get $X \in \fct_{d-1} \H/z$.

Conversely, let $X \in \fct_{d-1} \H/z$. Then $X \cup \{ z \} \in H$ and $|X| = d$. Hence there exists an enumeration $x_1,\ldots,x_{d+1}$ of the elements of $X \cup \{ z \}$ and a chain
$$F_0 \subset F_1 \subset \ldots \subset F_{d+1}$$
in $\flatx\H$ such that $x_i \in F_i \setminus F_{i-1}$ for each $i \in \un{d+1}$. In view of Lemma \ref{big}(i), we must have $F_{d+1} = V$, $F_d \in \L_{\H}$ and $|F_i| = i$ for $i = 0,\ldots,d-1$. Since $z \notin F_d$, we must have $z = x_{d+1}$ and so $X \subseteq F_d$. Thus $X \in H^*$. Since $|X| = d$, we get $X \in \fct \H^*$. Therefore (\ref{iso1}) holds.

Now $\H/z$ is shellable. By \cite{BW}, $\H/z$ admits a shelling where the dimension of the facets is not increasing. In view (\ref{iso1}), $\H/z$ admits a shelling which starts by an enumeration of the facets of $\H^*$. Therefore $\H^*$ is itself shellable.
\qed

\section{Going up (or not) for paving complexes}

Let $d \geq 1$ and let $\H\, = (V,H) \in \pav(d)$. Since $P_{\leq d-1}(V) \cup \{ V \} \subseteq \flatx\H \subseteq T(H)$ by Lemma \ref{propt}(ii), we have $\dim\J(T(H)) \geq d-1$. 

Note that equality may occur, take for instance $V = 1234$ and $H = P_{\leq 2}(V) \cup \{ 123 \}$, when $T(H) = P_{\leq 1}(V) \cup \{ V \}$. However, if $\H \in \tbpav(d)$, then $\dim\J(T(H)) \geq d$ because $H \subseteq J(T(H))$ by Theorem \ref{eqtr}(ii). We say that $\H \in \pav(d)$ {\em goes up} if $\dim\J(T(H)) > d$. Otherwise, we say that $\H$ {\em does not go up}. We denote by $\gu(d)$ (respectively $\ngu(d)$) the class of all $\H \in \pav(d)$ which go up (respectively, do not go up).

\bl
\label{orp}
Let $(V,H), (V,H') \in {\rm Pav}(d)$ with $H \subseteq H'$. Then:
\bi
\item[(i)] $T(H) \subseteq T(H')$;
\item[(ii)] $J(T(H)) \subseteq J(T(H'))$.
\ei
\el

\proof
(i) Let $T \in T(H)$. Let $X \in H' \cap P_{\leq d}(T)$ and $p \in V \setminus T$. Since $H' \cap P_{\leq d}(T) = P_{\leq d}(T) = H \cap P_{\leq d}(T)$ and $T \in T(H)$, we get $X \cup \{ p\} \in H$ and so $X \cup \{ p\} \in H'$. Thus $T \in T(H')$. 

(ii) Immediate.
\qed

\bc
\label{gungu}
Let $(V,H), (V,H') \in {\rm Pav}(d)$ with $H \subseteq H'$.
\bi
\item[(i)] If $(V,H)$ goes up, so does $(V,H')$.
\item[(ii)] If $(V,H')$ does not go up, neither does $(V,H)$.
\ei
\ec

\proof
(i) If $(V,H)$ goes up, then there exists a chain
\beq
\label{gungu1}
T_0 \subset T_1 \subset \ldots \subset T_{d+2}
\eeq
in $T(H)$. By Lemma \ref{orp}(i), this is also a chain in $T(H')$. Therefore $(V,H')$ goes up.

(ii) By part (i).
\qed

In view of Corollary \ref{gungu}, it is only natural to define the following two classes of complexes:
\bi
\item
If $(V,H) \in \gu(d)$ and $(V,H') \in \ngu(d)$ whenever $H \supset H' \supset P_{\leq d}(V)$, we say that $(V,H)$ is {\em minimal going up} and we write $(V,H) \in \mgu(d)$.
\item
If $(V,H) \in \ngu(d)$ and $(V,H') \in \gu(d)$ whenever $H \subset H' \subseteq P_{\leq d+1}(V)$, we say that $(V,H)$ is {\em minimal not going up} and we write $(V,H) \in \mngu(d)$.
\ei

Our main focus in this context is the class $\mngu(d)$.

Given $(V,H) \in \pav(d)$, we define a closure operator $\cl_T: 2^V \to T(H)$ by 
$$\cl_T(X) = \bigcap\{ T \in T(H) \mid X \subseteq T\}.$$
Note that $\cl_T$ is well defined since $V \in T(H)$ and $T(H)$ is closed under intersection by Lemma \ref{propt}(i).

\bl
\label{cltt}
Given $(V,H) \in {\rm Pav}(d)$, the following conditions are equivalent:
\bi
\item[(i)]
$(V,H)$ goes up;
\item[(ii)]
There exist $X \in P_{d+1}(V)$ and $Y \in P_d(X)$ such that ${\rm Cl}_T(Y) \subset {\rm Cl}_T(X) \subset V$.
\ei
\el

\proof
(i) $\Rw$ (ii). If $(V,H)$ goes up, there exists a chain of the form (\ref{gungu1}) in $T(H)$. Since $P_{\leq d-1}(V) \cup \{ V\} \subseteq \flatx\H \subseteq T(H)$, we may assume that (\ref{gungu1}) is of the form
$$\emptyset \subset a_1 \subset a_1a_2 \subset \ldots \subset a_1\ldots a_{d-1} \subset T_d \subset T_{d+1} \subset V.$$
If we pick $a_d \in T_d \setminus a_1\ldots a_{d-1}$, then $a_1\ldots a_{d-1} \subset \cl_T(a_1\ldots a_d) \subseteq T_d$. Choosing $a_{d+1} \in T_{d+1} \setminus T_d$, we get $\cl_T(a_1\ldots a_d) \subset \cl_T(a_1\ldots a_{d+1}) \subseteq T_{d+1} \subset V$, so condition (ii) is satisfied by $X = a_1\ldots a_{d+1}$ and $Y = a_1\ldots a_{d}$.

(ii) $\Rw$ (i). If $a_1,\ldots, a_{d-1}$ are distinct elements of $Y$, then we obtain a chain
$$\emptyset \subset a_1 \subset a_1a_2 \subset \ldots \subset a_1\ldots a_{d-1} \subset {\rm Cl}_T(Y) \subset {\rm Cl}_T(X) \subset V$$
in $T(H)$ and $(V,H)$ goes up.
\qed

\bl
\label{inco}
If $(V,H) \in {\rm MNGU}(d)$ and $W \in P_{d+2}(V)$, then $H \cap P_{d+1}(W) \neq \emptyset$.
\el

\proof
Suppose that $H \cap P_{d+1}(W) = \emptyset$. Fix $W' \in P_{d+1}(W)$ and let $H' = H \cup \{ W'\}$. Then $(V,H') \in {\rm Pav}(d)$ and $T(H) \subseteq T(H')$ by Lemma \ref{orp}(i). Suppose that $F \in T(H')$. Let $X \in H \cap P_{\leq d}(F)$ and $p \in V \setminus F$. Since $H \subseteq H'$ and  $F \in T(H')$, we get $X \cup \{ p \} \in H'$. Suppose that $X \cup \{ p \} = W'$. Then $W = X \cup \{ p,q \}$ for some other vertex $q$. Now $F \in T(H')$, $X \in H' \cap P_{\leq d}(F)$ and $X \cup \{ q \} \notin H'$ together yield $q \in F$. Take now $r \in W \setminus \{ p,q\}$. Then $F \in T(H')$, $W\setminus \{ r,p\} \in H' \cap P_{\leq d}(F)$ and $W\setminus \{ r \} \notin H'$ together yield $p \in F$, we reach a contradiction. Hence $X \cup \{ p \} \in H$ and so $T(H') = T(H)$. Since $\dim\J(T(H)) \leq \dim(V,H) = \dim(V,H')$, it follows that $(V,H') \in {\rm NGU}(d)$, contradicting $(V,H) \in {\rm MNGU}(d)$. Therefore $H \cap P_{d+1}(W) \neq \emptyset$.
\qed

Given $\H\, = (V,H) \in \pav(d)$, we define the {\em defect} of $\H$ through
$$\defe\H = P_{d+1}(V) \setminus H.$$
If $d = 1$, we can view $\defe\H$ as the set of edges of a graph with vertex set $V$, denoted by $\Defe\H$.

\bp
\label{mnguone}
The following conditions are equivalent for $\H \,\in {\rm Pav}(1)$:
\bi
\item[(i)] $\H\,\in {\rm GU}(1)$;
\item[(ii)] ${\rm Def}\H$ has more than two connected components.
\ei
\ep

\proof
Write $\H\, = (V,H)$. we start by showing that 
\beq
\label{mnguone1}
T(H)\mbox{ consists of unions of connected components of Def}\H.
\eeq

Indeed, let $C \subseteq V$ be a union of connected components of Def$\H$. Suppose that $X \in H \cap P_{\leq 1}(C)$ and $p \in V setminus C$. Since $P_1(V) \subseteq H$, we may assume that $X = \{ x \}$. If $xp \notin H$, then $xp$ would be an edge of $\Defe\H$, contradicting $p \notin C$. Hence $X \cup \{ p \} \in H$ and so $C \in T(H)$. 

Suppose now that $T \in T(H)$ and there exists some edge $pq$ in $\Defe\H$ with $q \in T$. Since $pq \notin H$ and $q \in H \cap P_{\leq 1}(T)$, it follows from $T \in T(H)$ that $p \in T$ as well. Therefore $T$ is a union of connected components of Def$\H$ and (\ref{mnguone1}) holds.

But $\H\,\in {\rm GU}(1)$ if and only if there exists some chain of the form
$$T_0 \subset T_1 \subset T_2 \subset T_3$$
in $T(H)$. In view of (\ref{mnguone1}), this is possible if and only if ${\rm Def}\H$ has at least three connected components. If $C_1, C_2,C_3$ are three distinct connected components, then $\emptyset \subset C_1 \subset C_1 \cup C_2 \subset C_1 \cup C_2 \cup C_3$ is the desired chain
in $T(H)$.
\qed

\bc
\label{mnone}
The following conditions are equivalent for $\H \,\in {\rm Pav}(1)$:
\bi
\item[(i)] $\H\,\in {\rm MNGU}(1)$;
\item[(ii)] ${\rm Def}\H$ is a forest with exactly two connected components.
\ei
\ec

\proof
By Proposition \ref{mnguone}, $\H\,\in {\rm NGU}(1)$ if and only if ${\rm Def}\H$ has at most two connected components. 

Suppose that ${\rm Def}\H$ is connected. Let $pq$ be an edge of ${\rm Def}\H$ and let $H' = H \cup \{ pq \}$. Then ${\rm Def}(V,H')$ is obtained by removing the edge $pq$ from $\Defe\H$, hence ${\rm Def}(V,H')$ has at most two connected components. It follows that $(V,H') \in {\rm NGU}(1)$ and so $\H\,\notin {\rm NGU}(1)$.

Suppose now that ${\rm Def}\H$ contains a cycle. Then we can remove an edge $pq$ from ${\rm Def}\H$ without changing the number of connected components. This corresponds to adding $pq$ to $H$ to obtain $(V,H \cup \{ pq \}) \in {\rm NGU}(1)$. Thus $\H\,\notin {\rm NGU}(1)$.

Therefore (i) implies (ii). Assume now condition (ii).
If we remove an edge from a forest with two connected components, we obtain three connected components. Hence any complex $(V,H')$ such that $H' \subset H$ is in ${\rm GU}(1)$ by Proposition \ref{mnguone}. Thus $\H\,\in {\rm MNGU}(1)$ as required.
\qed

\bc
\label{mone}
The following conditions are equivalent for $\H \,\in {\rm Pav}(1)$:
\bi
\item[(i)] $\H\,\in {\rm mGU}(1)$;
\item[(ii)] ${\rm Def}\H$ is a union of three disjoint cliques.
\ei
\ec

\proof
We can adapt the arguments used in the proof of Corollary \ref{mnone}. Indeed, it follows from
Proposition \ref{mnguone} that having at least three connected components in ${\rm Def}\H$ is a necessary condition for $\H\,\in {\rm mGU}(1)$. Now we are looking for graphs such that, by adding any edge, we obtain a graph with less than three connected components. This happens precisely when ${\rm Def}\H$ is a union of three disjoint cliques.
\qed

We compute next the complexes in $\mngu(2)$ with 4, 5 and 6 vertices.

\bp
\label{4mngu}
Up to isomorphism, there is only one $\H\, = (1234,H) \in {\rm MNGU}(2)$, which can be defined by ${\rm def}\H\, = \{ 123 \}$.
\ep

\proof
Suppose that $\H\, = (1234,H) \in \gu(2)$. By Lemma \ref{cltt}, there exist distinct $a,b,c \in 1234$ such that $ab, abc \in T(H)$. We may assume that $12,123 \in T(H)$. Since $123 \in T(H)$, 4 cannot occur in $\defe\H$. But $12 \in T(H)$, so $123 \in H$ and so $\defe\H\, = \emptyset$. Now the claim follows.
\qed

\bp
\label{5mngu}
For five vertices, there exist two isomorphism classes in ${\rm MNGU}(2)$, with representatives
$\H_i\, = (12345,H_i)$ $(i = 1,2)$ defined by ${\rm def}\H_1 = \{ 123, 124, 134 \}$ and ${\rm def}\H_2 = \{ 123, 345 \}$.
\ep

\proof
Suppose that $\H_1\,\in \gu(2)$. By Lemma \ref{cltt}, there exist distinct $a,b,c \in 12345$ such that $\cl_T(ab) \subset \cl_T(abc) \subset 12345$. But $|abc \cap 1234| \geq 2$, and it is easy to see that $\cl_T(X) = 1234$ for every $X \in P_2(1234)$. Hence $T(abc) = 1234$, yielding $T(ab) = 1234$ as well, a contradiction. Thus $\H_1\,\in \ngu(2)$.

Suppose now that we add a triangle to $H_1$. Without loss of generality, we may consider the complex $\H' = (12345,H')$ defined by $\defe\H' = \{ 123, 124, \}$. But then we get a chain
$$\emptyset \subset 3 \subset 34 \subset 345 \subset V$$
and so $\H'$ goes up. Thus $\H_1\,\in \mngu(2)$.

Suppose that $\H_2\,\in \gu(2)$. By Lemma \ref{cltt}, there exist distinct $a,b,c \in 12345$ such that $\cl_T(ab) \subset \cl_T(abc) \subset 12345$. If $abc \in \{ 123, 345\}$, then $T(abc) = abc = T(ab)$. Otherwise, we may assume out of symmetry that $|abc \cap 123| = 2$ and it follows easily that $T(abc) = 12345$, in any case a contradiction. Thus $\H_2\,\in \ngu(2)$.

Suppose now that we add a triangle to $H_2$. Without loss of generality, we may consider the complex $\H'' = (12345,H'')$ defined by $\defe\H'' = \{ 123 \}$. Since $\H''$ contains $\H_1 \in \mngu(2)$ as a proper subcomplex, it follows that $\H''$ goes up. Thus $\H_2\,\in \mngu(2)$. 

It is obvious that $\H_1$ and $\H_2$ are not isomorphic. Now let $(12345,H) \in \mngu(2)$ be arbitrary. 

If $|\defe\H| \leq 2$, then either $\H \cong \H_2$ or contains a proper subcomplex isomorphic to either $\H_1$ or $\H_2$. Thus we may assume that $|\defe\H| \geq 3$. We may also assume that 
\beq
\label{5mngu1}
\mbox{$|abc \cap def| \geq 2$ for any distinct $abc, def \in \defe\H$}
\eeq
otherwise $\H$ contains a proper subcomplex isomorphic to $\H_2$. 

We may now assume that $123, 124 \in \defe\H$. Let $xyz$ be a third element of $\defe\H$. If $5 \notin xyz$, then $\H$ is isomorphic to a subcomplex of $\H_1$, hence $\H \cong \H_1$. Thus, in view of (\ref{5mngu1}), we may assume that $125 \in \defe\H$ as well. 

If $\defe\H = \{ 123, 124, 125 \}$, then $\cl_T(34) = 34 \subset 345 = \cl_T(345)$, hence $\H$ goes up by Lemma \ref{cltt}. Thus we may assume that $\defe\H$ contains a fourth element $uvw$. But then $|uvw \cap 12k| = 1$ for some $k \in 345$ and so $\H$ is isomorphic to a proper subcomplex of $\H_2$. Therefore $\H$ is isomorphic to either $\H_1$ or $\H_2$.
\qed

\bp
\label{6mngu}
For six vertices, there exist ten isomorphism classes in ${\rm MNGU}(2)$, with representatives
$\M_i\, = (123456,M_i)$ $(i = 1,\ldots,10)$ defined by:
\bi
\item
${\rm def}\M_1 = \{ 124, 134, 234, 356 \}$,
\item
${\rm def}\M_2 = \{ 124, 134, 234, 456 \}$, 
\item
${\rm def}\M_3 = \{ 124, 134, 234, 135, 245 \}$, 
\item
${\rm def}\M_4 = \{ 124, 134, 234, 145, 245, 345 \}$, 
\item
${\rm def}\M_5 = \{ 123, 134, 125, 346 \}$, 
\item
${\rm def}\M_6 = \{ 123, 134, 256, 346 \}$, 
\item
${\rm def}\M_7 = \{ 123, 134, 245, 356 \}$, 
\item
${\rm def}\M_8 = \{ 123, 134, 235, 346, 356\}$, 
\item
${\rm def}\M_{9} = \{ 123, 134, 145, 235, 245 \}$. 
\item
${\rm def}\M_{10} = \{ 123, 146, 245, 356 \}$. 

\ei
\ep

\proof
Write $V = 123456$. If $V = a_1\ldots a_6$ and $\H = (V,H)$, we denote by $\H(a_1,\ldots,a_6)$ the complex obtained from $\H$ by replacing $j$ by $a_j$ for every $j \in V$. 

We start showing that $\M_i \in \mngu(2)$ for $i = 1,\ldots, 10$.

\smallskip

\noindent
$\M_1$: Suppose that there exist distinct $a,b,c \in V$ such that $\cl_T(ab) \subset \cl_T(abc) \subset V$. If $abc \subset 1234$, then $\cl_T(ab) = 1234 = \cl_T(abc)$. If $abc = 356$, then $\cl_T(ab) = 345 = \cl_T(abc)$. In any other case, it is easy to deduce first that $3 \in \cl_T(abc)$ and then $\cl_T(abc) = V$. Therefore we cannot have $\cl_T(ab) \subset \cl_T(abc) \subset V$ and so $\M_1 \in \ngu(2)$ by Lemma \ref{cltt}.

To show that $\M_1 \in \mngu(2)$, we must show that condition (ii) in Lemma \ref{cltt} holds when we remove one element from $\defe\M_1$:
\bi
\item
If we remove $356$, we get $\cl_T(15) = 15 \subset 156 = \cl_T(156)$.
\item
If we remove $124$, we get $\cl_T(15) = 15 \subset 125 = \cl_T(125)$.
\item
If we remove $134$, we get $\cl_T(15) = 15 \subset 1356 = \cl_T(156)$.
\ei
The remaining case follows by symmetry, thus $\M_1 \in \mngu(2)$.

\smallskip

\noindent
$\M_2$: Suppose that there exist distinct $a,b,c \in V$ such that $\cl_T(ab) \subset \cl_T(abc) \subset V$. If $abc \subset 1234$, then $\cl_T(ab) = 1234 = \cl_T(abc)$. If $abc = 356$, then $\cl_T(ab) = 345 = \cl_T(abc)$. In any other case, it is easy to deduce first that $4 \in \cl_T(abc)$ and then $\cl_T(abc) = V$. Therefore we cannot have $\cl_T(ab) \subset \cl_T(abc) \subset V$ and so $\M_2 \in \ngu(2)$ by Lemma \ref{cltt}.

To show that $\M_2 \in \mngu(2)$, we must show that condition (ii) in Lemma \ref{cltt} holds when we remove one element from $\defe\M_2$:
\bi
\item
If we remove $356$, we get $\cl_T(15) = 15 \subset 156 = \cl_T(156)$.
\item
If we remove $124$, we get $\cl_T(15) = 15 \subset 125 = \cl_T(125)$.
\ei
The remaining case follows by symmetry, thus $\M_2 \in \mngu(2)$.

\smallskip

\noindent
$\M_3$: If $X \in P_3(V)$, then $|X \cap 12345| \geq 2$. Since every $Y \in P_2(12345)$ is contained in some $Z \in \defe\M_3$, it is easy to check that $\cl_T(Y) = V$ and therefore $\cl_T(X) = V$. Thus $\M_3 \in \ngu(2)$ by Lemma \ref{cltt}.

To show that $\M_2 \in \mngu(2)$, we must show that condition (ii) in Lemma \ref{cltt} holds when we remove one element from $\defe\M_3$:
\bi
\item
If we remove $134$, we get $\cl_T(15) = 135 \subset 1356 = \cl_T(1356)$.
\item
Whichever other element we remove, we end up with some distinct $a,b \in 12345$ such that $\cl_T(ab) = ab \subset ab6 = \cl_T(ab6)$.
\ei
Thus $\M_3 \in \mngu(2)$.

\smallskip

\noindent
$\M_4$: It is easy to see that $\cl_T(X) = 12345$ for every $X \in P_2(12345)$. Suppose that there exist distinct $a,b,c \in V$ such that $\cl_T(ab) \subset \cl_T(abc) \subset V$. Then we may assume that $6 \in ab$, but then $\cl_T(abc) = V$. 
Thus $\M_4 \in \ngu(2)$ by Lemma \ref{cltt}.

To show that $\M_4 \in \mngu(2)$, we must show that condition (ii) in Lemma \ref{cltt} holds when we remove one element from $\defe\M_4$:
\bi
\item
If we remove $124$, we get $\cl_T(12) = 12 \subset 126 = \cl_T(126)$.
\item
If we remove $145$, we get $\cl_T(15) = 15 \subset 156 = \cl_T(156)$.
\ei
The remaining cases follows by symmetry, thus $\M_4 \in \mngu(2)$.

\smallskip

\noindent
$\M_5$: Suppose that $X \in P_3(V)$. It is easy to check that $|X \cap Y| \geq 2$ for some $Y \in \defe\M_5$, and this eventually yields $\cl_T(X) = V$. Therefore $\M_5 \in \ngu(2)$ by Lemma \ref{cltt}.

To show that $\M_5 \in \mngu(2)$, we must show that condition (ii) in Lemma \ref{cltt} holds when we remove one element from $\defe\M_5$:
\bi
\item
If we remove $125$, we get $\cl_T(15) = 15 \subset 156 = \cl_T(156)$.
\item
If we remove $123$, we get $\cl_T(15) = 125 \subset 1256 = \cl_T(156)$.
\ei
The remaining cases follows by symmetry, thus $\M_5 \in \mngu(2)$.

\smallskip

\noindent
$\M_6$: Suppose that there exist distinct $a,b,c \in V$ such that $\cl_T(ab) \subset \cl_T(abc) \subset V$. If $abc \subset 12346$, then $\cl_T(abc) = V$ in any case. If $5 \in abc$ and $abc \cap 134 \neq \emptyset$,  we also get $\cl_T(abc) = V$. So it remains the case $abc = 256$, but then $\cl_T(ab) = 256 =\cl_T(abc)$. Thus $\M_6 \in \ngu(2)$ by Lemma \ref{cltt}.

To show that $\M_6 \in \mngu(2)$, we must show that condition (ii) in Lemma \ref{cltt} holds when we remove one element from $\defe\M_6$:
\bi
\item
If we remove $123$, we get $\cl_T(15) = 15 \subset 1256 = \cl_T(156)$.
\item
If we remove $134$, we get $\cl_T(15) = 15 \subset 145 = \cl_T(145)$.
\item
If we remove $256$, we get $\cl_T(15) = 15 \subset 156 = \cl_T(156)$.
\ei
The remaining case follows by symmetry, thus $\M_6 \in \mngu(2)$.

\smallskip

\noindent
$\M_7$: Suppose that there exist distinct $a,b,c \in V$ such that $\cl_T(ab) \subset \cl_T(abc) \subset V$. Since $\cl_T(12) = \cl_T(13) = \cl_T(14) = \cl_T(23) = \cl_T(34) = V$, we may assume that $abc = 245$ or 246 or $i56$ with $i \in 1234$. If $abc = 245$, then $\cl_T(ab) = 245 = \cl_T(abc)$. If $abc = 246$, then we successively deduce $5,3 \in \cl_T(abc)$ and therefore $\cl_T(abc) = V$. If $abc = 356$, then $\cl_T(ab) = 356 = \cl_T(abc)$. Finally, if $abc = i56$ with $i \in 124$, we get $3 \in \cl_T(abc)$ and consequently $\cl_T(abc) = V$. Thus $\M_7 \in \ngu(2)$ by Lemma \ref{cltt}.

To show that $\M_7 \in \mngu(2)$, we must show that condition (ii) in Lemma \ref{cltt} holds when we remove one element from $\defe\M_7$:
\bi
\item
If we remove $123$, we get $\cl_T(12) = 12 \subset 126 = \cl_T(126)$.
\item
If we remove $245$, we get $\cl_T(25) = 25 \subset 245 = \cl_T(245)$.
\item
If we remove $356$, we get $\cl_T(15) = 15 \subset 156 = \cl_T(156)$.
\ei
The remaining case follows by symmetry, thus $\M_7 \in \mngu(2)$.

\smallskip

\noindent
$\M_8$: Suppose that $X \in P_3(V)$. It is easy to check that $|X \cap Y| \geq 2$ for some $Y \in \defe\M_5$, and this eventually yields $\cl_T(X) = V$. Therefore $\M_8 \in \ngu(2)$ by Lemma \ref{cltt}.

To show that $\M_8 \in \mngu(2)$, we must show that condition (ii) in Lemma \ref{cltt} holds when we remove one element from $\defe\M_8$. If we remove $123$, we get $\cl_T(12) = 12 \subset 126 = \cl_T(126)$.
Now we can picture $\defe\H$ building a planar graph from a cycle $125641$ and a central vertex 3 of degree 5. Hence the remaining cases follow by symmetry, and so $\M_8 \in \mngu(2)$.

\smallskip

\noindent
$\M_9$: It is easy to check that $\cl_T(X) = 12345$ for every $X \in P_2(12345)$. Suppose that there exist distinct $a,b,c \in V$ such that $\cl_T(ab) \subset \cl_T(abc) \subset V$. If $6 \in abc$, it follows that $\cl_T(abc) = V$. If $6 \notin abc$, we get $\cl_T(ab) = 12345 = \cl_T(abc)$, a contradiction in either case. Thus $\M_9 \in \ngu(2)$ by Lemma \ref{cltt}.

To show that $\M_9 \in \mngu(2)$, we must show that condition (ii) in Lemma \ref{cltt} holds when we remove one element from $\defe\M_9$:
\bi
\item
If we remove $145$, we get $\cl_T(15) = 15 \subset 156 = \cl_T(156)$.
\item
If we remove $123$, we get $\cl_T(12) = 12 \subset 126 = \cl_T(126)$.
\item
If we remove $134$, we get $\cl_T(34) = 34 \subset 346 = \cl_T(346)$.
\ei
The remaining case follows by symmetry, thus $\M_9 \in \mngu(2)$.

\smallskip

\noindent
$\M_{10}$: Notice that we can picture $\defe\M_{10}$ taking four nonadjacent faces of an octahedron, which ensures a remarkable degree of symmetry in $\defe\M_{10}$. Suppose that there exist distinct $a,b,c \in V$ such that $\cl_T(ab) \subset \cl_T(abc) \subset V$. If $abc \in \defe\M_{10}$, we get
$\cl_T(ab) = abc = \cl_T(abc)$, hence we may assume that $abc \notin \defe\M_{10}$. It is easy to check that $|abc \cap X| = 2$ for some $X \in \defe\M_{10}$, hence $X \subset \cl_T(abc)$, which must contain also a fourth element. This easily yields $\cl_T(abc) = V$, thus $\M_{10} \in \ngu(2)$ by Lemma \ref{cltt}.

To show that $\M_{10} \in \mngu(2)$, we must show that condition (ii) in Lemma \ref{cltt} holds when we remove one element from $\defe\M_{10}$. If we remove $146$, we get $\cl_T(14) = 14 \subset 146 = \cl_T(146)$. The remaining cases follow by symmetry, therefore $\M_{10} \in \mngu(2)$.

\smallskip

We show next that the complexes $\M_i$ $(i = 1,\ldots,10)$ are all mutually nonisomorphic.

Clearly, $\M_4$ cannot be isomorphic to any other complex because has larger defect. The same cardinality argument separates $\M_3$, $\M_8$ and $\M_9$ from the remaining cases. Now $\cup\, \defe\M_8 = V \supset \cup\, \defe\M_3 = \cup\, \defe\M_9$, hence $\M_3 \not\cong \M_8 \not\cong \M_9$. 
Now $\M_3 \not\cong \M_9$ because $4$ occurs in four elements of $\defe\M_3$, but no element occurs in four elements of $\defe\M_9$.

We also separate $\M_{10}$ from the others because $\defe\M_{10}$ is a PEG and this happens in no other case.

We still have to separate the cases $\M_i$ for $i \in 12567$. The property
$$\mbox{there exists $W \in P_4(V)$ such that $|P_3(W) \cap \defe\M_i| = 3$}$$
separates $\M_1$ and $\M_2$ from the others. Since $\cap \, \defe\M_1 = \emptyset \neq \cap \, \defe\M_2$, then $\M_1 \not\cong \M_2$. Now the number of vertices which appear in just a single element of $\defe\M_i$ is 2 when $i = 5$ and 1 when $i \in 67$. Thus we only have to separate $\M_6$ from $\M_7$. Since the property
$$\mbox{there exists some $a \in V$ such that, for every $b \in V \setminus \{ a \}$, $ab \subset X$ for some $X \in \defe\M_i$}$$
holds for $i = 7$ but not for $i = 6$, we get $\M_6 \not\cong \M_7$. Therefore the complexes $\M_i$ $(i = 1,\ldots,10)$ are all mutually nonisomorphic.

\smallskip

Finally, we prove that every $\H\, = (V,H) \in \mngu(2)$ is isomorphic to $\M_i$ for some $i \in \{ 1, \ldots, 10 \}$.  

We define
$$\omega(\H) = \min\{ |H \cap P_3(X)| : X \in P_4(V) \}.$$
Since $|P_3(X)| = 4$ for every $X \in P_4(V)$, we have $0 \leq \omega(\H) \leq 4$. However, we claim that values 0 and 4 can be excluded.
Indeed $\omega(\H) > 0$ by Lemma \ref{inco}, and $\omega(\H) = 4$ yields $\H = U_{3,6} \in \gu(2)$, thus we may split our discussion into four mutually exclusive cases.

\smallskip

\noindent
\underline{Case A}: $\omega(\H) = 1$.

\smallskip

Without loss of generality, we may assume that $H \cap P_3(1234) = \{ 123 \}$. This implies $\{ 124, 134, 234 \} \subseteq \defe\H$.

Suppose first that $a56 \in \defe\H$ for some $a \in 1234$. If $a = 4$, then $\defe\H \supseteq \defe\M_2$ , hence $H \subseteq M_2$ and so $\H = \M_2$ since  $\H\, \in \mngu(2)$. Thus, out of symmetry, we may assume that $a = 3$, yielding $\defe\H \supseteq \defe\M_1$ and consequently $\H = \M_1$.

Thus we may assume that $a56 \in H$ for every $i \in 1234$. 

Suppose now that $\cup\, \defe\H \subset V$. Note that we cannot have $\cup\, \defe\H = 1234$ since otherwise $\defe\H = \{ 124, 134, 234 \}$ and we have already established that $\M_1 \in \mngu(2)$. Hence we may assume that $\cup\, \defe\H = 12345$. It follows that $\defe\H = \{ 124, 134, 234 \} \cup Q$ for some nonempty $Q \subseteq \{ 125, 135, 235, 145, 245, 345 \}$. 

Write
$$\defe_2\H = \displaystyle\cup_{X \in \defe\H} P_2(X).$$
Suppose that $i5 \notin \defe_2\H$ for some $i \in 1234$. Then $\cl_T(i5) = i5 \subset i56 = \cl_T(i56)$ and so $\H \in \gu(2)$ by Lemma \ref{cltt}. Hence $15,25,35,45 \in \defe_2\H$. 

If $Q \supseteq \{ 145, 245, 345 \}$, the usual maximality argument yields $\H = \M_4$. Thus, out of symmetry, we may assume that $345 \notin Q$. Suppose that $145, 245 \in Q$. Since $35 \in \defe_2\H$, we get $135 \in Q$ or $235 \in Q$. If $135 \in Q$, then $\H$ is a subcomplex of $\M_3 \in \mngu(2)$, hence $\H = \M_3$. If $235 \in Q$, then $\H$ is a subcomplex of $\M_3(2,1,3,4,5,6) \in \mngu(2)$, hence $\H \cong \M_3$. Therefore we may assume that $245 \notin Q$. Since $45 \in \defe_2\H$, we get $145 \in Q$. 

If $235 \in Q$, then $\H$ is a subcomplex of $\M_3(2,1,3,4,5,6) \in \mngu(2)$, hence $\H \cong \M_3$. Hence we may assume that $235 \notin Q$. Since 
$25,35 \in \defe_2\H$, we get $Q =  \{ 125, 135, 145 \}$. Let $\H' = (V, H \cup \{ 134\})$. Then $\defe\H' = \{ 124, 234, 125, 135, 145 \}$, hence $\H' = \M_3(2,5,4,1,3,6) \in \mngu(2)$, contradicting $\H \in \mngu(2)$. 

Therefore $\H \cong \M_3$ or $\M_4$ if $\cup\, \defe\H \subset V$. Thus we assume now that $\cup\, \defe\H = V$. 

Suppose that $i5, i6 \notin \defe_2\H$ for some $i \in 1234$. Then $\cl_T(i5) = i5 \subset i56 = \cl_T(i56)$ and so $\H \in \gu(2)$ by Lemma \ref{cltt}. Hence, for each $i \in 1234$, at least one among $i5,i6$ belongs to $\defe_2\H$. 

Suppose next that 
\beq
\label{6mngu1}
\{ i45,i46 \} \cap \defe\H \neq \emptyset \mbox{ for every }i \in 123.
\eeq
Then we may assume that both 5 and 6 occur in these intersections, otherwise $\H$ would be a proper subcomplex of either $\M_4$ or $\M_4(1,2,3,4,6,5)$. Without loss of generality, we may assume that $145,345,246 \in \defe\H$. Let $\H' = (V,H \cup \{ 124,234\})$. Then $\defe\H' \supseteq \{ 134,145,345,246 \}$, hence $\H'$ is a subcomplex of $\M_2(1,5,3,4,2,6)$. This contradicts $\H \in \mngu(2)$, hence (\ref{6mngu1}) fails.

Thus, out of symmetry, we may assume that $345, 346 \in H$. Suppose that 
\beq
\label{6mngu2}
\{ i45,i46 \} \cap \defe\H \neq \emptyset \mbox{ for every }i \in 12.
\eeq
We know that at least one among $35,36$ belongs to $\defe_2\H$. Without loss of generality, we may assume that $35 \in \defe_2\H$. Since $345 \in H$, we may assume out of symmetry that $135 \in \defe\H$. But by (\ref{6mngu2}), at least one among $245,246$ belongs to $\defe\H$. If $245 \in \defe\H$, then $\H$ is a subcomplex of $\M_3$, a contradiction, hence we may assume that $246 \in \defe\H$. Let $\H' = (V,H\cup \{ 234\})$. Then $\defe\H' \supseteq \{ 124, 134, 135,246 \}$, hence $\H'$ is a subcomplex of $\M_5(1,3,4,2,5,6) \in \mngu(2)$, a contradiction. Therefore (\ref{6mngu2}) also fails, and we may assume that $245, 246 \in H$. 

Since at least one among $45,46$ belongs to $\defe_2\H$, we may assume out of symmetry that $145 \in \defe\H$. If $235$ or $236 \in \defe\H$, we proceed as in the discussion of (\ref{6mngu2}), hence we may assume that $235, 236 \in H$. Therefore we have
$$245, 246, 345, 346, 235, 236 \in H, \quad 124, 134, 234, 145 \in \defe\H.$$
On the other side, for each $i \in 23$, at least one among $i5,i6$ belongs to $\defe_2\H$. Thus
$$\{ 125, 126 \} \cap \defe\H \neq \emptyset, \quad \{ 135, 136 \} \cap \defe\H \neq \emptyset.$$
If $125, 135 \in \defe\H$, we get a subcomplex of one of the complexes studied in the subcase $\cup\, \defe\H \subset V$, hence we may assume out of symmetry that $126 \in \defe\H$. We are now left with a final alternative:

Suppose that $135 \in \defe\H$. Let $\H' = (V,H\cup \{ 124, 234\})$. Then $\defe\H' \supseteq \{ 134, 126, 135, 145 \}$, hence $\H'$ is a subcomplex of $\M_2(4,5,3,1,2,6) \in \mngu(2)$, a contradiction.

Suppose now that $136 \in \defe\H$. Let $\H' = (V,H\cup \{ 124, 134\})$. Then $\defe\H' \supseteq \{ 234, 126, 136, 145 \}$, hence $\H'$ is a subcomplex of $\M_7(6,2,1,3,4,5) \in \mngu(2)$, also a contradiction.

\smallskip

\noindent
\underline{Case B}: $\omega(\H) = 2$ and there exist distinct $a_1,\ldots,a_5 \in V$ such that $a_1a_2a_3, a_1a_3a_4, a_2a_3a_5 \in \defe\H$.

\smallskip

We may assume that $a_i = i$ for $i = 1,\ldots,5$. 

\smallskip

\noindent
\underline{Subcase B1}: $\defe\H \cap P_3(12345) \supset \{ 123, 134, 235 \}$.

\smallskip

Suppose first that $345 \in \defe\H$. Note that $\defe\H \supset \{ 123, 134, 235, 345 \}$, otherwise $\cl_T(15) = 15 \subset 12345 = \cl_T(125)$. Let $X \in \defe\H \setminus \{ 123, 134, 235, 345 \}$. If $6 \notin X$, we fall into Case A, hence we may assume that $6 \in X$. Out of symmetry, we may assume that $X \in \{ 256, 356, 246 \}$.

Suppose that $X = 256$. Let $\H' = (V,H\cup \{ 123\})$. Then $\defe\H' \supseteq \{ 134, 235, 345, 256 \}$, hence $\H'$ is a subcomplex of $\M_5(3,4,5,2,1,6) \in \mngu(2)$, a contradiction.

Suppose next that $X = 356$ and $126, 146, 256, 456 \notin \defe\H$ (out of symmetry). Note that $\defe\H \supset \{ 123, 134, 235, 345, 356 \}$, otherwise $\cl_T(24) = 24 \subset 246 = \cl_T(246)$. To avoid this inclusion, $\defe\H$ must contain some $Y$ such that  $6 \in Y$ and $|Y \cap 24| \neq \emptyset$. Out of symmetry, we may assume that $Y \in \{ 236, 246 \}$. Now $Y = 236$ makes us fall into Case A, and $Y = 246$ yields a subcomplex of $\M_7(1,2,3,4,6,5) \in \mngu(2)$.

Finally, suppose that $X = 246$ and $126, 146, 256, 456, 136, 236, 346, 356 \notin \defe\H$ (out of symmetry). Note that $\defe\H \supset \{ 123, 134, 246, 345 \}$, otherwise $\cl_T(15) = 15 \subset 156 = \cl_T(156)$. With all the exclusions stated, the only possibility left is now $\defe\H = \{ 123, 134, 156, 246, 345 \}$. But now $\H$ is a subcomplex of $M_7(3,2,1,4,6,5)$.

Hence we may assume that $345 \notin \defe\H$. To avoid falling back into Case A, and  since we are in Subcase B1, we may assume that $\{ 145,245 \} \cap \defe\H \neq \emptyset$. If $145, 245 \in \defe\H$, then $\defe\H \supseteq \{ 123, 134, 235, 145, 245 \}$ and so $\H = \M_9$ since $\H \in \mngu(2)$. 

On the other hand, if we define $\J = (V,J)$ and $\J' = (V,J')$ by $\defe\J\,  = \{ 123, 134, 235, 245 \}$ and $\defe\J'  = \{ 123, 134, 235, 145 \}$, then $\J' = \J(4,3,1,5,2,6)$, hence we may assume that $\defe\H\, = \{ 123, 134, 235, 245 \} \cup Q$, where $6 \in X$ for every $X \in Q$. 
\bi
\item
If $126 \in Q$, then $\H$ is a proper subcomplex of $\M_7(5,3,2,4,1,6)$, a contradiction.
\item
If $136 \in Q$, then $\H$ is a proper subcomplex of $\M_5(2,5,3,1,4,6)$, a contradiction.
\item
If $146 \in Q$, then $\H$ is a proper subcomplex of $\M_5(1,4,3,2,6,5)$, a contradiction.
\item
If $156 \in Q$, then $\H$ is a proper subcomplex of $\M_7(3,2,1,4,5,6)$, a contradiction.
\item
If $256 \in Q$, then $\H$ is a proper subcomplex of $\M_5(2,5,3,1,6,4)$, a contradiction.
\item
If $356 \in Q$, then $\H$ is a proper subcomplex of $\M_7$, a contradiction.
\item
If $456 \in Q$, then $\H$ is a proper subcomplex of $\M_6(1,4,3,2,6,5)$, a contradiction.
\ei
Therefore we must have $Q \supseteq \{ 236, 246, 346 \}$. But then we have in any case $\cl_T(15) = 15 \subset 156 = \cl_T(156)$, contradicting $\H \in \mngu(2)$.  

\smallskip

\noindent
\underline{Subcase B2}: $\defe\H \cap P_3(12345) = \{ 123, 134, 235 \}$.

\smallskip

Since $\{ 123, 134, 235 \} \subset \defe\M_9$, then $\defe\H$ must contain extra elements, and 6 must belong to each one of them.
\bi
\item
If $456 \in \defe\H$, then $\H$ is a subcomplex of $\M_6(1,4,3,2,6,5)$ and so $\H \cong \M_6$. 
\item
If $146 \in \defe\H$, then $\H$ is a subcomplex of $\M_5(1,4,3,2,6,5)$ and so $\H \cong \M_5$. 
\item
If $256 \in \defe\H$, then $\H \cong \M_5$ follows out of symmetry. 
\ei
Thus we may assume that $456, 146, 256 \notin \defe\H$. 

Suppose now that $126 \in \defe\H$. If $|X \cap 456| \leq 1$ for every $X \in \defe\H$, we get $\cl_T(56) = 56 \subset 456 = \cl_T(456)$, a contradiction in view of Lemma \ref{cltt}. Hence there exists some $X \in \defe\H$ such that $|X \cap 456| \geq 2$. If $X = 456$, then $\H$ is a proper subcomplex of $\M_6(1,5,2,3,4,6)$, a contradiction. Thus $|X \cap 456| = 2$, and we can assume out of symmetry that $X \cap 456 = 56$ and $X = 156$. But then we fall into Subcase B1 by applying the permutation $(12)(465)$ to $V$. 

Therefore we may also assume that $126 \notin \defe\H$. Then $\defe\H\, = \{ 123, 134, 235 \} \cup Q$ for some nonempty 
$$Q \subseteq \{ 136, 156, 236, 246, 346, 356 \}.$$
If 2 and 4 do not occur in $Q$, we get $\cl_T(24) = 24 \subset 246 = \cl_T(246)$, a contradiction. Similarly, 1 or 5 must occur in $Q$, and 4 or 6 must occur in $Q$ as well.  Hence
$$\mbox{(4 and 1 occur in $Q$) or (4 and 5 occur in $Q$) or (2 and 5 occur in $Q$).}$$
Out of symmetry, we may assume that one of the first two alternatives holds. Also out of symmetry, we can restrict the discussion to five cases:
\bi
\item
If $246, 136 \in \defe\H$, then $\H$ is a proper subcomplex of $\M_7(1,4,3,6,2,5)$, a contradiction.
\item
If $246, 156 \in \defe\H$, then $\H$ is a proper subcomplex of $\M_{10}(3,2,5,4,6,1)$, a contradiction.
\item
If $346, 136 \in \defe\H$, we fall into Case A since $H \cap P_3(1346) = \{ 146 \}$.
\item
If $346, 156 \in \defe\H$, then $\H$ is a proper subcomplex of $\M_7(2,5,3,1,6,4)$, a contradiction.
\item
If $246, 136 \in \defe\H$, then $\H = \M_8$.
\ei

\smallskip

\noindent
\underline{Case C}: $\omega(\H) = 2$ and there exist no distinct $a_1,\ldots,a_5 \in V$ such that $a_1a_2a_3, a_1a_3a_4, a_2a_3a_5 \in \defe\H$.

\smallskip

The assumptions made so far imply that $\defe\H\, = \{ 123, 134, 235 \} \cup Q$ with
$$Q \subseteq \{ 135, 136, 245, 246, 156, 256, 356, 456 \}.$$

Suppose first that $245, 246 \notin Q$. 
\bi
\item
If $256, 456 \notin Q$, then $\cl_T(24) = 24 \subset 245 = \cl_T(245)$, contradicting Lemma \ref{cltt}.
\item
If $135, 136 \notin Q$, then $\cl_T(24) = 24 \subset 1234 = \cl_T(124)$, contradicting Lemma \ref{cltt}.
\item
If $156, 356 \notin Q$, then $\cl_T(24) = 24 \subset \cl_T(245) \subseteq 2456$, contradicting Lemma \ref{cltt}.
\ei
Thus we may assume out of symmetry that $135, 156, 256 \in Q$. But then we fall into Case B.

Suppose next that $245, 246 \in Q$. To avoid Case A, it follows that $256, 456 \notin Q$.
Suppose that $156, 356 \notin Q$. Then $\cl_T(56) = 56 \subset 256 = \cl_T(2456)$, contradicting Lemma \ref{cltt}. Out of symmetry, we may assume that $156 \in Q$ and so $\defe\H \supseteq \{ 123, 134, 245, 246, 156 \}$. But then $\H$ is a subcomplex of $\M_7(2,5,4,6,1,3)$, thus $\H \cong \M_7$.

Therefore we may assume that $|Q \cap \{ 245,246\}| = 1$. Out of symmetry, we may assume that $245 \in Q$. Suppose that $156, 256, 356 \notin Q$. Then 
$\cl_T(26) = 26 \subset 2456 = \cl_T(246)$, contradicting Lemma \ref{cltt}. It follows that $Q \cap \{ 156, 256, 356 \} \neq \emptyset$.
\bi
\item
If $156 \in Q$, then $\H = \M_7(3,2,1,4,5,6)$.
\item
If $356 \in Q$, then $\H = \M_7$.
\ei
Thus we may assume that $Q \cap \{ 156, 256, 356 \} = \{ 256 \}$. Suppose that $456 \in Q$. Then we fall into Case A since $H \cap P_3(2456) = \{ 246 \}$. Hence $456 \notin Q$. But then $\cl_T(46) = 46 \subset 2456 = \cl_T(246)$, contradicting Lemma \ref{cltt}. 

\smallskip

\noindent
\underline{Case D}: $\omega(\H) = 3$.

\smallskip

Without loss of generality, we may assume that $123, 245 \in \defe\H$. Since $\{ 123, 245 \} \subset \defe\M_7$, there exists some $X \in \defe\H \setminus \{ 123, 245 \}$. Since $|X \cap 123|, |X \cap 245| \leq 1$, we may assume without loss of generality that $X = 356$. Since $\{ 123, 245, 356 \} \subset \defe\M_7$, there exists some $Y \in \defe\H \setminus \{ 123, 245, 356 \}$. Since $|Y \cap Z| \leq 1$ for every $Z \in \{ 123, 245, 356 \}$, it follows that $2,3,5 \notin Y$. Thus $Z = 146$ and so $\H = \M_{10}$.  
\qed

\subsection{Computing $\mgu(2)$}

Let $(V,H)$ and $(V',H')$ be two simplicial complexes.
Given a bijection $\p:V \to V'$, we write
$H\p =  \{ X\p \mid X \in H\}$. We will adopt the notation  $(V,H) \leq (V',H')$ to express that there exists a bijection $\p:V \to V'$ such that $H\p \subseteq H'$. 

In this section, we adopt the notation $V_n = 1\ldots n$. Given natural numbers $2 \leq i < j < n$, we define a simplicial complex 
$\J(i,j,n) = (V_n, J(i,j,n)$, where
$$J(i,j,n) = P_{\leq 2}(V_n) \cup \{ aa'b \mid 1 \leq a <  a' \leq i < b \leq j\} \cup \{ bb'c \mid 1 \leq b < b' \leq j < c \leq n\}.$$

\bl
\label{necess}
Let $\H = (V, H) \in {\rm GU}(2)$. 
Then $\H \geq \J(i,j,n)$ for some natural numbers $2 \leq i < j < n$.
\el

\proof
Since $\H$ goes up, it must have at least 4 vertices.
Since $\H \in {\rm mGU}(2)$, it follows from Lemma \ref{cltt} that there exists a chain $A \subset B \subset V$ in $T(H)$ with $|A| \geq 2$. Replacing $\H$ by an isomorphic image if needed, we may assume that $V = V_n$,  $A = 1\ldots i$ and $B = 1\ldots j$. Hence $2 \leq i < j < n$. 

Suppose that $1 \leq a <  a' \leq i < b \leq j$. Since $aa' \in H$, $aa' \subseteq A \in T(H)$ and $b \notin A$, we get $aa'b \in H$. Suppose now that $1 \leq b < b' \leq j < c \leq n$. Since $bb' \in H$, $bb' \subseteq B \in T(H)$ and $c \notin B$, we get $bb'c \in H$. 

Thus $J(i,j,n) \subseteq H$ and so $\H \geq \J(i,j,n)$.
\qed

\bl
\label{penul}
\bi
\item[(i)] If $2 \leq i < n-1$, then $\J(2,i+1,n) \leq \J(i,n-1,n)$.
\item[(ii)] If also $n > 4$, then $\J(2,i+1,n) \not\geq \J(i,n-1,n)$.
\ei
\el

\proof
(i) Let $\p: V_n \to V_n$ be the cyclic permutation $(n\ldots 3\, 2\, 1)$. We show that $(J(2,i+1,n))\p \subseteq J(i,n-1,n)$.

Let $X \in J(2,i+1,n) \cap P_3(V_n)$. Suppose first that $X = aa'b$ with $1 \leq a <  a' \leq 2 < b \leq i+1$. Then $X\p = 1(b-1)n \in J(i,n-1,n)$. 

Suppose now that $X = bb'c$ with $1 \leq b < b' \leq i+1 < c \leq n$. Consider first the case $b = 1$. Then $X\p = (b'-1)(c-1)n \in J(i,n-1,n)$. On the other hand, if $b > 1$, then $X\p = (b-1)(b'-1)(c-1) \in J(i,n-1,n)$ as well.

Therefore $(J(2,i+1,n))\p \subseteq J(i,n-1,n)$ and so $\J(2,i+1,n) \leq \J(i,n-1,n)$.

(ii) A cardinality argument will suffice. Let $(n-3)(n-2)n \in J(i,n-1,n)$. Then $((n-3)(n-2)n)\p\inv = 1(n-2)(n-1) \notin J(2,i+1,n)$. Therefore $(J(2,i+1,n))\p \subset J(i,n-1,n)$ and so $\J(2,i+1,n) \not\geq \J(i,n-1,n)$.
\qed

\bl
\label{cons}
If $3 \leq i < n-1$, then $\J(2,i+1,n) \leq \J(i,i+1,n)$ and $\J(2,i+1,n) \not\geq \J(i,i+1,n)$.
\el

\proof
Let $\p: V_n \to V_n$ be the cyclic permutation $((i+1)\, i \ldots 3\, 2\, 1)$. We show that $(J(2,i+1,n))\p \subseteq J(i,i+1,n)$.

Let $X \in J(2,i+1,n) \cap P_3(V_n)$. Suppose first that $X = aa'b$ with $1 \leq a <  a' \leq 2 < b \leq i+1$. Then $X\p = 1(b-1)(i+1) \in J(i,i+1,n)$. 

Suppose now that $X = bb'c$ with $1 \leq b < b' \leq i+1 < c \leq n$. Consider first the case $b = 1$. Then $X\p = (b'-1)(i+1)c \in J(i,i+1,n)$. On the other hand, if $b > 1$, then $X\p = (b-1)(b'-1)c \in J(i,i+1,n)$ as well.

Therefore $(J(2,i+1,n))\p \subseteq J(i,i+1,n)$ and so $\J(2,i+1,n) \leq \J(i,i+1,n)$.

Now let $(i-1)i(i+1) \in J(i,i+1,n)$. Then $((i-1)i(i+1))\p\inv = 1i(i+1) \notin J(2,i+1,n)$. Therefore $(J(2,i+1,n))\p \subset J(i,i+1,n)$ and so $\J(2,i+1,n) \not\geq \J(i,i+1,n)$.
\qed

For each $n \geq 4$, we define
$$Q_n = \{ (i,j) \in \N^2 \mid 2 \leq i < j-1 < n-2 \} \cup \{ (2,3) \}.$$

\bl
\label{tj}
Let $(i,j) \in Q_n$. 
\bi
\item[(i)] If $(i,j) \neq (2,3)$, then $T(J(i,j,n)) = P_{\leq 1}(V_n) \cup \{ 12\ldots i, 12\ldots j, V_n\}$.
\item[(ii)] If $(i,j) = (2,3)$ and $n > 4$, then $T(J(i,j,n)) = P_{\leq 1}(V_n) \cup \{ 12, 13, 23, 123, V_n\}$.
\item[(iii)] If $(i,j) = (2,3)$ and $n = 4$, then $T(J(i,j,n)) = 2^{V_n}$.
\ei
\el

\proof
(i) It follows from the definitions that the subsets in the right hand side are indeed in $T(J(i,j,n))$. 

If $X \in P_{\leq 1}(V_n)$, then $\cl_T(X) = X$. We consider now the case $X \in P_2(V_n)$.
We may assume that $X = ab$ with $a < b$.

Suppose first that $b \leq i$. Then $a < i$. For every $c \in X \setminus ab$, we have $abc \notin H$, hence $12\ldots i \subseteq \cl_T(X)$. Since $12\ldots i \in T(J(i,j,n))$, we get $\cl_T(X) = 12\ldots i$. 

Suppose next that $i < b \leq j$. If $i < a$, then $abc \notin H$ for every $c \in 12\ldots j \setminus ab$, hence $12\ldots j \subseteq \cl_T(X)$. Since $12\ldots j \in T(J(i,j,n))$, we get $\cl_T(X) = 12\ldots j$. If $a \leq i$, take some $c \neq b$ such that $i < c \leq j$ (there exists since $i < j-1$). Then $abc \notin H$, hence $c \in \cl_T(X)$ and so the previous case yields $12\ldots j = \cl_T(ac) \subseteq \cl_T(X) \subseteq 12\ldots j$. Therefore $\cl_T(X) = 12\ldots j$.

Thus we may assume that $j < b$. If $j < a$, then $abc \notin H$ for every $c \in V_n \setminus ab$, hence $\cl_T(X) = V_n$. Hence we may assume that
$a \leq j$. Take some $c > j$ distinct from $b$ (there exists since $j-1 < n-2$). Then $abc \notin H$, hence $c \in \cl_T(X)$ and by the preceding case we get
$V_n = \cl_T(bc) \subseteq \cl_T(X)$. Therefore $\cl_T(X) = V_n$ and so $\cl_T(X) \in \{ 12\ldots i, 12\ldots j, V_n\}$ for every $X \in P_2(V_n)$. 

So let $Y \in T(J(i,j,n))$. We may assume that $|Y| \geq 2$. If $Y \subseteq 12\ldots i$, it follows from the discussion of 2-set closures that $Y = \cl_T(Y) = 12\ldots i$. If $Y \subseteq 12\ldots j$ but $Y \not\subseteq 12\ldots i$, it follows from the discussion of 2-set closures that $Y = \cl_T(Y) = 12\ldots j$. If $Y \not\subseteq 12\ldots j$, it follows from the discussion of 2-set closures that $Y = \cl_T(Y) = V_n$.

(ii) It follows from the definitions that the subsets in the right hand side are indeed in $T(J(i,j,n))$. Similarly to the proof of (i), we compute $\cl_T(X)$ for every $X \in P_2(V)$, getting $\cl_T(X) \in \{ 12, 13, 23, V_n\}$. The claim now follows easily.

(ii) Since $\J(2,3,4) = U_{3,4}$.
\qed

\bl
\label{ordert}
Let $d \geq 2$ and $(V,H),(V,H') \in {\rm Pav}(d)$. Then $H \subseteq H'$ implies $T(H) \subseteq T(H')$.
\el

\proof
Since $P_d(V) \subseteq H$, we have
$$T(H) = \{ T \subseteq V \mid X \cup \{ p \} \in H \mbox{ for all $X \in P_d(T)$ and $p \in V \setminus T$}\}.$$
The analogous equality holds for $T(H')$, hence $H \subseteq H'$ implies $T(H) \subseteq T(H')$.
\qed

\bt
\label{computemgu}
Let $n \geq 4$. For $n$ vertices, there exist precisely $\frac{n^2-9n+22}{2}$ isomorphism classes in ${\rm mGU}(2)$, given by representatives
$\J(i,j,n)$, for $(i,j) \in Q_n$.
\et

\proof
Assume first that $n = 4$. By the proof of Proposition \ref{4mngu}, $(V_4,H) \in \ngu(2)$ if $H \subset P_{\leq 3}(V_4)$. On the other hand, it follows from Lemma \ref{tj}(iii) that $\J(2,3,4) = (V_4, P_{\leq 3}(V_4)) \in \gu(2)$. Thus, up to isomorphism, $\J(2,3,4)$ is the unique complex in $\mgu(2)$ with 4 vertices. Since $\frac{4^2-9\cdot 4 +22}{2} = 1$, the theorem holds for $n = 4$.

Assume from now on that $n \geq 5$. Suppose that $\J(i,j,n) \leq \J(i',j',n)$ for some distinct $(i,j),(i',j') \in Q_n$. Then there exists some permutation $\p$ of $V_n$ such that $(J(i,j,n))\p \subseteq J(i',j',n)$. By Lemma \ref{ordert}, we get $T((J(i,j,n))\p) \subseteq T(J(i',j',n))$. By Lemma \ref{tj}, the cardinalities of the elements of $T(J(i,j,n))$ and therefore of $T((J(i,j,n))\p)$) can be ordered by $0 < 1 < i < j < n$. On the other hands, the cardinalities of the elements of $T(J(i',j',n))$ can be ordered by $0 < 1 < i' < j' < n$. Since $(i,j) \neq (i',j')$, the two sets of cardinalities are incomparable, contradicting $T((J(i,j,n))\p) \subseteq T(J(i',j',n))$. Therefore $\J(i,j,n) \leq \J(i',j',n)$ implies $(i,j,n) = (i',j',n)$.

We show next that $\J(i,j,n) \in \mgu(2)$ for every $(i,j) \in Q_n$. Let $\H\, = (V_n, H)$ with $H = J(i,j,n) \setminus \{ X\}$, where $X \in P_3(V_n) \cap J(i,j,n)$. Since $H \subset J(i,j,n)$, it follows from Lemma \ref{ordert} that $T(H) \subseteq T(J(i,j,n))$. Note that $P_{\leq 1}(V_n) \cup \{ V_n\} \subseteq T(H)$ in any case since $P_2(V) \subseteq H$.

Assume first that $(i,j) \neq (2,3)$. Suppose that $X = aa'b$ with $1 \leq a < a' \leq i < b \leq j$. Since $aa'b \notin H$, we have $12\ldots i \notin T(H)$. Hence $T(H) \subseteq T(J(i,j,n))$ and Lemma \ref{tj}(i) yield $T(H) \subseteq P_{\leq 1}(V_n) \cup \{ 12\ldots j, V_n\}$. Therefore $\H\, \in \ngu(2)$ by Lemma \ref{cltt}. Suppose now that $X = bb'c$ with $1 \leq b < b' \leq j < c \leq n$. Since $bb'c \notin H$, we have $12\ldots j \notin T(H)$. Hence $T(H) \subseteq T(J(i,j,n))$ and Lemma \ref{tj}(i) yield $T(H) \subseteq P_{\leq 1}(V_n) \cup \{ 12\ldots i, V_n\}$. Therefore $\H\, \in \ngu(2)$ by Lemma \ref{cltt}.

Assume now that $(i,j) = (2,3)$. Suppose that $X = 123$. Since $123 \notin H$, we have $\{ 12, 13, 23 \} \cap T(H) = \emptyset$. Hence $T(H) \subseteq T(J(i,j,n))$ and Lemma \ref{tj}(ii) yield $T(H) \subseteq P_{\leq 1}(V_n) \cup \{ 123, V_n\}$. Therefore $\H\, \in \ngu(2)$ by Lemma \ref{cltt}. Suppose now that $X = bb'c$ with $1 \leq b < b' \leq 3 < c \leq n$. We assume that $b = 1$ and $b' = 2$, the other cases being similar. Since $12c \notin H$, we have $12, 123 \notin T(H)$. Hence $T(H) \subseteq T(J(i,j,n))$ and Lemma \ref{tj}(ii) yield $T(H) \subseteq P_{\leq 1}(V_n) \cup \{ 13,23, V_n\}$. Therefore $\H\, \in \ngu(2)$ by Lemma \ref{cltt}.

Therefore $\J(i,j,n) \in \mgu(2)$ for every $(i,j) \in Q_n$.
It follows from Lemmas \ref{necess}, \ref{penul} and \ref{cons} that every $\H \in \mgu(2)$ with $n$ vertices satisfies $\H \geq \J(i,j,n)$ for some $(i,j) \in Q_n$. Thus $\H \cong \J(i,j,n)$. Since we have already proved that the $\J(i,j,n)$ are nonisomorphic elements of $\ngu(2)$ with 6 vertices, it follows that they act as representatives of the isomorphism classes of ${\rm mGU}(2)$ with $n$ vertices.

Their number is 
$$\begin{array}{lll}
|Q_n|&=&(\sum_{i = 2}^{n-4} \sum_{j = i+2}^{n-2} 1) + 1 = (\sum_{i = 2}^{n-4} (n-i-3)) + 1 = (n-5)(n-3) - (\sum_{i = 1}^{n-4} 1) + 2\\ &&\\
&=&(n-5)(n-3) - \frac{(n-4)(n-3)}{2} + 2 = \frac{(n-3)(2n-10-n+4) +4}{2} = \frac{(n-3)(n-6) +4}{2} = \frac{n^2-9n+22}{2}
\end{array}$$
as claimed.
\qed

The next two corollaries help to classify the complexes in $\mgu(2)$:

\bc
\label{brj}
Let $n \geq 4$ and $(i,j) \in Q_n$. Then:
\bi
\item[(i)] $\J(i,j,n) \in {\rm TBPav}(2)$;
\item[(ii)] $\J(i,j,n) \in {\rm BPav}(2)$ if and only if $j = 3$.
\ei
\ec

\proof
(i) We have $J(i,j,n) = B_2(V_n, 12\ldots i) \cup B_2(V_n, 12\ldots j)$, hence $\J(i,j,n) \in {\rm TBPav}(2)$ by Proposition \ref{tpav}.

(ii) Assume first that $j = 3$. Then $i = 2$ and it follows from Lemma \ref{tj} that $P_2(123) \subseteq T(J(2,3,n))$. It follows easily that 
$$P_{\leq 1}(V_n) \cup P_2(123) \cup \{ V_n \} \subseteq \flatx\J(2,3,n).$$
Since every $X \in J(2,3,n) \cap P_3(V_n)$ has at least two elements in 123 (say $a$ and $b$), then $X$ is a transversal of the successive differences for the chain
$$\emptyset \subset a \subset ab \subset V_n$$
in $\flatx\J(2,3,n)$. Therefore $\J(2,3,n)$ is boolean representable.

Assume now that $j \neq 3$. Then $j > 3$ and so 
$$T(J(i,j,n)) = P_{\leq 1}(V_n) \cup \{ 12\ldots i, 12\ldots j, V_n\}$$
by Lemma \ref{tj}. Now $12j \in J(i,j,n)$ but $12jn \notin J(i,j,n)$, hence $12\ldots j \notin \flatx\J(i,j,n)$. In view of Lemma \ref{propt}(ii), we get
$$\flatx\J(i,j,n) \subseteq P_{\leq 1}(V_n) \cup \{ 12\ldots i, V_n\}.$$
It follows that $1jn \in J(i,j,n)$ is not a a transversal of the successive differences for any chain
in $\flatx\J(i,j,n)$. Therefore $\J(i,j,n)$ is not boolean representable.
\qed

\bc
\label{nomat}
Given $\H \, \in {\rm mGU}(2)$, the following conditions are equivalent:
\bi
\item[(i)] $\H$ is a matroid;
\item[(ii)] $\H$ is pure;
\item[(iii)] $\H$ has 4 vertices.
\ei
\ec

\proof
(i) $\Rw$ (ii). It follows from the exchange property that every matroid is pure.

(ii) $\Rw$ (iii). In view of Theorem \ref{computemgu}, we may assume that $\H\, = \J(i,j,n)$ with $n > 4$ and $(i,j) \in Q_n$. Since $j-1 < n-2$, it follows that $(n-1)n$ is a facet of $\H$. Therefore $\H$ is not pure.

(iii) $\Rw$ (i). In this case, it follows from Theorem \ref{computemgu} and Lemma \ref{tj}(iii) that $\H$ is the uniform matroid $U_{3,4}$.
\qed

\bt
\label{bigc}
Let $\H \, = (V,H) \in {\rm mGU}(2)$ have more than four vertices. Then there exists some $p \in V$ such that the restriction of $\H$ to $V \setminus \{ p \}$ is in ${\rm mGU}(2)$.
\et

\proof
By Theorem \ref{computemgu}, we may assume that $\H\, = \J(i,j,n)$ for some $(i,j) \in Q_n$. Given $p \in V_n$, let $\H'_p = (V_n \setminus \{ p \}, J(i,j,n) \cap 2^{V_n \setminus \{ p \}})$ denote the restriction of $\H$ to $V \setminus \{ p \}$.

Suppose that $i > 2$. Take $p = 1$. Then
$J(i,j,n) \cap 2^{V_n \setminus \{ p \}} = B_2(2\ldots i) \cup B_2(2\ldots j)$ and $\H'_p \cong \J(i-1,j-1,n-1) \in \mgu(2)$ by Theorem \ref{computemgu}. Thus we may assume that $i = 2$. 

Suppose that $j > 3$. Take $p = j$. Then
$J(i,j,n) \cap 2^{V_n \setminus \{ p \}} = B_2(12) \cup B_2(12\ldots (j-1))$ and $\H'_p \cong \J(2,j-1,n-1) \in \mgu(2)$ by Theorem \ref{computemgu}. Thus we may assume as the last final case that $j = 3$.

Take $p = n$. Then
$J(i,j,n) \cap 2^{V_n \setminus \{ p \}} = B_2(12) \cup B_2(123)$ and $\H'_p \cong \J(2,3,n-1) \in \mgu(2)$ by Theorem \ref{computemgu}. 
\qed

\bc
\label{everyres}
Let $Q'_n$ denote the isomorphism classes of complexes in ${\rm mGU}(2)$ with $n$ vertices where every restriction to $n-1$ vertices is still in ${\rm mGU}(2)$.
\bi
\item[(i)] These isomorphism classes are given by representatives
$\J(i,j,n)$, where $3 \leq i \leq j-3 \leq n-6$.
\item[(ii)] The smallest such representative is $\J(3,6,9)$.
\item[(iii)] $|Q'_n| = \frac{n^2-15n+56}{2}$ for every $n \geq 9$.
\ei
\ec

\proof
(i) In view of Theorem \ref{computemgu}, we can restrict our attentions to complexes of the form $\H\, = \J(i,j,n)$ with $(i,j) \in Q_n$. If $3 \leq i \leq j-3 \leq n-6$, a straightforward adaptation of the proof of the aforementioned theorem yields:
\bi
\item
if $p \leq i$, then $\H'_p \cong \J(i-1,j-1,n-1) \in \mgu(2)$;
\item
if $i < p \leq j$, then $\H'_p \cong \J(i,j-1,n-1) \in \mgu(2)$;
\item
if $j < p$, then $\H'_p \cong \J(i,j,n-1) \in \mgu(2)$.
\ei

Suppose now that $i = 2$. Then $\H'_1 = (V_n,B_2(2\ldots j)) \notin \mgu(2)$ in view of Theorem \ref{computemgu}.

Suppose next that $3 \leq i = j-2$. Then $\H'_j \cong \J(i,i+1,n-1) \notin \mgu(2)$ in view of Lemma \ref{cons}.

Thus we may assume that $3 \leq i \leq j-3 = n-5$. But then $\H'_n \cong \J(i,n-2,n-1) \notin \mgu(2)$ in view of Lemma \ref{penul}. 

(ii) By part (i).

(iii) By part (i), we have
$$\ba{lll}
|Q'_n|&=&\sum_{i = 3}^{n-6} \sum_{j = i+3}^{n-3} 1 = \sum_{i = 3}^{n-6} (n-i-5) = (n-8)(n-5) - \frac{(n-3)(n-8)}{2} = \frac{(n-8)(2n-10-n+3)}{2}\\ &&\\
&=&\frac{n^2-15n+56}{2}
\ea$$
as claimed.
\qed



\bigskip

{\sc Stuart Margolis, Department of Mathematics, Bar Ilan University,
  52900 Ramat Gan, Israel} 

{\em E-mail address:} margolis@math.biu.ac.il

\bigskip

{\sc John Rhodes, Department of Mathematics, University of California,
  Berkeley, California 94720, U.S.A.}

{\em E-mail addresses}: rhodes@math.berkeley.edu, BlvdBastille@aol.com

\bigskip

{\sc Pedro V. Silva, Centro de
Matem\'{a}tica, Faculdade de Ci\^{e}ncias, Universidade do
Porto, R. Campo Alegre 687, 4169-007 Porto, Portugal}

{\em E-mail address}: pvsilva@fc.up.pt

\end{document}